\documentclass[a4paper]{article}

\usepackage[english]{babel}
\usepackage[T1]{fontenc}
\usepackage{latexsym}
\usepackage{ae}
\usepackage{amsthm,amsfonts}
\usepackage{amsmath}
\usepackage{amssymb}
\usepackage{theorem}
\usepackage{fancyhdr}
\usepackage{bm}
\usepackage[mathscr]{eucal}
\usepackage[hmarginratio = 1:1]{geometry}
\usepackage{layout}
\usepackage{enumitem}
\usepackage[hidelinks]{hyperref}

\renewcommand{\thefootnote}{\alph{footnote}}

\begin{document}

\begin{center}

\textbf{\Large Exponentiation of Lie Algebras of Linear Operators on Locally Convex Spaces}\\[0.5cm]
\textbf{\large Rodrigo A. H. M. Cabral}\\[0.5cm]
Departamento de Matem\'atica Aplicada, Instituto de Matem\'atica e Estat\'istica, Universidade de S\~ao Paulo (IME-USP), BR-05508-090, S\~ao Paulo, SP, Brazil. \footnote{E-mails: rahmc@ime.usp.br; rodrigoahmc@gmail.com 


This work is part of a PhD thesis \cite{rodrigotese}, which was supported by CNPq (Conselho Nacional de Desenvolvimento Cient\'ifico e Tecnol\'ogico).

}\\[0.5cm]
\textbf{\large Abstract}

\end{center}

Necessary and sufficient conditions for the exponentiation of finite-dimensional real Lie algebras of linear operators on complete Hausdorff locally convex spaces are obtained, focused on the equicontinuous case - in particular, necessary conditions for exponentiation to compact Lie groups are established. Applications to complete locally convex algebras, with special attention to locally C$^*$-algebras, are given. The definition of a projective analytic vector is introduced, playing an important role in some of the exponentiation theorems.

\tableofcontents \vspace{0.5cm}

\renewcommand{\thefootnote}{\arabic{footnote}}
\setcounter{footnote}{0}

\section*{Introduction}
\addcontentsline{toc}{section}{Introduction}

\indent

Infinitesimal generators of semigroups or groups are mathematical objects which arise in various contexts. One of the most general ones is that of locally convex spaces, where they appear as particular kinds of linear operators, defined from an action of the additive semigroup $[0, +\infty)$ (in the case of semigroups) or the additive group $\mathbb{R}$ of real numbers (in the case of one-parameter groups) on the subjacent space. This action is usually implemented by continuous linear operators and subject to a condition of continuity: it must be continuous with respect to the usual topology of $[0, + \infty)$, or $\mathbb{R}$, and a fixed locally convex topology on the space of continuous linear operators (usually, at least in the infinite-dimensional framework, there are several options for this choice). One of the possible choices gives rise to strongly continuous semigroups and strongly continuous groups, and will be the topology of choice for the investigations in this work.

Within the context of normed spaces, there exist two very important classes of such strongly continuous actions: contraction semigroups and groups of isometries. In Hilbert space theory, for example, strongly continuous one-parameter groups of isometries are implemented by unitary operators and, according to the spectral theorem for self-adjoint operators and Stone's theorem (see Theorems VIII.7 and VIII.8, of \cite{reedsimon1}), their generators are precisely the anti-self-adjoint operators (in other words, operators such that $A^* = -A$). The Feller-Miyadera-Phillips Theorem characterizes the generators of strongly continuous one-parameter (semi)groups on Banach spaces - see \cite[Theorem 3.8, page 77]{engel} and \cite[Generation Theorem for Groups, page 79]{engel}. When specialized to the contractive and to the isometric cases, it yields the respective versions of the famous Hille-Yosida Theorem \cite[Theorem 3.5, page 73]{engel}. Also in this context there is the Lumer-Phillips Theorem (\cite[Theorem 3.15, page 83]{engel}) which, for groups of isometries, states that the infinitesimal generators must be conservative linear operators whose perturbations by non-zero multiples of the identity are surjective.
This is a direct generalization of Stone's Theorem to Banach spaces, since anti-symmetric (or anti-hermitian, if one prefers) operators are conservative and the self-adjointness property for such operators translates to an analogous surjectivity condition (see \cite[Theorem VIII.3]{reedsimon1}).

For more general complete locally convex spaces, there is a theorem which characterizes the generators of equicontinuous semigroups in an analogous way that the Lumer-Phillips Theorem does \cite[Theorem 3.14]{albanese}. Also in this more general setting, references \cite[Theorem 4.2]{babalola} and \cite[Corollary 4.5]{babalola} give characterizations in the same spirit that the Feller-Miyadera-Phillips and Hille-Yosida Theorems do, respectively.\\

When the locally convex space under consideration has an additional algebraic structure, turning it into an algebra (or a $*$-algebra), the actions of interest are by strongly continuous one-parameter groups of automorphisms (or $*$-automorphisms) and their generators then become derivations (or $*$-derivations), since they satisfy the Leibniz product rule. Such situations will be dealt with in the last part of this paper.\\

There is another very important direction for generalization, which consists in passing from one-parameter groups to more general groups. One framework that comes to mind here would be to consider abstract topological groups and their actions on locally convex spaces. But often, the study of smooth and analytic elements of a given action is very important and useful, so a natural requirement on the group is that it should be a Lie group. Also, it is very important that the group in question is not restricted to be commutative. Classical results in this direction, in the Hilbert and Banach space contexts, may be found in the influential paper \cite{nelson} of Edward Nelson. Two of the theorems there state that every strongly continuous Lie group representation has a dense subspace of analytic vectors when the representation is by unitary operators on a Hilbert space \cite[Theorem 3]{nelson} and, more generally, by bounded linear operators on a Banach space \cite[Theorem 4]{nelson}.\footnote{See also \cite[Theorem 5, page 56]{moore}.}\\

It should be mentioned that the present work is part of a PhD thesis, and that the study of reference \cite{nelson} was suggested by the author's advisor. This was the starting point for the research and, in one of the discussions, it was suggested that the author investigate what was the actual role that the generalized Laplacian defined in \cite{nelson} played in some of its results.\footnote{The author would like to thank his advisor, Michael Forger, for suggesting the study of Lie group representations in the context of locally convex $*$-algebras. The author would also like to thank Professor Severino T. Melo for the many discussions on the papers \cite{nelson} and \cite{harish-chandra} in the early stages of this work, and also for his teachings on the theory of pseudodifferential operators.}\\

Conversely, one may investigate the exponentiability of a (real finite-dimensional) Lie algebra $\mathcal{L}$ of linear operators: when can its elements be obtained (as generators of one-parameter groups) from a strongly continuous representation $V$ of a (connected) Lie group $G$, having a Lie algebra $\mathfrak{g}$ isomorphic to $\mathcal{L}$ via $\eta \colon \mathfrak{g} \longrightarrow \mathcal{L}$, according to the formula $$\frac{d}{dt} \, V(\exp tX) (x) \, \Big|_{t = 0} = \eta(X)(x), \qquad X \in \mathfrak{g}, \, x \in \mathcal{D},$$ where $\text{exp}\colon \mathfrak{g} \longrightarrow G$ is the corresponding exponential map? (See Definition 2.6 for the precise formulation of exponentiation)

For representations on finite-dimensional vector spaces, a classical theorem states that every Lie algebra of linear operators is exponentiable, provided only that one chooses $G$ to be simply connected\footnote{A topological space $\mathcal{X}$ is said to be simply connected if it is path-connected - so, in particular, it is connected - and if every continuous curve $\gamma \colon [0, 1] \longrightarrow \mathcal{X}$ satisfying $\gamma(0) = \gamma(1) =: x_0$ (this is called a loop based at $x_0$) is path homotopic to the constant curve $\gamma_0 \colon [0, 1] \ni t \longmapsto x_0$ (see \cite[page 333]{munkres}).} (see the brief discussion after Definition 2.6).

In the infinite-dimensional case, things are much more complicated, and when formulating the question some reasonable a priori requirements are usually assumed: the elements of $\mathcal{L}$ are linear operators defined on a common, dense domain $\mathcal{D}$ which they all map into itself. But there are counterexamples showing that these are not sufficient: one can find linear operators satisfying all these conditions whose closures generate strongly continuous one-parameter groups, but the closures of certain of their real linear combinations (including perturbations of one by the other), or of their commutator, do not - see the result stated at the end of Section 10 of \cite{nelson}; see also \cite[Example 4.1]{bratteliheat}. Finding sufficient criteria for exponentiability is thus an intricate problem because these must exclude such unpleasant situations. One of the first important results in this direction, for Lie algebras of anti-symmetric (anti-hermitian) operators on Hilbert spaces, can be found in \cite[Theorem 5]{nelson}.

As for representations by $*$-automorphisms of C$^*$-algebras, there exist several studies in the literature. For example, in \cite{bratteliu1} and \cite{bratteliu2}, generation theorems\footnote{In other words, theorems which give necessary and (or) sufficient conditions on a linear operator in order for it to be the generator of some one-parameter semigroup or group.} for unbounded $*$-derivations, satisfying different kinds of hypotheses, are investigated. Generation theorems for unbounded $*$-derivations on von Neumann algebras are also investigated in \cite{brattelirob}. Concerning more general Lie groups, the exponentiation theorem \cite[Theorem 3.9]{bratteliheat} is redirected to the more specific context of representations on C$^*$-algebras, in \cite{brattelidissipative}. It should be mentioned that references \cite{bratteliheat}, along with \cite{jorgensenmoore} and \cite{albanese}, have been three of the most inspiring works for this paper.\\

The main objective of this work is to report some new results regarding the exponentiation of (in general, noncommutative) finite-dimensional real Lie algebras of linear operators acting on complete Hausdorff locally convex spaces, focusing on the equicontinuous case, and to search for applications within the realm of locally convex algebras. To the knowledge of the author, there are very few theorems in the literature dealing with the exponentiability of Lie algebras of dimension $d > 1$ of linear operators on locally convex spaces beyond the context of Banach spaces (\cite[Theorem 9.1, page 196]{jorgensenmoore} would be one of them). There exist results for $d = 1$ (see \cite{albanese}, \cite{babalola}, \cite{komura}, \cite{ouchi} and \cite[page 234]{yosida}, for example), but the technical difficulties to prove exponentiability in this context for $d > 1$ are considerably more severe. The main exponentiation theorem of this paper is Theorem 2.14, whose proof is divided into three steps. The first one is to show how to construct a group invariant dense C$^\infty$ domain from a mere dense C$^\infty$ domain. To this end, techniques developed in Chapters 5, 6 and 7 of \cite{jorgensenmoore} in the Banach space context must be suitably adapted. The second step consists in formulating ``locally convex equicontinuous versions'' of three exponentiation theorems found in the literature - \cite[Theorem 9.2]{jorgensenmoore}, \cite[Theorem 3.1]{jorgensengoodman} and \cite[Theorem 3.9]{bratteliheat}. Considerable upgrades on the last reference are made: besides the extension to the vastly more general framework of complete Hausdorff locally convex spaces, the role of the generalized Laplacian employed in \cite{bratteliheat} is clarified - it is substituted by an \textbf{arbitrary strongly elliptic operator} in the (complexification of the) universal enveloping algebra of the operator Lie algebra under consideration. In the final step, necessary conditions for exponentiability are obtained in Theorem 2.13, and in the main theorem of the paper (Theorem 2.14), a \textbf{characterization of exponentiability} in complete Hausdorff locally convex spaces in the same spirit as in \cite[Theorem 3.9]{bratteliheat} is given, with the exception that, in the present work, \textbf{arbitrary strongly elliptic operators} are allowed. Since the subjacent Lie group representations are always locally equicontinuous, this theorem gives, in particular, necessary conditions for exponentiation with respect to \textbf{compact Lie groups}. Finally, in Section 3, some applications to complete locally convex algebras are given (Theorem 3.4), with special attention to locally C$^*$-algebras (Theorems 3.6 and 3.7).

\section{Some General Facts}

\subsection{One-Parameter Semigroups and Groups}

\indent

In the classical theory of one-parameter semigroups and groups, a strongly continuous one-parameter semigroup on a normed space $(\mathcal{Y}, \|\, \cdot \,\|)$ is a family of continuous linear operators $\left\{V(t)\right\}_{t \geq 0}$ satisfying $$V(0) = I, \, V(s + t) = V(s)V(t), \qquad s, t \geq 0$$ and $$\lim_{t \rightarrow t_0} \|V(t)y - V(t_0)y\| = 0, \qquad t_0 \geq 0, \, y \in \mathcal{Y}.$$ It is a well-known fact that in this setting there exist $M > 0$ and $a \geq 0$ satisfying $$(1.1.1) \qquad \|V(t)y\| \leq M \exp (at) \, \|y\|,$$ for all $y \in \mathcal{Y}$, as a consequence of the Uniform Boundedness Principle - see \cite[Proposition 5.5]{engel}. If $a$ and $M$ can be chosen to be, respectively, 0 and 1, then it is called a \textbf{contraction semigroup}. All definitions are analogous for groups, switching from ``$t \geq 0$'' to ``$t \in \mathbb{R}$'' and ``$\exp (at)$'' to ``$\exp (a|t|)$''. Also, in the context of one-parameter groups, if $a$ and $M$ can be chosen to be, respectively, 0 and 1, then it is called a \textbf{group of isometries}, or a representation by isometries.\\

The \textbf{type} of a strongly continuous semigroup $t \longmapsto V(t)$ on a Banach space $(\mathcal{Y}, \|\, \cdot \,\|)$ is defined as the number $$\inf_{t > 0} \, \frac{1}{t} \, \text{log }\|V(t)\|$$ - see \cite[page 306]{hille}; see also \cite[Definition 5.6, page 40]{engel}. Analogously, if $t \longmapsto V(t)$ is a strongly continuous group, its type is defined as $$\inf_{t \in \mathbb{R} \backslash \left\{0\right\}} \frac{1}{|t|} \, \text{log }\|V(t)\|.$$

The task, now, will be to formulate analogous concepts for locally convex spaces.

A family of seminorms $\Gamma := \left\{p_\lambda\right\}_{\lambda \in \Lambda}$ defined on a locally convex space $\mathcal{X}$ is said to be \textbf{saturated} if, for any given finite subset $F$ of $\Lambda$, the seminorm defined by $$p_F \colon x \longmapsto \max \left\{p(x): p \in F\right\}$$ also belongs to $\Gamma$. Every fundamental system of seminorms can always be enlarged to a saturated one by including the seminorms $p_F$, as defined above, in such a way that the resulting family generates the same topology. \textbf{Hence, it will always be assumed in this paper that the families of seminorms to be considered are already saturated, whenever convenient, without further notice.}\\

\textbf{Definition (Equicontinuous Sets):} \textit{Let $\mathcal{X}$ be a Hausdorff locally convex space with a fundamental system of seminorms $\Gamma$ and denote by $\mathcal{L}(\mathcal{X})$ the vector space of continuous linear operators defined on all of $\mathcal{X}$. A set $\Phi \subseteq \mathcal{L}(\mathcal{X})$ of linear operators is called \textbf{equicontinuous} if for every neighborhood $V$ of the origin of $\mathcal{X}$ there exists another neighborhood $U \subseteq \mathcal{X}$ of 0 such that $T[U] \subseteq V$, for every $T \in \Phi$ - equivalently, if for every $p \in \Gamma$ there exist $q \in \Gamma$ and $M_p > 0$ satisfying $$p(T(x)) \leq M_p \, q(x),$$ for all $T \in \Phi$ and $x \in \mathcal{X}$.}\\

\textbf{Definition (Equicontinuous and Exponentially Equicontinuous Semigroups and Groups):} \textit{A one-parameter semigroup on a Hausdorff locally convex space $(\mathcal{X}, \Gamma)$ is a family of continuous linear operators $\left\{V(t)\right\}_{t \geq 0}$ satisfying $$V(0) = I, \, V(s + t) = V(s)V(t), \qquad s, t \geq 0.$$ If, in addition, the semigroup satisfies the property that $$\lim_{t \rightarrow t_0} p(V(t)x - V(t_0)x) = 0, \qquad t_0 \geq 0, \, p \in \Gamma, \, x \in \mathcal{X},$$ then $\left\{V(t)\right\}_{t \geq 0}$ is called \textbf{strongly continuous} or, more explicitly, a \textbf{strongly continuous one-parameter semigroup}. Such a semigroup is called \textbf{exponentially equicontinuous} \cite[Definition 2.1]{albanese} if there exists $a \geq 0$ satisfying the following property: for all $p \in \Gamma$ there exist $q \in \Gamma$ and $M_p > 0$ such that $$p(V(t)x) \leq M_p \, \text{exp}(at) \, q(x), \qquad t \geq 0, \, x \in \mathcal{X}$$ or, in other words, if the rescaled family $$\left\{\exp (-at) \, V(t)\right\}_{t \geq 0} \subseteq \mathcal{L}(\mathcal{X})$$ is equicontinuous \cite[Definition 2.1]{albanese}. If $a$ can be chosen equal to 0, such a semigroup will be called \textbf{equicontinuous}. A strongly continuous semigroup $t \longmapsto V(t)$ is said to be \textbf{locally equicontinuous} if, for every compact $K \subseteq [0, + \infty)$, the set $\left\{V(t): t \in K\right\}$ is equicontinuous. (As was already mentioned before, every strongly continuous semigroup $V$ on a Banach space $\mathcal{Y}$ satisfies $\|V(t)y\| \leq M \exp (at) \, \|y\|$, for all $t \geq 0$ and $y \in \mathcal{Y}$, so $V$ is automatically locally equicontinuous) All definitions are analogous for one-parameter groups, switching from ``$\,t \geq 0$'' to ``$\,t \in \mathbb{R}$'', ``$\,\exp (at)$'' to ``$\,\exp (a|t|)$'' and $[0, + \infty)$ to $\mathbb{R}$.}\\

In the group case, much of the above terminology can be adapted from one-parameter groups to general Lie groups. For example, if $G$ is a Lie group with unit $e$, then a family of continuous linear operators $\left\{V(g)\right\}_{g \in G}$ satisfying $$V(e) = I, \, V(gh) = V(g)V(h), \qquad g, h \in G$$ and $$\lim_{g \rightarrow h} V(g)x = V(h)x, \qquad x \in \mathcal{X}, \, h \in G,$$ is called a \textbf{strongly continuous representation} of $G$. Such a group representation is called \textbf{locally equicontinuous} if for each compact $K \subseteq G$ the set $\left\{V(g): g \in K\right\}$ is equicontinuous.\\

Note that, in the definition of exponentially equicontinuous semigroups, the seminorm $q$ does not need to be equal to $p$. This phenomenon motivates the definition of $\Gamma$-semigroups:\\

\textbf{Definition ($\bm{\Gamma}$-semigroups and $\bm{\Gamma}$-groups):} \textit{Let $(\mathcal{X}, \Gamma)$ be a Hausdorff locally convex space - frequently, the notations $(\mathcal{X}, \Gamma)$ or $(\mathcal{X}, \tau)$ will be used to indicate, respectively, that $\Gamma$ is a fundamental system of seminorms for $\mathcal{X}$ or that $\tau$ is the locally convex topology of $\mathcal{X}$. For each $p \in \Gamma$, define $$V_p := \left\{x \in \mathcal{X}: p(x) \leq 1\right\}.$$ Following \cite{babalola}, the following conventions will be used:}

\begin{enumerate}

\item \textit{$\mathcal{L}_\Gamma(\mathcal{X})$ denotes the family of linear operators $T$ on $\mathcal{X}$ satisfying the property that, for all $p \in \Gamma$, there exists $\lambda(p, T) > 0$ such that $$T[V_p] \subseteq \lambda(p, T) \, V_p$$ or, equivalently, $$p(T(x)) \leq \lambda(p, T) \, p(x), \qquad x \in \mathcal{X}.$$}

\item \textit{A strongly continuous semigroup $t \longmapsto V(t)$ is said to be a $\bm{\Gamma}$\textbf{-semigroup}\footnote{In \cite{babalola} they are called (C$_0$, 1) semigroups - similarly for groups.} if $V(t) \in \mathcal{L}_\Gamma(\mathcal{X})$, for all $t \geq 0$ and, for every $p \in \Gamma$ and $\delta > 0$, there exists a number $\lambda = \lambda(p, \left\{V(t): 0 \leq t \leq \delta\right\}) > 0$ such that $$p(V(t)x) \leq \lambda \, p(x),$$ for all $0 \leq t \leq \delta$ and $x \in \mathcal{X}$. Analogously, a strongly continuous group $t \longmapsto V(t)$ is said to be a $\bm{\Gamma}$\textbf{-group} if $V(t) \in \mathcal{L}_\Gamma(\mathcal{X})$, for all $t \in \mathbb{R}$ and, for every $p \in \Gamma$ and $\delta > 0$, there exists a number $\lambda = \lambda(p, \left\{V(t): |t| \leq \delta\right\}) > 0$ such that $$p(V(t)x) \leq \lambda \, p(x),$$ for all $|t| \leq \delta$ and $x \in \mathcal{X}$ (see \cite[Definitions 2.1 and 2.2]{babalola} - note that a ``local equicontinuity-type'' requirement, which is present in Definition 2.1, is incorrectly missing in Definition 2.2). An equivalent way of defining then is the following: a strongly continuous one-parameter semigroup $t \longmapsto V(t)$ is said to be a $\bm{\Gamma}$\textbf{-semigroup} if, for each $p \in \Gamma$, there exist $M_p, \sigma_p \in \mathbb{R}$ such that $p(V(t)x) \leq M_p \, e^{\sigma_p t} p(x)$, for all $x \in \mathcal{X}$ and $t \geq 0$ - see \cite[Theorem 2.6]{babalola}. $\bm{\Gamma}$\textbf{-groups} are defined in an analogous way, but with $p(V(t)x) \leq M_p \, e^{\sigma_p |t|} p(x)$, for all $x \in \mathcal{X}$ and $t \in \mathbb{R}$. Note that these definitions automatically imply local equicontinuity of the one-parameter (semi)group $V$, and that the operators $V(t)$ have the (KIP) with respect to $\Gamma$.}

\item \textit{By a procedure which will become clearer in the future and which is closely related to the notion of Kernel Invariance Property - this will be properly defined, soon -, the author of \cite{babalola} associates to the $\Gamma$-semigroup $t \longmapsto V(t)$ a net of strongly continuous semigroups $\left\{\tilde{V}_p\right\}_{p \in \Gamma}$, each one of them being defined on a Banach space $(\mathcal{X}_p, \|\, \cdot \,\|_p)$. By what was already said above, for each $p \in \Gamma$, the number $$\inf_{t > 0} \, \frac{1}{t} \, \text{log }\|\tilde{V}_p(t)\|_p$$ is well-defined, and will be denoted by $w_p$. The family $\left\{w_p\right\}_{p \in \Gamma}$ is called \textbf{the type of $\bm{V}$}. If $w := \sup_{p \in \Gamma} w_p < \infty$, then $V$ is said to be of \textbf{bounded type} $\bm{w}$, following \cite[page 170]{babalola} - analogously, substituting ``$t > 0$'' by ``$t \in \mathbb{R} \backslash \left\{0\right\}$'' and ``$1/t$'' by ``$1/|t|$'', one defines $\Gamma$-groups of bounded type. A very nice aspect of $\Gamma$-semigroups of bounded type $w$ is that they satisfy the resolvent formula $$\text{R}(\lambda, A)x = \int_0^{+ \infty} e^{-\lambda t} V(t)x \, dt, \qquad x \in \mathcal{X},$$ for all $\lambda > w$, where $A$ is the infinitesimal generator of $V$, as is shown in \cite[Theorem 3.3]{babalola} - for the definition of infinitesimal generator, see the next subsection. This formula will be the key for many proofs to come, like the ones of Theorem 2.4 and Theorem 2.12.}

\end{enumerate}

\textbf{Note that all of the definitions regarding semigroups and groups, above, with the exception of the last one (regarding $\bm{\Gamma}$-semigroups and $\bm{\Gamma}$-groups), are \underline{independent} of the fixed fundamental system of seminorms $\bm{\Gamma}$.}\\

It should be emphasized, at this point, that the main object of study in this work is not the structure of locally convex spaces in itself, but rather families of linear operators defined on these spaces, like generators of one-parameter semigroups and groups of continuous operators and Lie algebras of linear operators, for example. These operators of interest are, in most cases, not assumed to be continuous, nor are they everywhere defined. Hence, if $\mathcal{X}$ is a locally convex space, a linear operator $T$ on $\mathcal{X}$ will always be assumed as being defined on a vector subspace $\text{Dom }T$ of $\mathcal{X}$, called its \textbf{domain}, and in case $\mathcal{X}$ is also an algebra (respectively, a $*$-algebra), $\text{Dom }T$ will, by definition, be a subalgebra (respectively, a $*$-subalgebra - in other words, a subalgebra which is closed under the operation $*$ of involution). When its domain is dense, the operator will be called \textbf{densely defined}. Also, if $S$ and $T$ are two linear operators then their sum will be defined by $$\text{Dom }(S + T) := \text{Dom }S \cap \text{Dom }T, \qquad (S + T)(x) := S(x) + T(x), \qquad x \in \text{Dom }(S + T)$$ and their composition by $$\text{Dom }(TS) := \left\{x \in \mathcal{X}: S(x) \in \text{Dom }T\right\}, \qquad (TS)(x) := T(S(x)), \qquad x \in \text{Dom }(TS),$$ following the usual conventions of the classical theory of unbounded linear operators on Hilbert and Banach spaces. The range of $T$ will be denoted by $\text{Ran }T$\\

A very important concept is that of a resolvent set, so let $T \colon \text{Dom }T \longrightarrow \mathcal{X}$ be a linear operator on $\mathcal{X}$. Then, the set $$\rho(T) := \left\{\lambda \in \mathbb{C}: (\lambda I - T) \colon \text{Dom }T \subseteq \mathcal{X} \longrightarrow \mathcal{X} \text{ is bijective and } (\lambda I - T)^{-1} \in \mathcal{L}(\mathcal{X})\right\}$$ is called the \textbf{resolvent set} of the operator $T$ (note that it may happen that $\rho(T) = \emptyset$, or even $\rho(T) = \mathbb{C}$), and the map $$\text{R}(\, \cdot \, , T) \colon \rho(T) \ni \lambda \longmapsto \text{R}(\lambda, T) := (\lambda I - T)^{-1} \in \mathcal{L}(\mathcal{X})$$ is called the \textbf{resolvent of $\bm{T}$}. When $\lambda \in \rho(T)$, it will be said that \textbf{the resolvent operator of $\bm{T}$ exists at $\bm{\lambda}$}, and $\text{R}(\lambda, T)$ is called the \textbf{resolvent operator of $\bm{T}$ at $\bm{\lambda}$} or the \textbf{$\bm{\lambda}$-resolvent} of $T$. Finally, the complementary set $$\sigma(T) := \mathbb{C} \, \backslash \rho(T)$$ is called the \textbf{spectrum} of the operator $T$ (see page 258 of \cite{albanesemontel}).\\

Now, a central notion for this work is going to be defined: that of an infinitesimal generator.

\subsection{Lie Group Representations and Infinitesimal Generators}

\indent

\textbf{Definition (Infinitesimal Generators):} \textit{Given a strongly continuous one-parameter semigroup $t \longmapsto V(t)$ on $(\mathcal{X}, \Gamma)$, consider the subspace of vectors $x \in \mathcal{X}$ such that the limit $$\lim_{t \rightarrow 0} \frac{V(t)x - x}{t}$$ exists. Then, the linear operator $dV$ defined by $$\text{Dom }dV := \left\{x \in \mathcal{X}: \lim_{t \rightarrow 0} \frac{V(t)x - x}{t} \text{ exists in } \mathcal{X}\right\}$$ and $dV(x) := \lim_{t \rightarrow 0} \frac{V(t)x - x}{t}$ is called the \textbf{infinitesimal generator} (or, simply, the generator) of the semigroup $t \longmapsto V(t)$ (the definition for groups is analogous). Also, if $G$ is a Lie group then, for each fixed element $X$ of its Lie algebra $\mathfrak{g}$, the infinitesimal generator of the one-parameter group $t \longmapsto V(\exp tX)$ will be denoted by $dV(X)$ ($\exp$ denotes the exponential map of the Lie group $G$).}\\

Two very important results regarding infinitesimal generators on locally convex spaces are Propositions 1.3 and 1.4 of \cite{komura}, which prove that infinitesimal generators of strongly continuous locally equicontinuous one-parameter semigroups on sequentially complete locally convex spaces are densely defined and closed.\footnote{A linear operator $T \colon \text{Dom T} \subseteq \mathcal{X} \longrightarrow \mathcal{X}$ is \textbf{closed} if its graph is a closed subspace of $\mathcal{X} \times \mathcal{X}$. An operator $S \colon \text{Dom S} \subseteq \mathcal{X} \longrightarrow \mathcal{X}$ is \textbf{closable} if it has a closed extension or, equivalently, if for every net $\left\{x_\alpha\right\}_{\alpha \in \mathcal{A}}$ in $\text{Dom }S$ such that $x_\alpha \longrightarrow 0$ and $S(x_\alpha) \longrightarrow y$, one has $y = 0$. If $S$ is closable, then it has a minimal closed extension, called the \textbf{closure} of $S$ and denoted by $\overline{S}$ - see \cite[page 250]{reedsimon1}.}

A closable linear operator with the property that its closure is an infinitesimal generator is called an \textbf{infinitesimal pregenerator}.\\

\textbf{Definition (Smooth and Analytic Vectors):} \textit{Returning to the original subject of this section, if $\mathcal{X}$ is a Hausdorff locally convex space and $V \colon G \longrightarrow \mathcal{L}(\mathcal{X})$ is a strongly continuous representation then a vector $x \in \mathcal{X}$ is called a \textbf{C$\bm{^\infty}$ vector for $\bm{V}$, or a smooth vector for $\bm{V}$}, if the map $G \ni g \longmapsto V(g)x$ is of class $C^\infty$: a map $f \colon G \longrightarrow \mathcal{X}$ is of class $C^\infty$ at $g \in G$ if it possesses continuous partial derivatives of all orders with respect to a chart around $g$. If $f$ is of class $C^\infty$ at all points $g \in G$, $f$ is said to be of class $C^\infty$ on $G$.\footnote{This definition of a smooth vector may be found in \cite[page 47]{moore}. It should be mentioned, however, that defining a smooth vector $x$ by asking that $G \ni g \longmapsto V(g)x$ must be a smooth map with respect to the \textbf{weak topology}, instead, would be ``more natural'', in a certain sense: the verification of chart independence of this definition, for example, amounts to the usual proof one finds in \textbf{finite dimensional} Differential Geometry. In a \textbf{wafer complete} Hausdorff locally convex space - this definition will be introduced, soon - these two notions of smoothness coincide, as it is shown in \cite[Lemma 1, page 47]{moore}.} The subspace of smooth vectors for $V$ will be denoted by $C^\infty(V)$. Moreover, following \cite[page 54]{moore}, a vector $x \in \mathcal{X}$ is called \textbf{analytic for $\bm{V}$} if $x \in C^\infty(V)$ and the map $F_x \colon G \ni g \longmapsto V(g)x$ is analytic: in other words, if $x \in C^\infty(V)$ and, for each $g \in G$ and every analytic chart $h \colon g' \longrightarrow (t_k(g'))_{1 \leq k \leq d}$ around $g$ sending it to 0 there exists $r_x > 0$ such that the series $$\sum_{\alpha \in \mathbb{N}^d} \frac{p(\partial^\alpha F_x(g'))}{\alpha!} \, t(g')^\alpha$$ is absolutely convergent to $F_x(g')$, for every $p \in \Gamma$, whenever $|t(g')| < r_x$, where $t(g')^\alpha = t_1(g')^{\alpha_1} \ldots t_d(g')^{\alpha_d}$ and $|t(g')| := \max_{1 \leq k \leq d} |t_k(g')|$ (note that $r_x$ is independent of $p \in \Gamma$).\footnote{Note, also, that \cite[Lemma 3, page 52]{moore} proves that the convergence mentioned holds in \text{any} Hausdorff locally convex topology between the weak and the strong ones, and with the same radius of convergence.} The subspace of analytic vectors for $V$ will be denoted by $C^\omega(V)$. If $\tau$ is the topology defined by $\Gamma$, then the elements of $C^\omega(V)$ will sometimes be called \textbf{$\bm{\tau}$-analytic}. Also, if the Lie group under consideration is $\mathbb{R}$, then the subspace of analytic vectors will be denoted by $C^\omega(T)$, where $T$ is the infinitesimal generator of $V$, and called the \textbf{the subspace of analytic vectors for $\bm{T}$}.}\\

An important observation is that, if $$f \colon g \longmapsto V(g)x$$ is of class C$^\infty$ at $e$, then it is of class C$^\infty$ on all of $G$. In fact, fixing $\sigma \in G$ and defining the function $\phi_\sigma \colon G \ni g \longmapsto \sigma^{-1}g \in G$, then $$f(g) = V(\sigma)[(f \circ \phi_\sigma)(g))],$$ for every $g \in G$. Since $V(\sigma)$ is a continuous linear operator on $\mathcal{X}$, $\phi_\sigma$ is of class C$^\infty$ at $\sigma$ and $f$ is of class C$^\infty$ at $e$, it follows that $f$ is of class C$^\infty$ at $\sigma$. Hence, in order to verify that a vector $x$ is smooth for the representation $V$ it suffices to show that $g \longmapsto V(g)x$ is of class C$^\infty$ at $e$. An analogous observation is valid for analytic vectors.\\

Consider a strongly continuous representation $V$ of the Lie group $G$, whose Lie algebra is $\mathfrak{g}$. Fix a left invariant Haar measure on $G$. There exists an important subspace of $C^\infty(V)$ defined by $$D_G(V) := \text{span}_{\mathbb{C}} \left\{x_V(\phi), \, \phi \in C_c^\infty(G), x \in \mathcal{X}\right\}, \quad \text{where} \quad x_V(\phi) := \int_G \phi(g) \, V(g) x \, dg,$$ known as the G\aa rding subspace of $V$ - see \cite{garding}. To begin with, the meaning of the integrals defining these elements must be specified. This can be accomplished by appropriately adapting a beautiful argument originally formulated within a Banach space context in \cite[page 219]{arveson},\footnote{This argument remains valid if $G$ is only a Hausdorff locally compact topological group.} but one additional hypothesis is needed: $\mathcal{X}$ must be a \textbf{wafer complete} locally convex space or, in other words, a space with the property that the weakly closed convex balanced hull of every weakly compact set is weakly compact (fortunately, by \cite[Corollary 8, page 19]{moore}, every complete locally convex space is wafer complete, so the hypothesis of completeness is sufficient). If $\Psi \colon G \ni g \longmapsto \Psi(g) \in \mathcal{X}$ is a compactly supported weakly continuous function, then the integral $\int_G \Psi(g) \, dg$ exists in $\mathcal{X}$, by definition, if there exists an element $y_\Psi \in \mathcal{X}$ such that, for all $f$ in the topological dual $\mathcal{X}'$ of $\mathcal{X}$, $$f(y_\Psi) = \int_G f(\Psi(g)) \, dg.\footnote{See, also, \cite{jorgensenmoore}, page 443.}$$ Under the hypotheses just stated, existence of the integral $\int_G \Psi(g) \, dg$ in this sense is guaranteed.

$D_G(V)$ is indeed a subspace of $C^\infty(V)$,\footnote{If $\mathcal{X}$ is a Fr\'echet space, then the G\aa rding subspace actually coincides with the space of smooth vectors. For a proof of this fact see, for example, Bill Casselman's essay ``The theorem of Dixmier \& Malliavin'', at his webpage.} as can be seen by an argument adapted from \cite{garding}. Moreover, fix a sequence $\left\{\phi_n\right\}_{n \in \mathbb{N}}$ of nonnegative functions in $C_c^\infty(G)$ with $$\int_G \phi_n \, dg = 1, \qquad n \in \mathbb{N},$$ and with the property that for every open neighborhood $U$ of the origin there exists $n_U \in \mathbb{N}$ such that $\text{supp }\phi_n \subseteq U$ for all $n \geq n_U$. Then, $x_V(\phi_n) \longrightarrow x$ in $\mathcal{X}$, when $n \longrightarrow + \infty$, since $V$ is strongly continuous. \textbf{This establishes that $\bm{D_G(V)}$ and hence also $\bm{C^\infty(V)}$ is dense in $\bm{\mathcal{X}}$.} Note that the above results regarding $D_G(V)$ remain valid if instead of a left invariant Haar measure one had chosen a right invariant one (see \cite[Corollary 8.30, page 534]{knapp}).\\

Fix $x \in C^\infty(V)$, $h \colon G \supseteq V \longrightarrow W \subseteq \mathbb{R}^d$ a chart of $G$ around $e$ and define the function $$f \colon g \longmapsto V(g)x.$$ If $X$ belongs to the Lie algebra of $G$, denoted by $\mathfrak{g}$, then $c \colon t \longmapsto \exp tX$ is an infinitely differentiable curve in $G$ which has $X$ as its tangent vector at $t = 0$. As $f \circ c = (f \circ h^{-1}) \circ (h \circ c)$, the limit $$\lim_{t \rightarrow 0} \frac{V(\exp tX)x - x}{t}$$ exists, by the chain rule for directional derivatives on locally convex spaces. 

It is clear that $$\mathcal{D}_\infty := \bigcap_{k = 1}^d \bigcap_{n = 1}^{+ \infty} \text{Dom }dV (X_k)^n \supseteq C^\infty(V).$$ In order to prove the reverse inclusion, an adapted argument of Theorem 1.1 of \cite{goodman} can be made: fix $x \in \mathcal{D}_\infty$, $f \in \mathcal{X}'$, a right invariant Haar measure $dg$ on $G$ and a $\phi \in C_c^\infty(G)$. Then, denoting by $\tilde{X}_k$ the complete, globally defined and left invariant vector field on $G$ such that $\tilde{X}_k(e) = X_k$, for all $1 \leq k \leq d$,\footnote{Remember that such field is defined by $$\tilde{X}_k(g) := D(L_g)(e)(X_k)$$ for every $g \in G$, where $L_g \colon h \longmapsto gh$ is the operator of left-translation by $g$.} one obtains $$\int_G \phi(g) \, f(V(g) \, dV(X_k)^n(x)) \, dg = \int_G \frac{d^n}{dt^n} \left. \left[\phi(g) \, f(V(g \exp tX_k)x)\right]\right|_{t = 0} \, dg$$ $$= \int_G \frac{d^n}{dt^n} \left. \left[\phi(g \, \text{exp}-tX_k) \, f(V(g)x)\right]\right|_{t = 0} \, dg = (-1)^n \int_G \tilde{X}_k^n(\phi)(g) \, f(V(g)x) \, dg, \qquad n \in \mathbb{N},$$ because $\tilde{X}_k(\phi)(g) = \frac{d}{dt} \left. \phi(g \exp tX_k) \right|_{t = 0}$. Hence, substituting $n$ by $2m$, $m \geq 1$, and repeating the argument for every $1 \leq k \leq d$, one finds that for every chart $h \colon G \supseteq U \longrightarrow W \subseteq \mathbb{R}^d$, where $U$ and $W$ are open sets in their respective topologies and $W$ contains the origin of $\mathbb{R}^d$, the function $\psi \colon W \ni w \longmapsto f(V(h^{-1}(w))x) \in \mathbb{R}$ is a continuous weak solution to the PDE $$\Delta_m(\psi)(w) = \sum_{k = 1}^d f(V(h^{-1}(w)) \, dV(X_k)^{2m}x), \qquad w \in W,$$ where $\Delta_m := \sum_{k = 1}^d dV(X_k)^{2m}$ is an elliptic operator. Since $x \in \mathcal{D}_\infty$, $m$ can be chosen arbitrarily large. Therefore, by \cite[Theorem 1, page 190]{bers}, $\psi$ is of class C$^\infty$ on $W$. By the arbitrariness of $h$, it follows that $g \longmapsto f(V(g)x)$ is of class C$^\infty$ on $G$. But, then again, since $f$ is arbitrary, one may use \cite[Lemma 1, page 47]{moore} to show that the map $g \longmapsto V(g)x$ is of class C$^\infty$. This shows the inclusion $$\mathcal{D}_\infty \subseteq C^\infty(V),$$ proving that these sets are equal. A straightforward consequence of this equality is that $C^\infty(V)$ is left invariant by all of the $dV(X)$: write $X = \sum_{k = 1}^d a_k X_k$, for some real numbers $a_k$. Then, $$dV(X)[C^\infty(V)] \subseteq \sum_{k = 1}^d a_k \, dV(X_k)[C^\infty(V)] = \sum_{k = 1}^d a_k \, dV(X_k)[\mathcal{D}_\infty] \subseteq \mathcal{D}_\infty = C^\infty(V).$$

Summarizing, for a strongly continuous Lie group representation $g \longmapsto V(g)$ on a Hausdorff locally convex space, the following are true:

\begin{enumerate}

\item the limit $$\lim_{t \rightarrow 0} \frac{V(\exp tX)x - x}{t}$$ exists, for every $x \in C^\infty(V)$ and $X \in \mathfrak{g}$;

\item if $V \colon G \longrightarrow \mathcal{L}(\mathcal{X})$ is a strongly continuous locally equicontinuous representation of the Lie group $G$ and $\mathcal{X}$ is sequentially complete, then the generator $dV(X)$ is a closed, densely defined linear operator on $\mathcal{X}$;

\item if $\mathcal{X}$ is wafer complete, $C^\infty(V)$ is dense in $\mathcal{X}$;

\item if $\mathcal{X}$ is wafer complete, $C^\infty(V)$ is left invariant by the operators $dV(X)$, for all $X \in \mathfrak{g}$.

\end{enumerate}

Therefore, when working with a strongly continuous Lie group representation $V$ on a Hausdorff locally convex space $\mathcal{X}$, some ``good'' basic general hypotheses which may be imposed a priori are that $\mathcal{X}$ is \textbf{complete} and that $V$ is \textbf{locally equicontinuous}.

\subsection{Lie Algebra Representations Induced by Group Representations}

\indent

Assume $\mathcal{X}$ is a complete Hausdorff locally convex space. Another important fact is that the application $$\partial V \colon \mathfrak{g} \ni X \longmapsto \partial V(X) \in \text{End}(C^\infty(V)),$$ where $\partial V(X) := dV(X)|_{C^\infty(V)}$ and $\text{End}(C^\infty(V))$ denotes the algebra of all endomorphisms on $C^\infty(V)$ or, in other words, the linear operators defined on $C^\infty(V)$ with ranges contained in $C^\infty(V)$, is a representation of Lie algebras, so a strongly continuous representation of a Lie group always induces a representation of its Lie algebra (note that items 1 and 4, above, show that this application is well-defined). This Lie algebra representation is called the \textbf{infinitesimal representation of $\bm{V}$}. To see that it preserves commutators, fix $f \in \mathcal{X}'$ and $x \in C^\infty(V)$. Then, $$f([\partial V(X) \partial V(Y) - \partial V(Y) \partial V(X)]x)$$ $$= \frac{d}{dt} \left.\left[\frac{d}{ds} \left.\left[f(V(\exp tX) \, V(\exp sY)x)\right]\right|_{s = 0} \right]\right|_{t = 0} - \, \frac{d}{ds} \left.\left[\frac{d}{dt} \left.\left[f(V(\exp sY) \, V(\exp tX)x)\right]\right|_{s = 0} \right]\right|_{t = 0}$$ $$= [\tilde{X} \, \tilde{Y} - \tilde{Y} \, \tilde{X}][f(V( \, \cdot \, )x)](e) = \widetilde{[X, Y]} \, [f(V( \, \cdot \, )x)](e)$$ $$= \frac{d}{dt} \left.f(V(\exp t[X, Y])x)\right|_{t = 0} = f(\partial V([X, Y])x).$$ Since $f$ is arbitrary, a corollary of Hahn-Banach's Theorem gives $$[\partial V(X) \partial V(Y) - \partial V(Y) \partial V(X)](x) = \partial V([X, Y])(x),$$ and arbitrariness of $x$ establishes the result.\footnote{This proof is an adaptation of the argument in \cite[Proposition 10.1.6, page 263]{schmudgen}.} The linearity follows by an analogous argument.

Actually, defining $$\partial V(Y_1 \ldots Y_n) := \partial V(Y_1) \ldots \partial V(Y_n),$$ for $n \in \mathbb{N}$ and $Y_i \in \mathfrak{g}$, and extending this definition by linearity, one can define a unital homomorphism between the universal enveloping algebra of $\mathfrak{g}$, $\mathfrak{U}(\mathfrak{g})$ (which is, roughly speaking, the unital associative algebra formed by the real ``noncommutative polynomials'' in the elements of $\mathfrak{g}$), and a unital associative subalgebra of $\text{End}(C^\infty(V))$ that extends the original Lie algebra homomorphism $X \longmapsto \partial V(X)$ (the same notation was used to represent both the original morphism and its extension to $\mathfrak{U}(\mathfrak{g})$). An analogous definition may be done with the operators $\partial V(Y_k)$ replaced by $dV(Y_k)$.\\

The results proved so far are motivators to define $C^\infty(V)$ as the \textit{domain of the representation $\partial V$}, and will be denoted by $\text{Dom }(\partial V)$. Also, it is clear that, if $C^\infty(V)$ is substituted by a subspace $\mathcal{D}$ of $C^\infty(V)$ which is left invariant by all of the generators $dV(X)$, $X \in \mathfrak{g}$, then it is also possible to define such infinitesimal representation by $\partial V(X) := dV(X)|_\mathcal{D}$.

\subsection{Group Invariance and Cores}

\indent

The domain $C^\infty(V)$ also has some nice properties with respect to the group representation $V$. For example, $C^\infty(V)$ is also left invariant by all of the operators $V(g)$, because if $x \in C^\infty(V)$, then the application $G \ni w \longmapsto V(w)V(g)x \in \mathcal{X}$ is the composition of two functions of class C$^\infty$: $w \longmapsto V(w)x$ and $w \longmapsto wg$.\\

\textbf{Definition (Projective C$\bm{^\infty}$-Topology on Space of Smooth Vectors):} \textit{Fix an ordered basis $\mathcal{B} := \left(X_k\right)_{1 \leq k \leq d}$ for $\mathfrak{g}$, choose a fundamental system of seminorms $\Gamma$ for $\mathcal{X}$ and equip $C^\infty(V)$ with the topology defined by the family $$\left\{\rho_{p, n}: p \in \Gamma, n \in \mathbb{N}\right\},$$ where $$\rho_{p, 0}(x) := p(x), \qquad dV(X_0) := I$$ and $$\rho_{p, n}(x) := \max \left\{p(dV(X_{i_1}) \ldots dV(X_{i_n})x): 0 \leq i_j \leq d\right\}$$ - this is called the \textbf{projective C$\bm{^\infty}$-topology} on $C^\infty(V)$, and it does not depend upon neither the fixed basis $\mathcal{B}$ nor on the particular $\Gamma$. If $\mathcal{X}$ is sequentially complete and $V$ is locally equicontinuous, then each generator $dV(X_k)$ is closed, as mentioned before, and $C^\infty(V)$ becomes a complete Hausdorff locally convex space - to prove this, adapt the argument of \cite[Corollary 1.1]{goodman}, exploring the closedness of the operators $dV(X_k)$, $1 \leq k \leq d$.}\\

Similarly, for each fixed $1 \leq k < \infty$, the subspace $$C^k(V) := \bigcap_{n = 1}^k \bigcap \left\{\text{Dom }[dV(X_{i_1}) \ldots dV(X_{i_j}) \ldots dV(X_{i_n})]: X_{i_j} \in \mathcal{B}\right\},$$ called the subspace of C$^k$ vectors for the representation $V$, is a sequentially complete (respectively, complete) Hausdorff locally convex space when equipped with the C$^k$-topology generated by the family $$\left\{\rho_{p, k}: p \in \Gamma\right\}$$ of seminorms, if $\mathcal{X}$ is sequentially complete (respectively, complete) and $V$ is locally equicontinuous. Moreover, the operators $$V_\infty(g) := V(g)|_{C^\infty(V)} \in \text{End}(C^\infty(V)), \qquad g \in G,$$ are continuous with respect to the projective topology (use an inductive argument together with the identity $$g \exp tX \, g^{-1} = \exp (\text{Ad}(g)(tX)), \qquad g \in G, \, X \in \mathfrak{g},$$ to prove this - see \cite[page 31]{schmudgen}).\\

\textbf{Definition (Cores):} \textit{If $T \colon \text{Dom }T \subseteq \mathcal{X} \longrightarrow \mathcal{X}$ is a closed linear operator and $\mathcal{D} \subseteq \text{Dom }T$ is a linear subspace of $\text{Dom }T$ such that $$\overline{T|_\mathcal{D}} = T,$$ then $\mathcal{D}$ is called a \textbf{core} for $T$.}\\

Before proceeding, a few words about the use of some notations must be said. Consider a Lie algebra $\mathfrak{g}$ with an ordered basis $\left(B_k\right)_{1 \leq k \leq d}$. Throughout the whole manuscript, whenever a monomial $B^\mathsf{u}$ of \textbf{size} $n \geq 1$ in the elements of $\left(B_k\right)_{1 \leq k \leq d}$ is considered, it is understood that $\mathsf{u} = (\mathsf{u}_j)_{j \in \left\{1, \ldots, n\right\}}$ is a function from $\left\{1, \ldots, n\right\}$ to $\left\{1, \ldots, d\right\}$ and $B^\mathsf{u} := B_{\mathsf{u}_1} \ldots B_{\mathsf{u}_n}$ is an element of the complexified universal enveloping algebra $(\mathfrak{U}[\mathfrak{g}])_\mathbb{C} := \mathfrak{U}[\mathfrak{g}] + i \, \mathfrak{U}[\mathfrak{g}]$. Therefore, $B^\mathsf{u}$ can be thought of as an unordered monomial, or a noncommutative monomial.\footnote{The expression ``unordered'' refers to the fact that the sequence $(\mathsf{u}_1, \ldots, \mathsf{u}_n)$ of numbers is not (necessarily) non-decreasing.} The size of the monomial is defined to be $|\mathsf{u}| = n$ so, roughly speaking, it is the size of the ``word'' $B_{\mathsf{u}_1} \ldots B_{\mathsf{u}_n}$. For a matter of convenience, the size of a linear combination of monomials is defined to be the size of the ``biggest'' monomial composing the sum. However, the following question arises: could this sum be written in another way, so that the associated size is different? It may very well be the case that the basis elements share an extra relation, thus causing this kind of ambiguity: for example, if $X$ and $Y$ are elements of a Lie algebra with basis $\left\{X, Y, [X, Y]\right\}$, then $[X, Y] = XY - YX$ has size 2, but also size 1, depending on how it is written. Therefore, the concept of size is not an intrinsic one, and is intimately related to how the sum was written. The same cannot be said about the definition of \textbf{order} of an element in $\mathfrak{U}[\mathfrak{g}]$: by the Poincar\'e-Birkhoff-Witt Theorem \cite[Theorem 3.8, page 217]{knapp} together with an inductive procedure on the argument sketched at the beginning of Subsection 1.7, one sees that every linear combination of monomials $P$ can be written in a unique way as a polynomial on the variables $\left(B_k\right)_{1 \leq k \leq d}$ of the form $$\sum_{|\alpha| \leq m} c_\alpha \, B_1^{\alpha_1} \ldots B_k^{\alpha_k} \ldots B_d^{\alpha_d},$$ where $c_{\alpha'} \neq 0$ for some multi-index $\alpha'$ of order $m$. In this case, the order of $P$ is defined to be $|\alpha'| = m$, the order of the multi-index $\alpha'$, and it does \textbf{not} depend on the particular choice of the basis. Note that the notations were carefully chosen, in order to avoid confusions. The symbols $\mathsf{u}$, $\mathsf{v}$ and $\mathsf{w}$ will always be employed to indicate the use of an unordered monomial, while the greek letters $\alpha$, $\beta$ and $\gamma$ will denote multi-indices, following the usual convention for ordered monomials. The use of the noncommutative monomial notation will occur in Subsection 1.7, Theorem 2.5 and in the proof of some of the exponentiation theorems. As an abbreviation, they will usually be referred to as ``monomials'', rather than ``noncommutative monomials''.\\

A statement which is similar to that of Theorem 1.4.1, below, is claimed in \cite[Theorem B.5, page 446]{jorgensenmoore} - see also \cite[Corollary 1.2]{poulsen} and \cite[Corollary 1.3]{poulsen}. Moreover, it is claimed there under the hypothesis of sequential completeness, but the author of the present work found it necessary to add the hypothesis that it is also wafer complete - due to a pragmatical reason, the hypothesis of completeness, instead of wafer-completeness + sequential completeness, was assumed by the author:\\

\textbf{Theorem 1.4.1:} \textit{Let $g \longmapsto V(g)$ be a strongly continuous locally equicontinuous representation of a Lie group $G$, with Lie algebra $\mathfrak{g}$, on a complete Hausdorff locally convex space $(\mathcal{X}, \Gamma)$, whose topology will be denoted by $\tau$. Suppose that $\mathcal{D}$ is a dense subspace of $\mathcal{X}$ which is a closed subspace of $C^\infty(V)$ when equipped with the projective C$^\infty$ topology $\tau_\infty$. Then, if $$V(g)[\mathcal{D}] \subseteq \mathcal{D}, \qquad g \in G,$$ that is, if $\mathcal{D}$ is group invariant, $\mathcal{D}$ must be a core for every closed operator $$\overline{dV(X_{i_1}) \ldots dV(X_{i_n})|_{C^n(V)}}$$ in the basis elements $\left\{dV(X_k)\right\}_{1 \leq k \leq d}$, for all $n \in \mathbb{N}$, $n \geq 1$. More generally, if $L \in \mathfrak{U}(\mathfrak{g})$ is an element of the universal enveloping algebra of size at most $n$, then $\mathcal{D}$ is a core for $\overline{dV(L)|_{C^n(V)}}$, and if $L \in \mathfrak{U}(\mathfrak{g})$ is any element, then $\mathcal{D}$ is a core for $\overline{dV(L)|_{C^\infty(V)}}$.}\\

\textbf{Proof of Theorem 1.4.1:} To see $g \longmapsto \phi(g) \, V(g)x$ is a compactly supported weakly-continuous function \textbf{from $\bm{G}$ to $\bm{\mathcal{D}}$}, for all $x \in \mathcal{D}$ and $\phi \in C_c^\infty(G)$, it will be proved that $g \longmapsto V(g)x$ is a smooth function from $G$ to $(\mathcal{D}, \tau_\infty)$. Following the idea of \cite[Proposition 1.2]{poulsen}, let $X_{i_1} \ldots X_{i_n}$ be a monomial in the basis elements. Then, it is sufficient to show $g \longmapsto \partial V(X_{i_1} \ldots X_{i_n}) V(g)x$ is a smooth function from $G$ to $(\mathcal{X}, \tau)$. Consider the analytic coordinate system $h(t) \longmapsto t$ around the identity $e$ of $G$ defined by $$h(t) := \exp (t_1 \, X_1) \ldots \exp (t_d \, X_d), \qquad t := (t_k)_{1 \leq k \leq d} \in \mathbb{R}^d.$$ The function $(s, t) \longmapsto V(h(s) h(t))x$ is C$^\infty$ from $\mathbb{R}^{2d}$ to $\mathcal{X}$, because $C^\infty(V)$ is left invariant by $V$ and by the generators $dV(X_k)$, $1 \leq k \leq d$. Therefore, since $$\partial V(X_{i_1} \ldots X_{i_n}) V(h(t))x = \left. \frac{\partial}{ds_{i_1}} \ldots \frac{\partial}{ds_{i_n}} V(h(s) h(t))x \right|_{s = 0}, \qquad t \in \mathbb{R}^d,$$ the result follows.

Hence, in particular, if $\left\{g_m\right\}_{m \in \mathbb{N}}$ is a sequence converging to $g_0 \in G$, the sequence $$\left\{dV(X_{i_1}) \ldots dV(X_{i_n}) \, V(g_m)x\right\}_{m \in \mathbb{N}}$$ converges to $dV(X_{i_1}) \ldots dV(X_{i_n}) \, V(g_0)x$, and so $$f'(V(g_\alpha)x) \longrightarrow f'(V(g_0)x),$$ for every $f' \in \mathcal{D}'$, by the definition of the projective C$^\infty$-topology. This means $g \longmapsto V(g)x$ is a weakly continuous function from $G$ to $\mathcal{D}$. Therefore, the observations made in the paragraphs where the G\"arding domain was defined guarantee the existence of the integral $\int_G \phi(g) \, V(g)x \, dg$ as an element of $\mathcal{D}$ (note that the hypotheses of group invariance and of wafer completeness - actually, completeness - of $\mathcal{D}$ with respect to the projective topology were essential).\\

For the next step, it is necessary to see that if $x \in \mathcal{X}$ and $\phi \in C_c^\infty(G)$, then $x(\phi) \in \mathcal{D}$. In fact, fixing $p \in \Gamma$ and using a corollary of the Hahn-Banach Theorem, one proves that $$p(y(\phi)) \leq vol(\text{supp }\phi) \left[\sup_{g \in \text{supp }\phi} |\phi(g)|\right] \left[\sup_{g \in \text{supp }\phi} p(V(g)y)\right], \qquad y \in \mathcal{X}.$$ Choose a net $\left\{x_\alpha\right\}_{\alpha \in \mathcal{A}}$ in $\mathcal{D}$ converging to $x$ in $\mathcal{X}$, whose existence is guaranteed by the denseness hypothesis. Then, it will be proved that $$x_\alpha(\phi) \longrightarrow x(\phi), \qquad \phi \in C_c^\infty(G),$$ with respect to $\tau_\infty$: by hypothesis, $V$ is locally equicontinuous, so making $K := \text{supp }\phi$, there exists $M_{p, K} > 0$ and $q \in \Gamma$ satisfying $$p(V(g)x) \leq M_{p, K} \, q(x), \qquad g \in K, \, x \in \mathcal{X}.$$ Therefore, $$p(x_\alpha(\phi) - x(\phi)) = p((x_\alpha - x)(\phi)) \leq vol(K) \left[\sup_{g \in K} |\phi(g)|\right] M_{p, K} \, q(x_\alpha - x) \longrightarrow 0,$$ for all $\phi \in C_c^\infty(G)$. Also, $$dV(X)(y(\phi)) = \lim_{t \rightarrow 0} \frac{V(\exp tX)(y(\phi)) - y(\phi)}{t}$$ $$= \int_G \lim_{t \rightarrow 0} \frac{\phi([\exp -tX] \, g) - \phi(g)}{t} \, V(g)y \, dg = - y(\tilde{X}(\phi)),$$ for all $y \in \mathcal{X}$, $\phi \in C_c^\infty(G)$ and $X \in \mathfrak{g}$. An iteration of this argument gives $$dV(X_{i_1}) \ldots dV(X_{i_n})(y(\phi)) = (-1)^n y(\phi_n), \qquad n \in \mathbb{N}, \, y \in \mathcal{X},$$ where $X_{i_k} \in \mathfrak{g}$, $1 \leq k \leq n$, and $$\phi_n(g) := \tilde{X}_{i_1} \ldots \tilde{X}_{i_n}(\phi)(g)$$ is a function in $C_c^\infty(G)$. Hence, if $\phi \in C_c^\infty(G)$ is fixed, $$p(dV(X_{i_1}) \ldots dV(X_{i_n})(x_\alpha(\phi) - x(\phi))) = p(dV(X_{i_1}) \ldots dV(X_{i_n})((x_\alpha - x)(\phi)))$$ $$= p((x_\alpha - x)(\phi_n)) \longrightarrow 0$$ for all $p \in \Gamma$ and $n \in \mathbb{N}$. Therefore, it follows that the net $\left\{x_\alpha(\phi)\right\}_{\alpha \in \mathcal{A}}$ in $\mathcal{D}$ actually converges to $x(\phi)$ with respect to the (stronger) projective C$^\infty$ topology. Since $\mathcal{D}$ is complete with respect to this topology, $x(\phi) \in \mathcal{D}$, proving the desired assertion.\\

Now, the main statement of the theorem can be proved: fix $n \in \mathbb{N}$, $$x \in C^n(V) \subseteq \text{Dom }dV(X_{i_1}) \ldots dV(X_{i_n})$$ and a left invariant Haar measure on $G$. Also, fix a sequence $\left\{\phi_m\right\}_{m \in \mathbb{N}}$ of nonnegative functions in $C_c^\infty(G)$ with $$\int_G \phi_m \, dg = 1, \qquad m \in \mathbb{N},$$ and with the property that, for every open neighborhood $U$ of the origin, there exists $m_U \in \mathbb{N}$ such that $\text{supp }\phi_m \subseteq U$, for all $m \geq m_U$. Then, $x(\phi_m)$ converges to $x$ with respect to $\tau$. The fact that the operators $dV(X)$, $X \in \mathfrak{g}$, are continuous with respect to the C$^\infty$-topology ensures that $$dV(X_{i_1}) \ldots dV(X_{i_n}) \, (x(\phi_m) - x) = \int_G \phi_m(g) \, dV(X_{i_1}) \ldots dV(X_{i_n}) \, [V(g)x - x] \, dg,$$ for all $m \in \mathbb{N}$. Also, $$g \longmapsto dV(X_{j_1}) \ldots dV(X_{j_n}) \, V(g)x$$ is a continuous function from $G$ to $\mathcal{X}$, where $X_{j_k}$ is a basis element, so $$p(dV(X_{i_1}) \ldots dV(X_{i_n}) \, (x(\phi_m) - x))$$ $$\leq \left[\int_G \phi_m(g) \, dg\right] \sup_{g \in \text{supp }\phi_m} p(dV(X_{i_1}) \ldots dV(X_{i_n}) \, [V(g)x - x]), \qquad p \in \Gamma,$$ and $p(dV(X_{i_1}) \ldots dV(X_{i_n}) \, (x(\phi_m) - x)) \longrightarrow 0$.\\

Because $x(\phi_m)$ belongs to $\mathcal{D}$, for all $m \in \mathbb{N}$, this proves that $$dV(X_{i_1}) \ldots dV(X_{i_n})|_{C^n(V)} \subset \overline{\left. dV(X_{i_1}) \ldots dV(X_{i_n}) \right|_{\mathcal{D}}},$$ so $$\overline{dV(X_{i_1}) \ldots dV(X_{i_n})|_{C^n(V)}} \subset \overline{\left. dV(X_{i_1}) \ldots dV(X_{i_n}) \right|_{\mathcal{D}}}.$$ Since the other inclusion is immediate, the result is proved.\\

The proofs for a general $L \in \mathfrak{U}(\mathfrak{g})$ and for $C^\infty(V)$ are analogous, and need simple adaptations in the last two paragraphs.\hfill $\blacksquare$\\

\textbf{A corollary of the theorem above is that, under these hypotheses, $\bm{\mathcal{D}}$ is a core for $\bm{\overline{dV(X_k)^n}}$,} for all $n \in \mathbb{N}$, $n \geq 1$, and $1 \leq k \leq d$. To see this, just repeat the proof of Theorem 1.4.1 with $C^n(V)$ replaced by $\text{Dom }dV(X_k)^n$ and $\overline{dV(X_{i_1}) \ldots dV(X_{i_n})|_{C^n(V)}}$ replaced by $\overline{dV(X_k)^n}$. \textbf{In particular, $\bm{C^\infty(V)}$ is a core for $\bm{dV(X_k)}$, for all $\bm{1 \leq k \leq d}$.} Also, under the hypotheses of Theorem 1.4.1, if $\mathcal{D}$ is not assumed to be complete with respect to $\tau_\infty$, then its closure with respect to $\tau_\infty$ will also be left invariant by the operators $V(g)|_{C^\infty(V)}$, because they are all continuous with respect to $\tau_\infty$, as was already noted at the beginning of this subsection. \textbf{Hence, the hypothesis of completeness on $\bm{\mathcal{D}}$ may be dropped, without affecting the resulting conclusions}.\\

The next lemma is an adaptation of \cite[Corollary 3.1.7, page 167]{bratteli1} to $\Gamma$-semigroups of bounded type on locally convex spaces, and will be used in the proof of Theorem 2.13:\\

\textbf{Lemma 1.4.2:} \textit{Let $(\mathcal{X}, \Gamma)$ be a complete Hausdorff locally convex space and $t \longmapsto S(t)$ a $\Gamma$-semigroup of bounded type $w$, whose generator is $T$. Suppose $\mathcal{D} \subseteq \text{Dom }T$ is a dense subspace of $\mathcal{X}$ which is left invariant by $S$ - in other words, $S(t)[\mathcal{D}] \subseteq \mathcal{D}$, for all $t \geq 0$. Then, $\mathcal{D}$ is a core for $T$.}\\

\textbf{Proof of Lemma 1.4.2:} Define $\tilde{T} := \overline{T|_\mathcal{D}}$, so that $\tilde{T} \subset T$ (this holds because $S$ is locally equicontinuous, so $T$ is closed - see \cite[Proposition 1.4]{komura}), and fix a real $\lambda$ satisfying $\lambda > w$. Then, \cite[Theorem 3.3, page 172]{babalola} shows that formula $$\text{R}(\lambda, T)x = \int_0^{+ \infty} e^{-\lambda t} S(t) x \, dt$$ holds for all $x \in \mathcal{X}$, where the integral is strongly convergent with respect to the topology induced by $\Gamma$. Fix $x \in \mathcal{D}$. There exist Riemann sums of the form $$\sum_{k = 1}^N e^{-\lambda t_k} S(t_k)x \, (t_{k + 1} - t_k)$$ for the integral above (which are all elements of $\mathcal{D}$, by the invariance hypothesis) which converge to $\text{R}(\lambda, T)x$ and possess the property that the Riemann sums $$\sum_{k = 1}^N e^{-\lambda t_k} S(t_k) [(\lambda I - T) x] \, (t_{k + 1} - t_k)$$ converge to $\text{R}(\lambda, T)[(\lambda I - T) x] = x$. Closedness of $\tilde{T}$ imply $\text{R}(\lambda, T)x \in \text{Dom }\tilde{T}$ and $(\lambda I - \tilde{T})\text{R}(\lambda, T)x = x$, so $\mathcal{D} \subseteq \text{Ran }(\lambda I - \tilde{T})$. Therefore, $\text{Ran }(\lambda I - \tilde{T})$ is dense in $\mathcal{X}$. Since $\text{R}(\lambda, \tilde{T})$ is continuous and $\tilde{T}$ is closed, $\text{Ran }(\lambda I - \tilde{T})$ is also closed, so $\text{Ran }(\lambda I - \tilde{T}) = \mathcal{X}$.

Now, let $x \in \text{Dom }T$. By what was just proved, there exists $y \in \text{Dom }\tilde{T}$ such that $$(\lambda I - \tilde{T})x = (\lambda I - T)x = (\lambda I - \tilde{T})y,$$ which implies $(\lambda I - \tilde{T})(x - y) = 0$. But the operator $\lambda I - \tilde{T}$, being a restriction of an injective operator, is itself injective, showing that $x = y \in \text{Dom }\tilde{T}$. This proves $\tilde{T} = T$ and ends the proof. \hfill $\blacksquare$

\subsection{Dissipative and Conservative Operators}

\indent

\textbf{Definitions (Dissipative and Conservative Operators) \cite[Definition 3.9]{albanese}:} \textit{Let $\mathcal{X}$ be a Hausdorff locally convex space and $\Gamma$ a fundamental system of seminorms for $\mathcal{X}$. A linear operator $T \colon \text{Dom }T \subseteq \mathcal{X} \longrightarrow \mathcal{X}$ is called \textbf{$\bm{\Gamma}$-dissipative} if, for every $p \in \Gamma$, $\mu > 0$ and $x \in \text{Dom }T$, $$p((\mu I - T)x) \geq \mu \, p(x).$$ Since $\mathcal{X}$ is Hausdorff, this implies $\mu I - T$ is an injective linear operator, for every $\mu > 0$. If both $T$ and $-T$ are $\Gamma$-dissipative, $T$ is called \textbf{$\bm{\Gamma}$-conservative} or, equivalently, if the inequality $$p((\mu I - T)x) \geq |\mu| \, p(x)$$ holds for all $p \in \Gamma$, $\mu \in \mathbb{R}$ and $x \in \text{Dom }T$. Note that every $\Gamma$-conservative operator is $\Gamma$-dissipative and that the definitions of dissipativity and conservativity both depend on the particular choice of the fundamental system of seminorms - see Remark 3.10 of \cite{albanese}, for an illustration of this fact.}\\

\textbf{Definition ($\bm{\Gamma}$-Contractively Equicontinuous Semigroups and $\bm{\Gamma}$-Isometrically Equicontinuous Groups):} \textit{Suppose $\mathcal{X}$ is sequentially complete, $t \longmapsto V(t)$ is an equicontinuous one-parameter semigroup on $\mathcal{X}$ (as in the definition given in Subsection 1.1) with generator $T$ and define, for each $p \in \Gamma$, the seminorm $$\tilde{p}(x) := \sup_{t \geq 0} p(V(t)x), \qquad x \in \mathcal{X}.$$ Then, their very definitions show that $$p(x) \leq \tilde{p}(x) \leq M_p \, q(x) \leq M_p \, \tilde{q}(x), \qquad x \in \mathcal{X},$$ proving that the families $\Gamma$ and $\tilde{\Gamma} := \left\{\tilde{p}: p \in \Gamma\right\}$ are equivalent, in the sense that they generate the same topology of $\mathcal{X}$ (see also \cite[Remark 2.2(i)]{albanesemontel}). It follows also that $$\tilde{p}(V(t)x) \leq \tilde{p}(x), \qquad \tilde{p} \in \tilde{\Gamma}, \, x \in \mathcal{X},$$ which means $t \longmapsto V(t)$ is a \textbf{$\bm{\tilde{\Gamma}}$-contractively equicontinuous semigroup}, according to the terminology introduced at the bottom of page 935 of \cite{albanese}.}\\

\textit{It is a straightforward corollary of \cite[Corollary 1, page 241]{yosida} that $$\mathbb{C} \backslash i \mathbb{R} = \mathbb{C} \, \backslash \left\{\lambda \in \mathbb{C}: \text{Re }\lambda = 0\right\} \subseteq \rho(T).$$ Therefore, formula $$\text{R}(\lambda, T)x = \int_0^{+ \infty} \text{exp}(- \lambda t) \, V(t)x \, dt, \qquad \text{Re }\lambda > 0, \, x \in \mathcal{X}$$ holds (see \cite[Remark 3.12]{albanese}), so it follows that $$\tilde{p}(\text{R}(\lambda, T)x) \leq \frac{1}{\text{Re }\lambda} \, \tilde{p}(x), \qquad \tilde{p} \in \tilde{\Gamma}, \, \text{Re }\lambda > 0, \, x \in \mathcal{X}.$$ In particular, $$\tilde{p}(\text{R}(\lambda, T)x) \leq \frac{1}{\lambda} \, \tilde{p}(x), \qquad \tilde{p} \in \tilde{\Gamma}, \, \lambda > 0, \, x \in \mathcal{X},$$ meaning $T$ is a $\tilde{\Gamma}$-dissipative operator.}\\

\textit{If $t \longmapsto V(t)$ is an equicontinuous one-parameter group these conclusions also follow in perfect analogy by making $$\tilde{p}(x) := \sup_{t \in \mathbb{R}} p(V(t)x), \qquad x \in \mathcal{X},$$ for all $p \in \Gamma$. In this case, $$\tilde{p}(V(t)x) = \tilde{p}(x), \qquad \tilde{p} \in \tilde{\Gamma}, \, x \in \mathcal{X},$$ so $V$ will be called a \textbf{$\bm{\tilde{\Gamma}}$-isometrically equicontinuous group}. Moreover, the semigroups $V_+$ and $V_-$ defined by $V_+(t) := V(t)$ and $V_-(t) := V(-t)$, for all $t \geq 0$, are equicontinuous. Applying the above results to $V_+$ and $V_-$ gives $$\text{R}(\lambda, \pm T)x = \int_0^{+ \infty} \text{exp}(-\lambda t) \, V_{\pm}(t)x \, dt, \qquad x \in \mathcal{X},$$ so $$\tilde{p}(\text{R}(\lambda, \pm T)x) \leq \frac{1}{\text{Re }\lambda} \, \tilde{p}(x), \qquad \tilde{p} \in \tilde{\Gamma}, \, \text{Re }\lambda > 0, \, x \in \mathcal{X}.$$ On the other hand, if $\text{Re }\lambda < 0$, applying the relation $\text{R}(\lambda, T) = -\text{R}(-\lambda, -T)$ yields $$\tilde{p}(\text{R}(\lambda, T)x) = \tilde{p}(\text{R}(-\lambda, -T)x) \leq \frac{1}{-\text{Re }\lambda} \, \tilde{p}(x), \qquad p \in \Gamma, \, x \in \mathcal{X}.$$ Hence, $$\tilde{p}(\text{R}(\lambda, T)x) \leq \frac{1}{|\text{Re }\lambda|} \, \tilde{p}(x), \qquad p \in \Gamma, \, \text{Re }\lambda \neq 0, \, x \in \mathcal{X}.$$ In particular, $$\tilde{p}(\text{R}(\lambda, T)x) \leq \frac{1}{|\lambda|} \, \tilde{p}(x), \qquad p \in \Gamma, \, \lambda \in \mathbb{R} \, \backslash \left\{0\right\}, \, x \in \mathcal{X},$$ proving that $T$ is $\tilde{\Gamma}$-conservative, where $\tilde{\Gamma} := \left\{\tilde{p}: p \in \Gamma\right\}$.}\\

These particular choices of fundamental systems of seminorms will be very useful for the future proofs.

\subsection{The Kernel Invariance Property (KIP),\\ Projective Analytic Vectors}

\indent

Now a concept which will appear very frequently throughout this manuscript is going to be introduced.\\

\textbf{Definition (The Kernel Invariance Property (KIP)):} \textit{If $(\mathcal{X}, \Gamma)$ is a Hausdorff locally convex space, define for each $p \in \Gamma$ the closed subspace $$N_p := \left\{x \in \mathcal{X}: p(x) = 0\right\},$$ often referred to as the \textbf{kernel} of the seminorm $p$, and the quotient map $\pi_p \colon \mathcal{X} \ni x \longmapsto [x]_p \in \mathcal{X}/N_p$. Then, $\mathcal{X}/N_p$ is a normed space with respect to the norm $\|[x]_p\|_p := p(x)$, and is not necessarily complete. Denote its completion by $\mathcal{X}_p := \overline{\mathcal{X}/N_p}$. A densely defined linear operator $T \colon \text{Dom }T \subseteq \mathcal{X} \longrightarrow \mathcal{X}$ is said to possess the \textbf{kernel invariance property (KIP)} with respect to $\Gamma$ if it leaves their seminorms' kernels invariant, that is, $$T \, [\text{Dom }T \cap N_p] \subseteq N_p, \qquad p \in \Gamma.\footnote{\cite{babalola} calls them ``compartmentalized operators''.}$$ If this property is fulfilled, then the linear operators $$T_p \colon \pi_p[\text{Dom }T] \subseteq \mathcal{X}_p \longrightarrow \mathcal{X}_p, \qquad T_p \colon [x]_p \longmapsto [T(x)]_p, \qquad p \in \Gamma, \, x \in \text{Dom }T$$ on the quotients are well-defined, and their domains are dense in each $\mathcal{X}_p$. If a one-parameter semigroup $t \longmapsto V(t)$ is such that all of the operators $V(t)$ possess the (KIP) with respect to $\Gamma$, it will be said that the semigroup $V$ possesses the (KIP) with respect to $\Gamma$. The terminology is analogous for the representation $g \longmapsto V(g)$ of a Lie group.}\\

Note that, if $\Gamma$ is a fundamental system of seminorms with respect to which $t \longmapsto V(t)$ is $\Gamma$-contractively equicontinuous, then $t \longmapsto V(t)$ leaves all the $N_p$'s invariant. Hence, $T$ possesses the kernel invariance property with respect to $\Gamma$, by the definition of infinitesimal generator and by closedness of the $N_p$'s. \textbf{As a corollary of what was proved in Subsection 1.5 and this last observation, it follows that for any generator $\bm{T}$ of a one-parameter equicontinuous semigroup $\bm{t \longmapsto V(t)}$ acting on a sequentially complete Hausdorff locally convex space $\bm{\mathcal{X}}$ there exists a fundamental system of seminorms $\bm{\tilde{\Gamma}}$ for $\bm{\mathcal{X}}$ such that $\bm{T}$ is $\bm{\tilde{\Gamma}}$-dissipative, $\bm{T}$ has the (KIP) with respect to $\bm{\tilde{\Gamma}}$ and $\bm{V}$ is a $\bm{\tilde{\Gamma}}$-contractively equicontinuous semigroup - analogously, this also holds for one-parameter equicontinuous groups.} Also, an obvious fact is that if a linear operator $T$ possesses an extension having the (KIP) with respect to some fundamental system of seminorms, then $T$ also possesses this property (in particular, infinitesimal pregenerators of equicontinuous semigroups always possess the (KIP) with respect to a certain fundamental system of seminorms). In a locally C$^*$-algebra (they will be properly defined later), for example, it is proved inside \cite[Proposition 2]{becker} that an everywhere defined $*$-derivation satisfies the kernel invariance property with respect to \textbf{any} saturated fundamental system of seminorms. Hence, a $*$-derivation $\delta$ which is the pointwise limit of a net of globally defined $*$-derivations - that is, if there exists a net of $*$-derivations $\left\{\delta_j \colon \mathcal{A} \longrightarrow \mathcal{A}\right\}_{j \in J}$ such that $$\delta_j(a) \longrightarrow \delta(a), \qquad a \in \text{Dom }\delta$$ - also possesses this property, since all $N_p$'s are closed $*$-ideals of $\mathcal{A}$ - see also \cite{phillipsd} and \cite{fragoulopoulou2}.\\

Now, a very useful result involving the (KIP), and which will be invoked in Theorem 2.3, is going to be proved:\\

\textbf{Lemma 1.6.1:} \textit{Let $\mathcal{X}$ be a sequentially complete Hausdorff locally convex space and let $T$ be the generator of an equicontinuous semigroup $t \longmapsto V(t)$ on $\mathcal{X}$. If $\Gamma$ is a fundamental system of seminorms for $\mathcal{X}$ with respect to which $T$ has the (KIP), is $\Gamma$-dissipative and $V$ is $\Gamma$-contractively equicontinuous (by what was just observed, it is always possible to arrange such $\Gamma$), then $T_p$ is an infinitesimal pregenerator of a contraction semigroup on the Banach space $\mathcal{X}_p$, for all $p \in \Gamma$.}\\

\textbf{Proof of Lemma 1.6.1:} Since $\mathcal{X}$ is sequentially complete, \cite[Proposition 1.3]{komura} implies that $T$ is a densely defined linear operator. Also, being the generator of an equicontinuous semigroup, it satisfies $$\text{Ran }(\lambda I - T) = \text{Ran }(\lambda I - \overline{T}) = \mathcal{X},$$ for all $\lambda > 0$, by \cite[Corollary 1, page 241]{yosida}. Hence, given $p \in \Gamma$, the induced linear operator $T_p$ on the quotient $\mathcal{X}_p$ is densely defined, dissipative and satisfies $$(\lambda I - T_p)[\pi_p[\text{Dom }T]] = \mathcal{X}/N_p, \qquad \lambda > 0.$$ Therefore, $\overline{\text{Ran }(\lambda I - T_p)} = \overline{\mathcal{X}/N_p} = \mathcal{X}_p$, showing that $\text{Ran }(\lambda I - T_p)$ is dense in $\mathcal{X}_p$, for all $p \in \Gamma$. By Lumer-Phillips Theorem on Banach spaces, $T_p$ is an infinitesimal pregenerator of a contraction semigroup on the Banach space $\mathcal{X}_p$. \hfill $\blacksquare$\\

\textbf{Observation 1.6.1.1:} \textit{If $\mathcal{X}_p = \overline{\mathcal{X}/N_p} = \mathcal{X}/N_p$, for all $p \in \Gamma$, then the stronger conclusion that $T_p$ is the generator, and not only a pregenerator of a contraction semigroup, may be obtained, for all $p \in \Gamma$.}\footnote{This is the case, for example, when $\mathcal{X} \equiv \mathcal{A}$ is a locally C$^*$-algebra - see Section 3.}\\

\textbf{Observation 1.6.1.2:} \textit{A different version of Lemma 1.6.1 may be given, if $\mathcal{X}$ is assumed to be \textbf{complete}. If $T$ is only assumed to be a \textbf{pregenerator} of an equicontinuous group, then the same conclusion can be obtained: since $\overline{T}$ is a densely defined operator, so is $T$. Also, by the Lumer-Phillips Theorem for locally convex spaces \cite[Theorem 3.14]{albanese}, $\text{Ran }(\lambda I - T)$ is a dense subspace of $\mathcal{X}$, so $\text{Ran }(\lambda I - T_p)$ is a dense subspace of $\mathcal{X}_p$, for all $p \in \Gamma$. Hence, the conclusion follows, just as in Lemma 1.6.1.}\\

If $\mathcal{X}$ is a Banach space and $T$ is a linear operator defined on $\mathcal{X}$, then a vector $x \in \mathcal{X}$ is called \textbf{analytic for $\bm{T}$} if $$x \in C^\infty(T) := \bigcap_{n = 1}^{+ \infty} \text{Dom }T^n$$ and there exists $r_x > 0$ such that $$\sum_{n \geq 0} \frac{\|T^n(x)\|}{n!} |u|^n < \infty, \qquad |u| < r_x.$$

The next task will be to define projective analytic vectors on locally convex spaces, so that some useful theorems become available:\\

\textbf{Definition 1.6.2 (Projective Analytic Vectors):} \textit{Let $(\mathcal{X}, \tau)$ be a locally convex space with a fundamental system of seminorms $\Gamma$ and $T$ a linear operator defined on $\mathcal{X}$. An element $x \in \mathcal{X}$ is called a \textbf{$\bm{\tau}$-projective analytic vector} for $T$ if $$x \in C^\infty(T) := \bigcap_{n = 1}^{+ \infty} \text{Dom }T^n$$ and, for every $p \in \Gamma$, there exists $r_{x, p} > 0$ such that $$\sum_{n \geq 0} \frac{p\,(T^n(x))}{n!} |u|^n < \infty, \qquad |u| < r_{x, p}.$$} Note that, for this definition to make sense, it is necessary to show that it does not depend on the choice of the particular system of seminorms: if $x$ is $\tau$-projective analytic with respect to $\Gamma$ and $\Gamma'$ is another saturated family of seminorms generating the topology of $\mathcal{X}$ then, for each $q' \in \Gamma'$, there exists $C_{q'} > 0$ and $q \in \Gamma$ such that $q'(y) \leq C_{q'} \, q(y)$, for all $y \in \mathcal{X}$. Therefore, making $r_{x, q'} := r_{x, q}$, one obtains for every $u \in \mathbb{C}$ satisfying $|u| < r_{x, q'}$ that $$\sum_{n \geq 0} \frac{q'(T^n(x))}{n!} |u|^n \leq C_{q'} \sum_{n \geq 0} \frac{q \,(T^n(x))}{n!} |u|^n < \infty.$$ By symmetry, the assertion is proved. This motivates the use of the notation ``$\tau$-projective analytic'' to indicate that the projective analytic vector in question is related to the topology $\tau$. Sometimes it will be necessary to make explicit which is the topology under consideration to talk about projective analytic vectors, since in some occasions it will be necessary to deal with more than one topology at once - see Theorem 2.9 for a concrete illustration of this situation. When there is no danger of confusion, the symbol $\tau$ will be omitted. The subspace formed by all of the projective analytic vectors for $T$ is going to be denoted by $C^\omega_\leftarrow(T)$. The prefix ``projective'' stands for the fact that $C^\omega_\leftarrow(T)$ can be seen as a dense subspace of the projective limit $\varprojlim \pi_p[C^\omega_\leftarrow(T)]$ via the canonical map $x \longmapsto ([x]_p)_{p \in \Gamma}$ and, if $T$ has the (KIP) with respect to $\Gamma$, then $\pi_p[C^\omega_\leftarrow(T)]$ consists entirely of analytic vectors for $T_p$, for every $p \in \Gamma$ (for the definition of a projective limit of locally convex spaces, see the discussion preceding Theorem 2.13). Note that the projective limit is well-defined, since the family $\left\{\pi_p[C^\omega_\leftarrow(T)]\right\}_{p \in \Gamma}$ gives rise to a canonical projective system.\\

Differently of what is required from the usual definition of analytic vectors, \textbf{no uniformity in $\bm{p}$} is asked in the above definition. Indeed, using the definition of analytic vectors given before, it is possible to adapt the proof of \cite[Theorem 2, page 209]{harish-chandra} and conclude that for every analytic vector $x \in \mathcal{X}$ for $T$ there exists $r_x > 0$ such that, whenever $|u| < r_x$, the series $$\sum_{n = 0}^{+ \infty} \frac{p(T^n (x))}{n!} \, |u|^n$$ is convergent, for every $p \in \Gamma$. Hence, the definition just given is weaker that the usual one, so the subspace of \textbf{$\bm{\tau}$-analytic vectors} $C^\omega(T)$ satisfies $C^\omega(T) \subseteq C^\omega_\leftarrow(T)$.\\

The next theorem, which is a ``locally convex version'' of \cite[Theorem 3]{rusinek}, will play an important role in Theorems 2.9 and 2.12:\\

\textbf{Lemma 1.6.3:} \textit{Let $(\mathcal{X}, \Gamma)$ be a complete Hausdorff locally convex space and let $T$ be a $\Gamma$-conservative linear operator on $\mathcal{X}$ having the (KIP) with respect to $\Gamma$. If $T$ has a dense set of projective analytic vectors, then it is closable and $\overline{T}$ is the generator of a $\Gamma$-isometrically equicontinuous group.}\\

\textbf{Proof of Lemma 1.6.3:} For each $p \in \Gamma$, the densely defined linear operator $T_p$ induced on the quotient $\mathcal{X}/N_p$ possesses a dense subspace $\pi_p[C^\omega_\leftarrow(T)]$ of analytic vectors and is conservative. Therefore, by \cite[Theorem 2]{rusinek}, $T_p$ is an infinitesimal pregenerator of a group of isometries $t \longmapsto V_p(t)$. This implies $\text{Ran }(\lambda I - T_p)$ is dense in $\mathcal{X}_p$ for all $\lambda \in \mathbb{R} \, \backslash \left\{0\right\}$, by Lumer-Phillips Theorem. Fix $\lambda \in \mathbb{R} \, \backslash \left\{0\right\}$. The idea, now, is to prove $\text{Ran }(\lambda I - T)$ is dense in $\mathcal{X}$. So fix $y \in \mathcal{X}$, $\epsilon > 0$, $p \in \Gamma$ and $$V := \left\{x \in \mathcal{X}: p(x - y) < \epsilon\right\}$$ an open neighborhood of $y$. Denseness of $\text{Ran }(\lambda I - T_p)$ in $\mathcal{X}_p$ implies the existence of $x_0 \in \text{Dom }T$ such that $p((\lambda I - T)(x_0) - y) < \epsilon$. Hence, $V \cap \text{Ran }(\lambda I - T) \neq \emptyset$. By the arbitrariness of $p$, it follows that $y$ must belong to $\overline{\text{Ran }(\lambda I - T)}$, so $\overline{\text{Ran }(\lambda I - T)} = \mathcal{X}$ and $\overline{T}$ must be the infinitesimal generator of an equicontinuous group $t \longmapsto V(t)$, by \cite[Proposition 3.13]{albanese} and by the Lumer-Phillips Theorem for locally convex spaces \cite[Theorem 3.14]{albanese}.\\

To see that $\overline{T}$ is actually the generator of a $\Gamma$-isometrically equicontinuous group, first note that formula (7) on \cite[page 248]{yosida} says that $$V(t)x = \lim_{n \rightarrow + \infty} \text{exp}\left(t \, \overline{T} \left(I - \frac{1}{n} \overline{T}\right)^{-1}\right)x, \qquad x \in \mathcal{X}, \, t \in [0, + \infty).$$ Since $\overline{T}$ is $\Gamma$-dissipative, $$p((I - \lambda \overline{T})^{-1}(x)) \leq p(x), \qquad p \in \Gamma, \, \lambda > 0,$$ so the identity $$\overline{T} \left(I - \frac{1}{n} \, \overline{T} \right)^{-1} = n \left(\left(I - \frac{1}{n} \, \overline{T}\right)^{-1} - I\right)$$ shows, for every fixed $t \geq 0$, $p \in \Gamma$ and $x \in \mathcal{X}$ satisfying $p(x) \leq 1$, that $$p (\text{exp}(t \, \overline{T} (I - (1/n) \, \overline{T})^{-1})x) = p (\text{exp}(t \, n (I - (1/n) \, \overline{T})^{-1} - t \, n )x )$$ $$= p (\text{exp}(t \, n \, (I - (1/n) \, \overline{T})^{-1}) \exp (-t \, n)x ) \leq \exp (t \, n) p(x) \, \exp (-t \, n) \leq 1.$$ Hence, $$p(V(t)x) \leq p(x),$$ for every $t \geq 0$, $p \in \Gamma$ and $x \in \mathcal{X}$. But $-\overline{T}$ also generates an equicontinuous semigroup (more precisely, it generates the semigroup $V_- \colon t \longmapsto V(-t)$), so formula $$V(-t)x = V_-(t)x = \lim_{n \rightarrow + \infty} \text{exp}\left(- t \, \overline{T} \left(I + \frac{1}{n} \overline{T}\right)^{-1}\right)x, \qquad x \in \mathcal{X}, \, t \in [0, + \infty),$$ is also valid. Therefore, for every fixed $t \geq 0$, $p \in \Gamma$ and $x \in \mathcal{X}$ satisfying $p(x) \leq 1$, $$p (\text{exp}(-t \, \overline{T} (I + (1/n) \, \overline{T})^{-1})x) = p (\text{exp}(t \, n (I + (1/n) \, \overline{T})^{-1} - t \, n )x )$$ $$= p (\text{exp}(t \, n \, (I + (1/n) \, \overline{T})^{-1}) \exp (-t \, n)x ) \leq \exp (-t \, n) \exp (t \, n \, p(x)) \leq 1,$$ so taking the limit $n \longrightarrow + \infty$ on both sides of this inequality shows $p(V(t)x) \leq p(x)$, for all $t \leq 0$, $p \in \Gamma$ and $x \in \mathcal{X}$. This proves $p(V(t)x) \leq p(x)$, whenever $t \in \mathbb{R}$, $p \in \Gamma$ and $x \in \mathcal{X}$. But this also shows that $p(x) = p(V(-t) V(t)x) \leq p(V(t)x)$, for all $t \in \mathbb{R}$, $p \in \Gamma$ and $x \in \mathcal{X}$, so $V$ is a $\Gamma$-isometrically equicontinuous group. \hfill $\blacksquare$

\subsection{Some Estimates Involving Lie Algebras}

\indent

Let $(\mathcal{X}, \Gamma)$ be a locally convex space, $\mathcal{D} \subseteq \mathcal{X}$ and $\mathcal{L} \subseteq \text{End}(\mathcal{D})$ be a real finite-dimensional Lie algebra of linear operators acting on $\mathcal{X}$, with an ordered basis $\left(B_k\right)_{1 \leq k \leq d}$. Then, for each $1 \leq i, j \leq d$, one has $$(*) \qquad (\text{ad }B_i)(B_j) := [B_i, B_j] = \sum_{k = 1}^d c_{ij}^{(k)} \, B_k,$$ for some constants $c_{ij}^{(k)} \in \mathbb{R}$. Therefore, if $B^\mathsf{u}$ and $B^\mathsf{v}$ are two monomials in the variables $\left(B_k\right)_{1 \leq k \leq d}$, applying recursively the identity above it is possible to prove that the element $(\text{ad }B^\mathsf{u})(B^\mathsf{v}) :=  B^\mathsf{u} B^\mathsf{v} - B^\mathsf{v} B^\mathsf{u}$, which has size $|\mathsf{u}| + |\mathsf{v}|$, is actually a sum of terms of size at most $|\mathsf{u}| + |\mathsf{v}| - 1$. To see this, one must proceed recursively: if $B^\mathsf{u} := B_{\mathsf{u}_1} \ldots B_{\mathsf{u}_{|\mathsf{u}|}}$ and $B^\mathsf{v} := B_{\mathsf{v}_1} \ldots B_{\mathsf{v}_{|\mathsf{v}|}}$, then switching the last element of $B^\mathsf{u}$ with the first element of $B^\mathsf{v}$, and using relation ($*$), gives $$B^\mathsf{u} B^\mathsf{v} = B_{\mathsf{u}_1} \ldots B_{\mathsf{v}_1} B_{\mathsf{u}_{|\mathsf{u}|}} \ldots B_{\mathsf{v}_{|\mathsf{v}|}} + B_{\mathsf{u}_1} \ldots \text{ad }(B_{\mathsf{u}_{|\mathsf{u}|}})(B_{\mathsf{v}_1}) \ldots B_{\mathsf{v}_{|\mathsf{v}|}}$$ $$= B_{\mathsf{u}_1} \ldots B_{\mathsf{v}_1} B_{\mathsf{u}_{|\mathsf{u}|}} \ldots B_{\mathsf{v}_{|\mathsf{v}|}} + \text{ (linear combination of } d \text{ terms of size at most } |\mathsf{u}| + |\mathsf{v}| - 1 \text{)}.$$ Therefore, by iterating this process one concludes that, after repeating this step $|\mathsf{u}| \cdot |\mathsf{v}|$ times, the identity $$(**) \qquad (\text{ad }B^\mathsf{u})(B^\mathsf{v}) :=  B^\mathsf{u} B^\mathsf{v} - B^\mathsf{v} B^\mathsf{u} = \sum_{|\mathsf{w}| \leq |\mathsf{u}| + |\mathsf{v}| + 1} c_\mathsf{w} \, B^\mathsf{w}$$ arises, where the sum has, at most, $d \cdot |\mathsf{u}| \cdot |\mathsf{v}|$ summands. With all of this in mind, one can obtain the estimates $$p((\text{ad }B^\mathsf{u})(B^\mathsf{v})(x)) \leq k \, |\mathsf{u}| \, |\mathsf{v}| \, \rho_{p, |\mathsf{u}| + |\mathsf{v}| - 1}(x), \qquad x \in \mathcal{D}, \, p \in \Gamma,$$ where $k$ is a non-negative constant which does not depend on $\mathsf{u}$ or $\mathsf{v}$, and is related only with the coefficients that come from ($*$) - remind that the seminorm $\rho_{p, |\mathsf{u}| + |\mathsf{v}| - 1}$ was introduced in Subsection 1.4.\\

These estimates will be extremely useful in Section 2.\footnote{See also \cite[page 78]{robinson}.}\\

Lemma 1.7.1 yields a key formula for the proofs of Theorems 2.11 and 2.12. This formula appears in the proof of \cite[Lemma 2.3]{bratteliheat} (under the name of ``Duhamel formula'') and in \cite[page 80]{robinson}, but in both cases, without proof. A product rule-type theorem for locally convex spaces will be necessary, so its statement will be written for the sake of completeness:\\

\textbf{The Product Rule, \cite[Theorem A.1, page 440]{jorgensenmoore}:} \textit{Let $E$, $F$ be Hausdorff locally convex spaces and $I \subseteq \mathbb{R}$ an open interval. Let $K \colon I \longrightarrow \mathcal{L}_s(E, F)$ be a differentiable locally equicontinuous mapping (with $\mathcal{L}_s(E, F)$ being the space of all continuous linear maps from $E$ to $F$ equipped with the strong operator topology) and let $f \colon I \longrightarrow E$ be differentiable. Then, the product mapping $H \colon I \longrightarrow F$ defined by $H(t) := K(t)(f(t))$, for $t \in I$, is differentiable. The first-order derivative is given by $$H'(t) = K(t) f'(t) + K'(t) f(t), \qquad t \in I.$$}

\textbf{Lemma 1.7.1:} \textit{Let $B^\mathsf{u}$ and $B^\mathsf{v}$ be two monomials in the basis elements $\left(B_k\right)_{1 \leq k \leq d}$. Suppose $H_m$ is an element of order $m$ in the complexification $\mathfrak{U}(\mathcal{L})_\mathbb{C}$ of the universal enveloping algebra of $\mathcal{L}$ such that $-H_m$ is a pregenerator of a strongly continuous locally equicontinuous semigroup $t \longmapsto S_t$ satisfying $S_t[\mathcal{X}] \subseteq \mathcal{D}$, for all $t > 0$, and $$(1.7.1.1) \qquad \rho_{p, n}(x) \leq \epsilon^{m - n} \, p(H_m(x)) + \frac{N_p}{\epsilon^n} \, p(x), \qquad p \in \Gamma, \, x \in \mathcal{D},$$ for all $0 < n \leq m - 1$, $0 < \epsilon \leq 1$ and a constant $N_p > 0$ (alternatively, the symbol $S(t)$ will sometimes be used to denote the operator $S_t$). Define the element $(\text{ad }B^\mathsf{u})(H_m)$ of $\mathfrak{U}(\mathcal{L})_\mathbb{C}$ via an extension of the operator $\text{ad }B^\mathsf{u}$ just defined in $(**)$, by linearity. Then, the identity\footnote{In the case where $H_m = - \sum_{k = 1}^d B_k^2$, \cite[page 356]{bratteliheat} refers to this identity as Duhamel formula.} $$\int_0^s B^\mathsf{v} S_r \, [(\text{ad }B^\mathsf{u})(H_m)] \, S_{t - r}(x) \, dr = B^\mathsf{v} [(\text{ad }S_s)(B^\mathsf{u})] \, S_{t - s}(x), \qquad x \in \mathcal{D}, \, 0 \leq s \leq t$$ is valid, for all $t > 0$, where $(\text{ad }S_s)(B^\mathsf{u})$ is defined to be the operator $S_s \, B^\mathsf{u} - B^\mathsf{u} S_s$ on $\mathcal{D}$.}\\

\textbf{Proof of Lemma 1.7.1:} Fix $t > 0$. The first task will be to establish continuity of the function $r \longmapsto B^\mathsf{v} S_r \, B^\mathsf{u} S_{t - r}(x)$ on $[0, t]$, for all $x \in \mathcal{D}$. To this purpose, it will first be obtained the differentiability of $r \longmapsto B^\mathsf{w} S_r(x)$ at $r_0 \in [0, t]$, for all monomials $B^\mathsf{w}$ of size $(q - 1)(m - 1) \leq |\mathsf{w}| \leq q(m - 1)$ in the elements of $\left(B_k\right)_{1 \leq k \leq d}$ and $x \in \mathcal{D}$, a fact which will be proved by induction on $q \geq 1$. To deal with the case $q = 1$, first note that $$\frac{S_r - S_{r_0}}{r - r_0} \, x = S_{r_0} \, \frac{S_{r - r_0} - I}{r - r_0} \, x, \qquad \text{if } r > r_0$$ and $$\frac{S_r - S_{r_0}}{r - r_0} \, x = S_r \, \frac{I - S_{r_0 - r}}{r - r_0} \, x, \qquad \text{if } r < r_0,$$ so the fact that $H_m \, S_s = S_s \, H_m$ on $\mathcal{D}$, for all $s \in [0, + \infty)$, together with (1.7.1.1), gives $$p \left(B^\mathsf{w} S_{r_0} \left(\frac{S_{r - r_0} - I}{r - r_0} \, x + H_m(x) \right) \right) \leq \epsilon^{m - n} \, p \left( S_{r_0} \left(\frac{S_{r - r_0} - I}{r - r_0} \, H_m(x) + H_m^2(x)\right) \right)$$ $$+ \, \frac{N_p}{\epsilon^n} \, p \left(S_{r_0} \left( \frac{S_{r - r_0} - I}{r - r_0} \, x + H_m(x) \right) \right), \qquad x \in \mathcal{D}, \, p \in \Gamma,$$ if $r > r_0$, and $$p \left(B^\mathsf{w} S_r \left( \frac{I - S_{r_0 - r}}{r - r_0} \, x + H_m(x)\right) \right) \leq \epsilon^{m - n} \, p \left( S_r \left( \frac{I - S_{r_0 - r}}{r - r_0} \, H_m(x) + H_m^2(x)\right) \right)$$ $$+ \, \frac{N_p}{\epsilon^n} \, p \left(S_r \left( \frac{I - S_{r_0 - r}}{r - r_0} \, x + H_m(x) \right) \right), \qquad x \in \mathcal{D}, \, p \in \Gamma,$$ if $r < r_0$, for all $0 < \epsilon \leq 1$. Also, $$\frac{S_r - S_{r_0}}{r - r_0} \, x + S_{r_0} \, H_m(x) = S_r \left( \frac{I - S_{r_0 - r}}{r - r_0} \, x + H_m(x)\right)$$ $$+ (S_{r_0} - S_r) \, H_m(x), \qquad r < r_0.$$ Joining all these informations together yields the desired differentiability.

Now, suppose $r \longmapsto B^\mathsf{w} S_r(x)$ is differentiable at $r_0 \in [0, t]$, for all $x \in \mathcal{D}$ and all monomials $B^\mathsf{w}$ in the elements of $\left(B_k\right)_{1 \leq k \leq d}$, with $(q - 1)(m - 1) \leq |\mathsf{w}| \leq q(m - 1)$, for some $q \geq 1$. A fixed monomial $B^\mathsf{w}$ of size $q(m - 1) \leq |\mathsf{w}| := n \leq (q + 1)(m - 1)$ may be decomposed as $B^\mathsf{w} = B^{\mathsf{w}_0} B^{\mathsf{w}'}$, with $|\mathsf{w}_0| = m - 1$ and $|\mathsf{w}'| = n - (m - 1)$. Using $(1.7.1.1)$, one obtains $$p(B^\mathsf{w} S(r)x) = p(B^{\mathsf{w}_0}B^{\mathsf{w}'}S(r)x)$$ $$\leq \epsilon \, p(H_m \, B^{\mathsf{w}'} S(r)x) + \frac{N_p}{\epsilon^{m - 1}} \, p(B^{\mathsf{w}'}S(r)x), \qquad p \in \Gamma, \, x \in \mathcal{D},$$ for all $r \in [0, + \infty)$ and $0 < \epsilon \leq 1$. On the other hand, $$p(H_m \, B^{\mathsf{w}'} S(r)x) \leq p([\text{ad }(H_m)(B^{\mathsf{w}'})] \, S(r)x) + p(B^{\mathsf{w}'} H_m \, S(r)x)$$ $$\leq k \, (n - (m - 1)) \, m \, \rho_{p, n}(S(r)x) + p(B^{\mathsf{w}'} H_m \, S(r)x),$$ $k$ being a constant depending only on $d$, $H_m$ and on the numbers $c_{ij}^{(k)}$ defined by $[B_i, B_j] = \sum_{k = 1}^d c_{ij}^{(k)} B_k$. Hence, $$p(B^\mathsf{w} S(r)x) \leq \epsilon \, [k \, (n - (m - 1)) \, m \, \rho_{p, n}(S(r)x) + p(B^{\mathsf{w}'} H_m \, S(r)x)]$$ $$+ \, \frac{N_p}{\epsilon^{m - 1}} \, p(B^{\mathsf{w}'}S(r)x), \qquad p \in \Gamma, \, x \in \mathcal{D},$$ for all $r \in [0, + \infty)$ and $0 < \epsilon \leq 1$. Choosing an $1 \geq \epsilon_0 > 0$ such that $\epsilon_0 \, k \, (n - (m - 1)) \, m < 1$, and taking the maximum over all monomials $B^\mathsf{w}$ of size $n$, gives $$p(B^\mathsf{w} S(r)x) \leq \rho_{p, n}(S(r)x) \leq \frac{\epsilon_0 \, \rho_{p, n - (m - 1)}(H_m \, S(r)x) + \frac{N_p}{\epsilon^{m - 1}} \, \rho_{p, n - (m - 1)}(S(r)x)}{1 - \epsilon_0 \, k \, (n - (m - 1))},$$ so the induction hypothesis together with a similar argument made in the case $q = 1$ ends the induction proof.\\

Fix $t > 0$, $x \in \mathcal{D}$, $g \in \mathcal{X}'$ and define on $[0, t]$ the functions $$g_1 \colon s \longmapsto \int_0^s g(B^\mathsf{v} S_r [(\text{ad }B^\mathsf{u})(H_m)] S_{t - r}(x)) \, dr$$ and $$g_2 \colon s \longmapsto g(B^\mathsf{v} [(\text{ad }S_s)(B^\mathsf{u})] \, S_{t - s}(x)) = g(B^\mathsf{v} S_s \, B^\mathsf{u} S_{t - s}(x)) - g(B^\mathsf{v} B^\mathsf{u} S_t(x)).$$ An application of \cite[Theorem A.1, page 440]{jorgensenmoore} together with what was just proved gives, in particular, that $r \longmapsto g(B^\mathsf{v} S_r (\text{ad }B^\mathsf{u})(H_m) S_{t - r}(x))$ is continuous on $[0, t]$, so the integral $$\int_0^s g(B^\mathsf{v} S_r \, [(\text{ad }B^\mathsf{u})(H_m)] \, S_{t - r}(x)) \, dr$$ defines a differentiable function on $[0, t]$ with $g_1'(s) = g(B^\mathsf{v} S_s \, [(\text{ad }B^\mathsf{u})(H_m)] \, S_{t - s}(x))$. If $\mathcal{D}$ is equipped with the C$^\infty$ projective topology, then it is clear by the induction proof above that $f \colon s \longmapsto B^\mathsf{u} S_{t - s}(x)$ is a differentiable map from $[0, t]$ to $\mathcal{D}$ and $f'(s) = B^\mathsf{u} S_{t - s} \, H_m(x)$. The same applies for the function $s \longmapsto B^\mathsf{v} S_s(y)$, for every fixed $y \in \mathcal{D}$. Therefore, by the product rule in \cite[Theorem A.1, page 440]{jorgensenmoore}, one sees that $H \colon s \longmapsto B^\mathsf{v} S_s \, B^\mathsf{u} S_{t - s}(x)$ is differentiable on $(0, t)$ and, defining $K \colon s \longmapsto B^\mathsf{v} S_s$, it follows that $$H'(s) = [K(s)(f(s))]' = K'(s)f(s) + K(s)f'(s)$$ $$= -B^\mathsf{v} S_s \, H_m \, B^\mathsf{u} S_{t - s}(x) + B^\mathsf{v} S_s \, B^\mathsf{u} S_{t - s} \, H_m(x) = B^\mathsf{v} S_s \, [(\text{ad }B^\mathsf{u})(H_m)] \, S_{t - s}(x),$$ where the derivative $K'$ is taken with respect to the strong operator topology. This implies $g_1' = g_2'$ on $(0, t)$ and, since $g_1(0) = g_2(0) = 0$, $g_1$ must be equal to $g_2$ on $[0, t]$. By a corollary of the Hahn-Banach Theorem, the equality $$\int_0^s B^\mathsf{v} S_r \, [(\text{ad }B^\mathsf{u})(H_m)] \, S_{t - r}(x) \, dr = B^\mathsf{v} [(\text{ad }S_s)(B^\mathsf{u})] \, S_{t - s}(x), \qquad x \in \mathcal{D}, \, 0 \leq s \leq t,$$ must hold. \hfill $\blacksquare$\\

To close this first section, a lemma which will be invoked in this work (in Theorem 2.2, for example) and which regards extensions of continuous linear maps between locally convex spaces is going to be enunciated (but not proved). The reason for this enunciation is that this theorem is not as well-known as its ``normed counterpart''. Its proof is a straightforward adaptation of the proof of its ``normed version'' (see ``the BLT Theorem'' \cite[Theorem I.7, page 9]{reedsimon1}):\\

\textit{Let $(X, \Gamma_\mathcal{X})$ and $(Y, \Gamma_\mathcal{Y})$ be locally convex spaces, $\mathcal{Y}$ being Hausdorff and complete. If $\mathcal{D} \subseteq \mathcal{X}$ is a linear subspace of $\mathcal{X}$ and $T \colon \mathcal{D} \longrightarrow \mathcal{Y}$ is a continuous linear map then there exists a unique linear map $\tilde{T} \colon \overline{\mathcal{D}} \longrightarrow \mathcal{Y}$ such that $\tilde{T}|_{\mathcal{D}} = T$. In particular, if $\mathcal{D}$ is dense in $\mathcal{X}$ then there exists a unique continuous linear map $\tilde{T} \colon \mathcal{X} \longrightarrow \mathcal{Y}$ such that $\tilde{T}|_{\mathcal{D}} = T$.}

\section{Group Invariance and Exponentiation}

\indent

This section is divided into two big subsections. The first one deals with the problem of constructing a larger, group invariant dense C$^\infty$ domain from a given dense C$^\infty$ domain associated to a finite-dimensional real Lie algebra of linear operators. After giving a more explicit characterization of the maximal C$^\infty$ domain for a finite set of closed linear operators and introducing some technical and very important definitions - like the \textbf{augmented spectrum} and the \textbf{diminished resolvent of a linear operator} and the \textbf{C$\bm{^1}$-closure $\bm{\mathcal{D}_1}$ of a dense C$\bm{^\infty}$ domain $\bm{\mathcal{D}}$} - the construction of the desired group invariant domain begins as an escalade which is divided into five theorems, ending in Theorem 2.5.\\

The other subsection contains the core results of this work: the exponentiability theorems for finite-dimensional real Lie algebras of linear operators on complete Hausdorff locally convex spaces, in which equicontinuity plays a central role.\\

Theorems 2.7 and 2.9 are essential results that will be used to obtain the main exponentiation theorems of this subsection. Among their hypotheses is the imposition that the basis elements must satisfy a kind of \textbf{``joint (KIP)'' condition} with respect to a certain fundamental system of seminorms associated to a generating set (in the sense of Lie algebras) of the operator Lie algebra, which is composed of pregenerators of equicontinuous groups. This is an aesthetically ugly hypothesis, but it will naturally disappear in the statements of the main exponentiation theorems of the paper.

For didactic reasons, the characterization of exponentiable real finite-dimensional Lie algebras of linear operators is broken down into three parts: Theorem 2.11 establishes sufficient conditions for exponentiation in \textbf{Banach spaces}, with an \textbf{arbitrary dense core domain} and a \textbf{general strongly elliptic operator}. Then, Theorem 2.12 reobtains this same result for \textbf{complete Hausdorff locally convex spaces}, and Theorem 2.13 gives necessary conditions for exponentiability in the same context - Observation 2.13.1 slightly strengthens it for \textbf{compact} Lie groups -, but restricted to the case where the dense core domain is the maximal C$^\infty$-domain. Finally, in Theorem 2.14, a characterization of exponentiability in complete Hausdorff locally convex spaces in the same spirit as in \cite[Theorem 3.9]{bratteliheat} is given, but with general strongly elliptic operators being considered.\\

The first step, which is to construct the maximal C$^\infty$ domain for a fixed finite set of closed linear operators on a complete Hausdorff locally convex space (not necessarily Fr\'echet), is based on \cite[Section 7A, page 151]{jorgensenmoore}. The canonical example one should have in mind as an illustration for the construction to follow is the inductive definition, starting from the initial space $C(U)$ of continuous functions on a open set $U \subseteq \mathbb{R}^d$, of the spaces $C^k(U)$, $k \in \mathbb{N}$, of continuously differentiable functions on $U$, and of the space $C^\infty(U)$ of smooth functions on $U$. There, the fixed set of operators is just the set of the usual partial differentiation operators.

\subsection*{Constructing a Group Invariant Domain}
\addcontentsline{toc}{subsection}{Constructing a Group Invariant Domain}

\indent

\textbf{Definition (C$\bm{^\infty}$ Domain for a Set of Linear Operators):} \textit{Let $(\mathcal{X}, \tau)$ be a Hausdorff locally convex space and $\mathcal{S}$ a collection of closable linear operators on $\mathcal{X}$. A subspace $\mathcal{D} \subseteq \mathcal{X}$ satisfying $\mathcal{D} \subseteq \text{Dom }A$ and $A[\mathcal{D}] \subseteq \mathcal{D}$, for all $A \in \mathcal{S}$, will be called a C$^\infty$ domain for $\mathcal{S}$; moreover, the subspace defined by $$\mathcal{X}_\infty(\mathcal{S}) := \bigcup \left\{\mathcal{D} \subseteq \mathcal{X}: \mathcal{D} \text{ is a C$^\infty$ domain for $\mathcal{S}$} \right\}$$ will be called \textbf{the space of C$\bm{^\infty}$ vectors for $\bm{\mathcal{S}}$}. Note that $\mathcal{X}_\infty(\mathcal{S})$ is indeed a vector subspace of $\mathcal{X}$, because if $\mathcal{D}_1$ and $\mathcal{D}_2$ are C$^\infty$ domains for $\mathcal{S}$, then so are $\mathcal{D}_1 + \mathcal{D}_2$ and $\lambda \mathcal{D}_1$, for all $\lambda \in \mathbb{C}$. It is clear that $\mathcal{X}_\infty(\mathcal{S})$ is also a C$^\infty$ domain for $\mathcal{S}$ - it is the maximal C$^\infty$ domain for $\mathcal{S}$ with respect to the partial order defined by inclusion.}\\

It may happen, in some cases, that the subspace $\mathcal{X}_\infty(\mathcal{S})$ reduces to $\left\{0\right\}$, so it is natural to search for sufficient conditions in order to guarantee that this space is ``big enough''. For the objectives of this manuscript, this ``largeness'' will amount to two requirements:

\begin{enumerate}

\item $\mathcal{X}_\infty(\mathcal{S})$ must be dense in $\mathcal{X}$;
\item $\mathcal{X}_\infty(\mathcal{S})$ must be a core for each operator $\overline{A}$, where $A \in \mathcal{S}$.

\end{enumerate}

Suppose, now, that $\mathcal{S} := \left\{B_j\right\}_{1 \leq j \leq d}$ is finite and $\mathcal{D} \subseteq \mathcal{X}$ is a C$^\infty$ domain for $\mathcal{S}$, with $\mathcal{X}$ being a \textbf{complete} Hausdorff locally convex space. It will be shown how a more concrete realization of $\mathcal{X}_\infty(\overline{\mathcal{S}})$, where $$\overline{S} := \left\{\overline{B}: B \in \mathcal{S}\right\} = \left\{\overline{B_j}\right\}_{1 \leq j \leq d},$$ may be obtained, and how to define a natural topology on it. The idea is to proceed via an inductive process to construct a sequence of complete Hausdorff locally convex spaces: let $\Gamma_0$ be a family of seminorms on $\mathcal{X}$ which generates its topology and define $\mathcal{X}_0 := \mathcal{X}$, $\mathcal{X}_1 := \cap_{1 \leq j \leq d} \text{Dom }\overline{B_j}$ and $B_0 := I$. Suppose that there already were constructed $n \geq 1$ Hausdorff locally convex spaces $(\mathcal{X}_k, \Gamma_k)$, $0 \leq k \leq n$, $\Gamma_k := \left\{\rho_{p, k}: p \in \Gamma\right\}$ ($\rho_{p, 0} := p$, for all $p \in \Gamma$), in such a way that:

\begin{enumerate}

\item $\mathcal{X}_{k + 1} \subseteq \mathcal{X}_k \subseteq \text{Dom }\overline{B_j}$, for $1 \leq j \leq d$, $1 \leq k < n$;
\item each $\overline{B_j}$ maps $\mathcal{X}_{k + 1}$ into $\mathcal{X}_k$, for $1 \leq j \leq d$, $0 \leq k < n$;
\item $\Gamma_{k + 1} := \left\{\rho_{p, k + 1}\right\}_{p \in \Gamma}$, with $\rho_{p, k + 1} := \max \left\{\rho_{p, k}(\overline{B_j}(\, \cdot \,)): 0 \leq j \leq d\right\}$, for all $0 \leq k < n$ and all $p \in \Gamma$;\footnote{Here, letting $j$ assume the value 0 is crucial for relating the topologies on $\mathcal{X}_k$ and $\mathcal{X}_{k + 1}$: the topology of the latter space is finer than the topology of the former.}
\item each $\mathcal{X}_k$ is complete with respect to the topology $\tau_k$ defined by the family $\Gamma_k$ of seminorms, for all $0 \leq k \leq n$ ($\tau_0$ is just the initial topology $\tau$ of $\mathcal{X}$).

\end{enumerate}

Define the subspace $\mathcal{X}_{n + 1} := \left\{x \in \mathcal{X}_n: \overline{B_j}(x) \in \mathcal{X}_n, 1 \leq j \leq d\right\}$ of $\mathcal{X}_n$ and a topology $\tau_{n + 1}$ on it via a family of seminorms defined by $\Gamma_{n + 1} := \left\{\rho_{p, n + 1}\right\}_{p \in \Gamma}$, where $$\rho_{p, n + 1} := \max \left\{\rho_{p, n}(\overline{B_j}(\, \cdot \,)): 0 \leq j \leq d\right\},$$ for all $p \in \Gamma$. To see $(\mathcal{X}_{n + 1}, \tau_{n + 1})$ is complete, let $\left\{x_\alpha\right\}_{\alpha \in \mathcal{A}}$ be a $\tau_{n + 1}$-Cauchy net in $\mathcal{X}_{n + 1}$. Then, there exist $x, y \in \mathcal{X}_n$ such that $\left\{x_\alpha\right\}_{\alpha \in \mathcal{A}}$ $\tau_n$-converges to $x$ and $\left\{\overline{B_j}(x_\alpha)\right\}_{\alpha \in \mathcal{A}}$ $\tau_n$-converges to $y$, for all $1 \leq j \leq d$. In particular, $\left\{x_\alpha\right\}_{\alpha \in \mathcal{A}}$ $\tau_0$-converges to $x$ and $\left\{\overline{B_j}(x_\alpha)\right\}_{\alpha \in \mathcal{A}}$ $\tau_0$-converges to $y$, for all $1 \leq j \leq d$, so $(\tau_0 \times \tau_0)$-closedness\footnote{The notation ``$(\tau \times \tau')$'' just used, where $\tau$ is the topology of the domain and $\tau'$ is the topology of the codomain, was employed to emphasize which topologies are being considered to address the question of closability. A similar notation will be used to deal with continuity of the operators. Sometimes, when the topologies under consideration are the same one ($\tau$, for example), the more economic terminology ``$\tau$-continuity'' shall be employed.} of $\overline{B_j}$, $1 \leq j \leq d$, implies $y = \overline{B_j}x$, for all $1 \leq j \leq d$. This proves $x$ belongs to $\left\{z \in \mathcal{X}_n: \overline{B_j}z \in \mathcal{X}_n, 1 \leq j \leq d\right\} = \mathcal{X}_{n + 1}$, so $(\mathcal{X}_{n + 1}, \tau_{n + 1})$ is complete.\\

The next step will be to show that $\mathcal{X}_{n + 1}$ is actually the intersection of the domains of $d$ $(\tau_n \times \tau_n)$-closed operators on $\mathcal{X}_n$. For each $1 \leq j \leq d$, set $B_j^{(n)} := \left(\overline{B_j}\right)|_{\mathcal{X}_{n + 1}}$ (so $\text{Dom }B_j^{(n)} := \mathcal{X}_{n + 1}$) and denote by $\overline{B_j^{(n)}}^{(n)}$ its $(\tau_n \times \tau_n)$-closure.\footnote{Note that it is $(\tau_n \times \tau_n)$-closable.} Since the $\tau_n$-topology on $\mathcal{X}_n$ is finer than the restricted $\tau_0$-topology, $$\text{Dom }\overline{B_j^{(n)}}^{(n)} \subseteq \left\{x \in \mathcal{X}_n: \overline{B_j}x \in \mathcal{X}_n\right\},$$ because $(\mathcal{X}_n, \tau_n)$ is complete. Therefore, $$\bigcap_{j = 1}^d \text{Dom }\overline{B_j^{(n)}}^{(n)} \subseteq \bigcap_{j = 1}^d \left\{x \in \mathcal{X}_n: \overline{B_j}x \in \mathcal{X}_n\right\} = \mathcal{X}_{n + 1}$$ and $$\overline{B_j^{(n)}}^{(n)} \subset \overline{\left(\overline{B_j}\right)|_{\mathcal{X}_{n + 1}}} \subset \overline{B_j}.$$ On the other hand, the definition $\text{Dom }B_j^{(n)} := \mathcal{X}_{n + 1}$, $1 \leq j \leq d$, trivially implies that $\mathcal{X}_{n + 1}$ is contained in $$\bigcap_{j = 1}^d \text{Dom }\overline{B_j^{(n)}}^{(n)}.$$ Hence, the equality $$\bigcap_{j = 1}^d \text{Dom }\overline{B_j^{(n)}}^{(n)} = \mathcal{X}_{n + 1}$$ holds.\\

This concludes the induction step on the construction of a sequence $\left\{(\mathcal{X}_n, \Gamma_n)\right\}_{n \in \mathbb{N}}$ of complete Hausdorff locally convex spaces.\\

Next, consider the subspace $$\mathcal{X}_\infty := \bigcap_{n = 1}^{+ \infty} \mathcal{X}_n$$ equipped with the complete Hausdorff locally convex topology induced by the family of seminorms $$\Gamma_\infty := \left\{\rho_{p, n}\right\}_{p \in \Gamma, n \in \mathbb{N}}.$$ First, note that $\mathcal{X}_\infty$ is a C$^\infty$ domain for $\overline{S}$, because $\mathcal{X}_\infty \subseteq \text{Dom }\overline{B_j}$, for every $1 \leq j \leq d$ and, if $x \in \mathcal{X}_\infty$, then in particular $x$ belongs to $\mathcal{X}_{n + 1}$, for every $n \in \mathbb{N}$, meaning exactly that $\overline{B_j}x \in \mathcal{X}_\infty$, where $1 \leq j \leq d$. As a consequence, $$\mathcal{X}_\infty \subseteq \mathcal{X}_\infty(\overline{\mathcal{S}}).$$ But, if $x \in \mathcal{X}_\infty(\overline{\mathcal{S}})$, then there exists a C$^\infty$ domain $F \subseteq \mathcal{X}$ for $\overline{\mathcal{S}}$ such that $x \in F$. Then, $F \subseteq \text{Dom }\overline{B_j}$ and $\overline{B_j}[F] \subseteq F$, for all $1 \leq j \leq d$, which implies $F \subseteq \mathcal{X}_1$ and $$F \subseteq \left\{y \in \mathcal{X}_1: \overline{B_j}(y) \in \mathcal{X}_1, 1 \leq j \leq d\right\} = \mathcal{X}_2.$$ Using the fact that $\overline{B_j}[F] \subseteq F$, for all $1 \leq j \leq d$, it is seen by the definition of the subspaces $\mathcal{X}_n$ and an inductive argument that $F \subseteq \mathcal{X}_n$, for all $n \in \mathbb{N}$. Thus $$x \in F \subseteq \mathcal{X}_\infty,$$ proving the inclusion $$\mathcal{X}_\infty(\overline{\mathcal{S}}) \subseteq \mathcal{X}_\infty.$$ This proves $$\mathcal{X}_\infty = \mathcal{X}_\infty(\overline{\mathcal{S}}),$$ which gives a concrete realization of $\mathcal{X}_\infty(\overline{\mathcal{S}})$ as a projective limit of complete Hausdorff locally convex spaces whose topology is furnished by the family $\Gamma_\infty$ of seminorms.\\

The next objective will be to show how to obtain, under the hypotheses fixed above, \textbf{together with the additional hypothesis of denseness of $\bm{\mathcal{D}}$ in $\bm{\mathcal{X}}$}, a group invariant dense domain from the fixed dense C$^\infty$ domain $\mathcal{D}$ for $\overline{\mathcal{S}}$. Moreover, it will be proved that this group invariant domain is also a C$^\infty$ domain for $\overline{{\mathcal{S}}}$.\\

In the following five theorems some results coming from the theory of linear operators over finite-dimensional complex vector spaces and from the theory of complex Banach algebras will be used. Therefore, sometimes it will be necessary to embed a real Lie algebra $\mathcal{L}$ into its complexification, $\mathcal{L}_\mathbb{C} := \mathcal{L} + i \mathcal{L}$. Define a basis-dependent norm on $\mathcal{L}$ by $$\left|\sum_{1 \leq k \leq d} b_k B_k\right|_1 := \sum_{1 \leq k \leq d} |b_k|.$$\\ Denoting by $|\, \cdot \,|_{1, \mathbb{C}}$ the norm $\left| \sum_{j = 1}^d b_j B_j + c_j i B_j \right|_{1, \mathbb{C}} := \sum_{j = 1}^d |b_j + i c_j|$ on the complexified Lie algebra $\mathcal{L}_\mathbb{C}$ and by $\|\, \cdot \,\|$ the induced operator norm on $\mathcal{L}(\mathcal{L}_\mathbb{C})$ one sees that $|\, \cdot \,|_{1, \mathbb{C}}$ coincides with $|\, \cdot \,|_1$ on $\mathcal{L}$. Moreover, if $A \in \mathcal{L}$ and $\tilde{\text{ad}}\text{ }A$ denotes the extended linear operator $$\tilde{\text{ad}}\text{ }A \colon \mathcal{L}_\mathbb{C} \longrightarrow \mathcal{L}_\mathbb{C}, \qquad (\tilde{\text{ad}}\text{ }A) (B + iC) := (\text{ad }A)(B) + i (\text{ad }A)(C),$$ then $$\|\tilde{\text{ad}}\text{ }A\| \geq \|(\tilde{\text{ad}}\text{ }A)|_\mathcal{L}\| = \|\text{ad }A\|.$$

Some other definitions which will be of great pertinence are those of \textbf{augmented spectrum} and of \textbf{diminished resolvent} of the closure of a closable operator:\\

\textbf{Definition (Augmented Spectrum and Diminished Resolvent):} \textit{Let $\mathcal{M} \subseteq \text{End}(\mathcal{D})$ be a $d$-dimensional \textbf{complex} Lie algebra. For each closable $T \in \mathcal{M}$, the linear operator on $\mathcal{M}$ given by $\text{ad }T \colon L \longmapsto (\text{ad }T)(L) := [T, L]$ has $1 \leq q \leq d$ eigenvalues. They constitute the spectrum $\sigma_\mathcal{M}(\text{ad }T) := \left\{\mu_k\right\}_{1 \leq k \leq q}$ of the linear operator $\text{ad }T$ when it is considered as a member of the complex (finite-dimensional) Banach algebra $\mathcal{L}(\mathcal{M})$. In this context, the \textbf{augmented spectrum of $\bm{\overline{T}}$} is defined as $$\sigma(\overline{T}; \mathcal{M}) := \sigma(\overline{T}) \cup \left\{\lambda - \mu_k: \lambda \in \sigma(\overline{T}), \mu_k \in \sigma_\mathcal{M}(\text{ad }T)\right\}.$$ Accordingly, the \textbf{diminished resolvent of $\bm{\overline{T}}$} is defined as $$\rho(\overline{T}; \mathcal{M}) := \mathbb{C} \backslash \sigma(\overline{T}; \mathcal{M}).$$ If $\mathcal{M}$ is a \textbf{real} Lie algebra and $T \in \mathcal{M}$, then the definitions of augmented spectrum and of diminished resolvent for $\overline{T}$ become, respectively, $$\sigma(\overline{T}; \mathcal{M}) := \sigma(\overline{T}) \cup \left\{\lambda - \mu_k: \lambda \in \sigma(\overline{T}), \mu_k \in \sigma_{\mathcal{M}_\mathbb{C}}(\tilde{\text{ad}}\text{ }T)\right\}$$ and $$\rho(\overline{T}; \mathcal{M}) := \mathbb{C} \backslash \sigma(\overline{T}; \mathcal{M}).$$}

\textbf{In order to avoid the repetition of words, the following will be assumed as hypotheses in all of the next four auxiliary theorems (Theorems 2.1, 2.2, 2.3 and 2.4):}\\

\textit{Let $(\mathcal{X}, \tau)$ be a complete Hausdorff locally convex space, $\mathcal{D}$ a dense subspace of $\mathcal{X}$ and $\mathcal{L} \subseteq \text{End}(\mathcal{D})$ - remember $\text{End}(\mathcal{D})$ denotes the algebra of all endomorphisms of $\mathcal{D}$ - a finite dimensional \textbf{real} Lie algebra of linear operators, so in particular $\text{Dom }B := \mathcal{D}$, for all $B \in \mathcal{L}$. Fix an ordered basis $\left(B_k\right)_{1 \leq k \leq d}$ of $\mathcal{L}$ and define a basis-dependent norm on $\mathcal{L}$ by $$\left|\sum_{1 \leq k \leq d} b_k B_k\right|_1 := \sum_{1 \leq k \leq d} |b_k|.$$ Suppose $\mathcal{L}$ is generated, as a Lie algebra, by a finite set $\mathcal{S}$ of infinitesimal pregenerators of \textbf{equicontinuous} groups and denote by $t \longmapsto V(t, \overline{A})$ the respective one-parameter group generated by $\overline{A}$, for all $A \in \mathcal{S}$. For each $A \in \mathcal{S}$, choose a fundamental system of seminorms $\Gamma_A$ for $\mathcal{X}$ with respect to which the operator $\overline{A}$ has the (KIP), is $\Gamma_A$-conservative and $V( \, \cdot \, , \overline{A})$ is $\Gamma_A$-isometrically equicontinuous (the arguments at Subsection 1.5 guarantee the existence of such $\Gamma_A$). Define $\Gamma_{A, 1} := \left\{\rho_{p, 1}\right\}_{p \in \Gamma}$, with $\rho_{p, 1} := \max \left\{p(\overline{B_j}(\, \cdot \,)): 0 \leq j \leq d\right\}$, for all $p \in \Gamma$. Now, assume that the following hypotheses are verified:}

\begin{enumerate}

\item \textit{the basis $\left(B_k\right)_{1 \leq k \leq d}$ is formed by closable elements;}
\item \textit{for each $A \in \mathcal{S}$, there exist two complex numbers $\lambda_{(-, A)}$, $\, \lambda_{(+, A)}$ satisfying $\text{Re }\lambda_{(-, A)} < - \|\tilde{\text{ad}}\text{ }A\|$ and $\text{Re }\lambda_{(+, A)} > \|\tilde{\text{ad}}\text{ }A\|$, where $\|\tilde{\text{ad}}\text{ }A\|$ denotes the usual operator norm of $\tilde{\text{ad}}\text{ }A \colon \mathcal{L}_\mathbb{C} \longrightarrow \mathcal{L}_\mathbb{C}$, such that, for all $\lambda \in \mathbb{C}$ satisfying $$|\text{Re }\lambda| \geq M_A := \max \left\{|\text{Re }\lambda_{(-, A)}|, |\text{Re }\lambda_{(+, A)}|\right\},$$ the subspace $\text{Ran }(\lambda I - A)$ is dense in $\mathcal{D}$ with respect to the topology $\tau_1$ induced by the family $\Gamma_{A, 1}$ (note that asking for this denseness, for a fixed $A \in \mathcal{S}$, is equivalent to asking for $\Gamma_{A', 1}$-denseness, for any other $A' \in \mathcal{S}$; this happens because of the fact that $\Gamma_A$ is a fundamental system of seminorms for $\mathcal{X}$ - for a more careful argument, see the next paragraph).}

\end{enumerate}

As the basis $\left(B_k\right)_{1 \leq k \leq d}$ is formed by closable elements, the chain of inclusions $$\mathcal{D} \subseteq \bigcap_{k = 1}^d \text{Dom }\overline{B_k} \subseteq \mathcal{X}$$ is clear. Define $\mathcal{D}_1$ as the closure of $\mathcal{D}$ in $\cap_{k = 1}^d \text{Dom }\overline{B_k}$ with respect to the $\tau_1$-topology. Note that the topology of $\mathcal{D}_1$ is independent of the choice of the element $A \in \mathcal{S}$: let $A_1, A_2 \in \mathcal{S}$ and $\Gamma_{A_1}$, $\Gamma_{A_2}$ be their respective families of seminorms, with $\tau_{A_1, 1}$ and $\tau_{A_2, 1}$ being the corresponding generated topologies. If $x \in \mathcal{D}_1$ is a fixed element in the $\tau_{A_1, 1}$-closure of $\mathcal{D}$ and $\left\{x_\alpha\right\}$ is a net in $\mathcal{D}$ such that $\rho_{p, 1}(x_\alpha - x) \rightarrow 0$, for all $\rho_{p, 1} \in \Gamma_{A_1, 1}$, then $$p(B_k \, (x_\alpha - x)) \longrightarrow 0, \qquad p \in \Gamma_{A_1}, \, 0 \leq k \leq d.$$ Since $\Gamma_{A_1}$ is a fundamental system of seminorms for $\mathcal{X}$, this is equivalent to the convergence $B_k \, (x_\alpha - x) \rightarrow 0$ in $\mathcal{X}$, for each $0 \leq k \leq d$. But, again, $\Gamma_{A_2}$ is also a fundamental system of seminorms for $\mathcal{X}$, so $$p'(B_k \, (x_\alpha - x)) \longrightarrow 0, \qquad p' \in \Gamma_{A_2}, \, 0 \leq k \leq d,$$ must hold. But this implies $x$ belongs to the $\tau_{A_2, 1}$-closure of $\mathcal{D}$, showing that the $\tau_{A_1, 1}$-closure of $\mathcal{D}$ is contained in the $\tau_{A_2, 1}$-closure. By symmetry, this implies that they are equal, showing the asserted independence.\\

\textbf{Definition (C$\bm{^1}$-closure of $\bm{\mathcal{D}}$):} \textit{For the rest of this section the topology on $\mathcal{D}_1$ generated by the family $\Gamma_{A, 1}$, for some $A \in \mathcal{S}$, will be denoted by $\bm{\tau_1}$.}\\

Note, also, the inclusions $$\mathcal{D} \subseteq \mathcal{D}_1 \subseteq \bigcap_{1 \leq k \leq d} \text{Dom }\overline{B_k}.$$

The next theorem is a straightforward adaptation of \cite[Theorem 5.1, page 112]{jorgensenmoore}.\\

\textbf{Theorem 2.1 (Commutation Relations Involving the Resolvent Operators):} \textit{Let $A \in \mathcal{S}$ with $\sigma_{\mathcal{L}_\mathbb{C}}(\tilde{\text{ad}}\text{ }A) := \left\{\mu_j\right\}_{1 \leq j \leq q}$. If $\lambda \in \rho(\overline{A}; \mathcal{L})$ and $\mu \in \mathbb{C}$ is a number such that $\lambda + \mu \in \rho(\overline{A})$, then $$(2.1.1) \qquad B \, \text{R}(\lambda, \overline{A})(x) = \sum_{k = 0}^n (-1)^k \text{R}(\lambda + \mu, \overline{A})^{k + 1}(\text{ad }A - \mu)^k(B)(x)$$ $$+ \, (-1)^{n + 1} \text{R}(\lambda + \mu, \overline{A})^{n + 1} (\text{ad }A - \mu)^{n + 1}(B) \, \text{R}(\lambda, \overline{A})(x),$$ where $x \in \text{Ran }(\lambda I - A)$, $B \in \mathcal{L}$ and $n \in \mathbb{N}$. Also, $$(2.1.2) \qquad B \, \text{R}(\lambda, \overline{A})(x)$$ $$= \sum_{1 \leq j \leq q, \, 0 \leq k \leq s_j} (-1)^k \, \text{R}(\lambda + \mu_j, \overline{A})^{k + 1} \, (\text{ad }A - \mu_j I)^k \, (P_j B)(x), \qquad x \in \text{Ran }(\lambda I - A), \, B \in \mathcal{L},$$ with $P_j$ being the projection over the generalized eigenspace $$(\mathcal{L}_\mathbb{C})_j := \left\{B \in \mathcal{L}_\mathbb{C}: (\tilde{\text{ad}}\text{ }A - \mu_j I)^s(B) = 0, \text{ for some } s\right\}$$ and $s_j$ being the non-negative integer satisfying $(\tilde{\text{ad}}\text{ }A - \mu_j I)^{s_j} \neq 0 = (\tilde{\text{ad}}\text{ }A - \mu_j I)^{s_j + 1}$ on $(\mathcal{L}_\mathbb{C})_j$.\footnote{The existence of the projections $P_j$ and the integers $s_j$ is guaranteed by the Primary Decomposition Theorem - see \cite[Theorem 12, page 220]{hoffman}.}}\\

\textbf{Proof of Theorem 2.1:} To prove (2.1.1), first note that if $x = (\lambda I - A)(y)$ and $B \in \mathcal{L}$, $$\text{R}(\lambda + \mu, \overline{A}) \, B(x) = \text{R}(\lambda + \mu, \overline{A}) \, B \, (\lambda I - A)(y)$$ $$= \text{R}(\lambda + \mu, \overline{A}) \, [(\lambda + \mu - A)B + [A, B] - \mu B](y) = B(y) + \text{R}(\lambda + \mu, \overline{A}) \, (\text{ad }A - \mu)(B)(y)$$ $$= B \, \text{R}(\lambda, \overline{A})(x) + \text{R}(\lambda + \mu, \overline{A}) \, (\text{ad }A - \mu)(B) \, \text{R}(\lambda, \overline{A})(x).$$ Now, suppose that for a fixed $n \in \mathbb{N}$, the equality $$C \, \text{R}(\lambda, \overline{A})(x) = \sum_{k = 0}^n (-1)^k \text{R}(\lambda + \mu, \overline{A})^{k + 1}(\text{ad }A - \mu)^k(C)(x)$$ $$+ \, (-1)^{n + 1} \text{R}(\lambda + \mu, \overline{A})^{n + 1} (\text{ad }A - \mu)^{n + 1}(C) \, \text{R}(\lambda, \overline{A})(x), \qquad x \in \text{Ran }(\lambda I - A),$$ holds, for all $C \in \mathcal{L}$. Substituting $C$ by $(\text{ad }A - \mu)(B)$, with $B \in \mathcal{L}$, one obtains $$(\text{ad }A - \mu)(B) \, \text{R}(\lambda, \overline{A})(x) = \sum_{k = 0}^n (-1)^k \text{R}(\lambda + \mu, \overline{A})^{k + 1}(\text{ad }A - \mu)^{k + 1}(B)(x)$$ $$+ \, (-1)^{n + 1} \text{R}(\lambda + \mu, \overline{A})^{n + 1} (\text{ad }A - \mu)^{n + 2}(B) \, \text{R}(\lambda, \overline{A})(x), \qquad x \in \text{Ran }(\lambda I - A).$$ Applying the operator $\text{R}(\lambda + \mu, \overline{A})$ on both members and using the first equality obtained in the proof gives $$\text{R}(\lambda + \mu, \overline{A}) \, B(x) - B \, \text{R}(\lambda, \overline{A})(x) = \text{R}(\lambda + \mu, \overline{A}) \, (\text{ad }A - \mu)(B) \, \text{R}(\lambda, \overline{A})(x)$$ $$= \sum_{k = 0}^n (-1)^k \text{R}(\lambda + \mu, \overline{A})^{k + 2}(\text{ad }A - \mu)^{k + 1}(B)(x)$$ $$+ \, (-1)^{n + 1} \text{R}(\lambda + \mu, \overline{A})^{n + 2} (\text{ad }A - \mu)^{n + 2}(B) \, \text{R}(\lambda, \overline{A})(x), \qquad x \in \text{Ran }(\lambda I - A), \, B \in \mathcal{L}.$$ This yields $$B \, \text{R}(\lambda, \overline{A})(x) = \sum_{k = 0}^{n + 1} (-1)^k \text{R}(\lambda + \mu, \overline{A})^{k + 1}(\text{ad }A - \mu)^k(B)(x)$$ $$+ \, (-1)^{n + 2} \text{R}(\lambda + \mu, \overline{A})^{n + 2} (\text{ad }A - \mu)^{n + 2} \, \text{R}(\lambda, \overline{A})(x),$$ for all $x \in \text{Ran }(\lambda I - A)$ and $B \in \mathcal{L}$, which finishes the induction proof and establishes (2.1.1) for all $n \in \mathbb{N}$.\\

In order to prove (2.1.2), first note that $\lambda + \mu$ belongs to $\rho(\overline{A})$ whenever $\mu$ belongs to $\sigma_{\mathcal{L}_\mathbb{C}}(\tilde{\text{ad}}\text{ }A)$. Indeed, if $\lambda + \mu$ belonged to $\sigma(\overline{A})$, then $\lambda = (\lambda + \mu) - \mu \in \sigma(\overline{A}; \mathcal{L}_\mathbb{C})$, which is absurd. Therefore, if $\mu \in \sigma_{\mathcal{L}_\mathbb{C}}(\tilde{\text{ad}}\text{ }A)$, (2.1.1) may be applied for $\lambda + \mu$. Substituting $B \in \mathcal{L}$ by $P_j B \in \mathcal{L}$,\footnote{$P_j B \in \mathcal{L}$ follows from the fact that $P_j$ is a polynomial in $\text{ad }A$.} $\mu$ by $\mu_j \in \sigma_{\mathcal{L}_\mathbb{C}}(\tilde{\text{ad}}\text{ }A)$ and $n$ by $s_j$, (2.1.1) gives $$P_j B \, \text{R}(\lambda, \overline{A})(x) = \sum_{k = 0}^{s_j} (-1)^k \text{R}(\lambda + \mu_j, \overline{A})^{k + 1}(\text{ad }A - \mu_j)^k(P_j B)(x)$$ $$+ \, (-1)^{s_j + 1} \text{R}(\lambda + \mu_j, \overline{A})^{s_j + 1} (\text{ad }A - \mu_j)^{s_j + 1}(P_j B) \, \text{R}(\lambda, \overline{A})(x)$$ $$= \sum_{k = 0}^{s_j} (-1)^k \text{R}(\lambda + \mu_j, \overline{A})^{k + 1}(\text{ad }A - \mu_j)^k(P_j B)(x), \qquad x \in \text{Ran }(\lambda I - A), \, B \in \mathcal{L},$$ for all $\mu_j \in \sigma_{\mathcal{L}_\mathbb{C}}(\tilde{\text{ad}}\text{ }A)$. Since $\sum_{1 \leq j \leq q} P_j B = B$, summation of the last equation over all $1 \leq j \leq q$ establishes the desired result. \hfill $\blacksquare$\\

\textbf{From Theorem 2.2 until Theorem 2.4, consider fixed an $\bm{A \in \mathcal{S}}$ with $$\bm{\sigma_{\mathcal{L}_\mathbb{\textbf{C}}}(\tilde{\text{ad}}\text{ }A) := \left\{\mu_j\right\}_{1 \leq j \leq q}}.$$ It is also important to have in mind that the results in all of the next three theorems - and even in Theorem 2.1 - remain valid if one substitutes $\bm{A}$ by $\bm{-A}$. This observation will be particularly useful in Theorem 2.4.}\\

Theorem 2.2 below is based on \cite[Theorem 5.4(1), page 119]{jorgensenmoore}, but here it is assumed that $A$ is not just closable, but that it is also $\Gamma_A$-conservative. This conservativity hypothesis (along with the (KIP)) will be, just as in Theorems 2.3 and 2.4, a key hypothesis to circumvent the fact that the resolvent operator $\text{R}(\lambda, \overline{A})$ does not belong, in general, to a Banach algebra of operators on $\mathcal{X}$. This is a technical obstruction which makes it difficult to extend the results beyond the normed case, as noted in \cite[page 113]{jorgensenmoore}. In the more general setting of complete Hausdorff locally convex spaces, one has the following:\\

\textbf{Theorem 2.2 (C$\bm{^1}$-Continuity of the Resolvent Operators):} \textit{Let $\lambda \in \mathbb{C}$ satisfy $|\text{Re }\lambda| \geq M_A := \max \left\{|\text{Re }\lambda_{(-, A)}|, |\text{Re }\lambda_{(+, A)}|\right\} > 0$ so that, in particular, $\lambda \in \rho(\overline{A})$. Then, $\text{R}(\lambda, \overline{A})$ leaves $\mathcal{D}_1$ invariant and restricts there as a $\tau_1$-continuous linear operator $\text{R}_1(\lambda, \overline{A})$. Moreover, $\text{R}_1(\lambda, \overline{A}) = (\lambda I - \overline{A}^{\tau_1})^{-1} = \text{R}(\lambda, \overline{A}^{\tau_1})$, where $\overline{A}^{\tau_1}$ denotes the $(\tau_1 \times \tau_1)$-closure of the operator $A \colon \mathcal{D} \subseteq \mathcal{D}_1 \longrightarrow \mathcal{D}_1$ when seen as a densely defined linear operator acting on $\mathcal{D}_1$.}\\

\textbf{Proof of Theorem 2.2:} First note that the condition $$|\text{Re }\lambda| \geq M_A := \max \left\{|\text{Re }\lambda_{(-, A)}|, |\text{Re }\lambda_{(+, A)}|\right\}$$ implies $\lambda \in \rho(\overline{A}; \mathcal{L})$, so Theorem 2.1 is applicable. By equation (2.1.2) of Theorem 2.1, $$B_k \, \text{R}(\lambda, \overline{A})(x)$$ $$= \sum_{1 \leq j \leq q, \, 0 \leq i \leq s_j} (-1)^i \, \text{R}(\lambda + \mu_j, \overline{A})^{i + 1} \, (\text{ad }A - \mu_j I)^i \, (P_j B_k)(x),$$ for every $1 \leq k \leq d$ and $x \in \text{Ran }(\lambda I - A)$, with $P_j$ being the projection over the generalized eigenspace $$(\mathcal{L}_\mathbb{C})_j := \left\{B \in \mathcal{L}_\mathbb{C}: (\tilde{\text{ad}}\text{ }A - \mu_j I)^s(B) = 0, \text{ for some } s\right\}$$ and $s_j$ being the non-negative integer satisfying $(\tilde{\text{ad}}\text{ }A - \mu_j I)^{s_j} \neq 0 = (\tilde{\text{ad}}\text{ }A - \mu_j I)^{s_j + 1}$ on $(\mathcal{L}_\mathbb{C})_j$. The restriction $|\text{Re }\lambda| > \|\tilde{\text{ad}}\text{ }A\|_\mathbb{C}$ of the hypothesis implies $\lambda + \mu_j \neq 0$, for all $1 \leq j \leq q$, because $|\mu_j| \leq \|\tilde{\text{ad}}\text{ }A\|_\mathbb{C}$, so $\text{Re }\lambda + \text{Re }\mu_j \neq 0$. Since $\overline{A}$ is a $\Gamma_A$-conservative linear operator, it follows that $$p(((\lambda + \mu_j) I - \overline{A})(x)) \geq |\lambda + \mu_j| \, p(x), \qquad p \in \Gamma_A, \, x \in \text{Dom }\overline{A}.$$ Therefore, in particular, $$p(\text{R}(\lambda + \mu_j, \overline{A})(x)) \leq \frac{1}{|\lambda + \mu_j|} \, p(x), \qquad p \in \Gamma_A, \, x \in \mathcal{X}.$$ Hence, $\text{R}(\lambda + \mu_j, \overline{A})$ possesses the (KIP) with respect to $\Gamma_A$, for all $1 \leq j \leq q$. This allows one to define for each $p \in \Gamma_A$ and each $1 \leq j \leq q$ the unique everywhere-defined continuous extension of the densely defined $\| \, \cdot \, \|_p$-continuous linear operator $(\text{R}(\lambda + \mu_j, \overline{A}))_p$ induced by $\text{R}(\lambda + \mu_j, \overline{A})$ on the completion $\mathcal{X}_p$. Such extension will be denoted by $\text{R}_p(\lambda + \mu_j, \overline{A})$, for every $p \in \Gamma_A$ and $1 \leq j \leq q$. Analogously, the $\| \, \cdot \, \|_p$-continuous linear operator $(\text{R}(\lambda, \overline{A}))_p$ also has such an extension $\text{R}_p(\lambda, \overline{A})$ on $\mathcal{X}_p$, by the $\Gamma_A$-conservativity hypothesis.\\

For each $1 \leq k \leq d$, $1 \leq j \leq q$ and $0 \leq i \leq s_j$ write $$(\text{ad }A - \mu_j I)^i \, (P_j B_k) = \sum_{l = 1}^d c_l B_l, \qquad c_l \in \mathbb{R}.\footnote{This is possible because $A$ and $P_j B_k$ both belong to $\mathcal{L}$.}$$ Then, for all $x \in \mathcal{D}$ the estimates $$p((\text{ad }A - \mu_j I)^i \, (P_j B_k)(x)) \leq \sum_{l = 1}^d |c_l| \, p(B_l(x))$$ $$\leq |(\text{ad }A - \mu_j I)^i \, (P_j B_k)|_1 \, \max \left\{p(B_k(x)): 0 \leq k \leq d\right\}$$ $$\leq \|\tilde{\text{ad}}\text{ }A - \mu_j I\|^i \, |P_j B_k|_1 \, \rho_{p, 1}(x), \qquad p \in \Gamma_A,$$ are legitimate. Therefore, denoting by $\|\, \cdot \,\|_p$ the usual operator norm on $\mathcal{L}(\mathcal{X}_p)$, the inequality $$p(B_k \, \text{R}(\lambda, \overline{A})(x)) = p\left(\sum_{1 \leq j \leq q, \, 0 \leq i \leq s_j} (-1)^i \, \text{R}(\lambda + \mu_j, \overline{A})^{i + 1} \, (\text{ad }A - \mu_j I)^i \, (P_j B_k)(x)\right)$$ $$\leq \sum_{1 \leq j \leq q, \, 0 \leq i \leq s_j} p(\text{R}_p(\lambda + \mu_j, \overline{A})^{i + 1}[(\text{ad }A - \mu_j I)^i \, (P_j B_k)(x)]_p)$$ $$\leq \sum_{1 \leq j \leq q, \, 0 \leq i \leq s_j} \|\text{R}_p(\lambda + \mu_j, \overline{A})\|_p^{i + 1} p((\text{ad }A - \mu_j I)^i \, (P_j B_k)(x))$$ $$\leq \sum_{1 \leq j \leq q, \, 0 \leq i \leq s_j} \|\text{R}_p(\lambda + \mu_j, \overline{A})\|_p^{i + 1} \|\tilde{\text{ad}}\text{ }A - \mu_j I\|^i \, |P_j B_k|_1 \, \rho_{p, 1}(x)$$ $$\leq C_p(k) \, \rho_{p, 1}(x), \qquad p \in \Gamma_A, \, 1 \leq k \leq d, \, x \in \text{Ran }(\lambda I - A)$$ is obtained, for some $C_p(k) > 0$. For $k = 0$ the estimate with $B_0 = I$ is simpler: $$p(B_0 \, \text{R}(\lambda, \overline{A})(x)) \leq \|\text{R}_p(\lambda, \overline{A})\|_p \, p(x) \leq \|\text{R}_p(\lambda, \overline{A})\|_p \, \rho_{p, 1}(x), \qquad p \in \Gamma_A,$$ where $\text{R}_p(\lambda, \overline{A})$ is defined in an analogous way as the operators $\text{R}_p(\lambda + \mu_j, \overline{A})$ were. This implies, upon taking the maximum over $0 \leq k \leq d$, that $$\rho_{p, 1}(\text{R}(\lambda, \overline{A})(x))$$ $$\leq \max \left\{\|\text{R}_p(\lambda, \overline{A})\|_p, \max \left\{C_p(k): 1 \leq k \leq d\right\} \right\} \rho_{p, 1}(x), \qquad p \in \Gamma_A, \, x \in \text{Ran }(\lambda I - A),$$ showing that $\text{R}(\lambda, \overline{A}) \colon \text{Ran }(\lambda I - A) \longrightarrow \mathcal{D} \subseteq \mathcal{D}_1$ is a $\tau_1$-continuous linear operator on $\mathcal{D}_1$. Hence, one might extend it via limits to a linear $\tau_1$-continuous everywhere defined linear operator $\text{R}_1(\lambda, \overline{A})$ on $\mathcal{D}_1$, as a consequence of the $\tau_1$-denseness hypothesis of $\text{Ran }(\lambda I - A)$ in $\mathcal{D}$. Since $\tau_1$ is finer than $\tau$ and $\text{R}(\lambda, \overline{A})$ is $\tau$-continuous, the equality $$\text{R}_1(\lambda, \overline{A}) = \left. \text{R}(\lambda, \overline{A})\right|_{\mathcal{D}_1}$$ is true and, moreover, $\text{R}(\lambda, \overline{A})$ leaves $\mathcal{D}_1$ invariant (to obtain this conclusion, it was also used that $\mathcal{D}_1$ is Hausdorff with respect to its $\tau_1$-topology).\\

To prove $\text{R}_1(\lambda, \overline{A}) = (\lambda I - \overline{A}^{\tau_1})^{-1}$ first note that, in fact, $A$ is a $(\tau_1 \times \tau_1)$-closable linear operator: indeed, if $\left\{x_\alpha\right\}$ is a net in $\mathcal{D}$ which is $\tau_1$-convergent to 0 and $A(x_\alpha) \xrightarrow{\alpha} y \in \mathcal{D}_1$, then both nets $\left\{x_\alpha\right\}$ and $\left\{A(x_\alpha)\right\}$ are $\tau$-convergent to 0 and $y$, respectively, since the $\tau_1$-topology is stronger than the $\tau$-topology. Since $A$ is $(\tau \times \tau)$-closable, it follows that $y = 0$, which proves $A$ is $(\tau_1 \times \tau_1)$-closable. Hence, its $(\tau_1 \times \tau_1)$-closure, $\overline{A}^{\tau_1}$, is well-defined, and $\overline{A}^{\tau_1} \subseteq \overline{A}$.

Now, fix $x \in \text{Dom }\overline{A}^{\tau_1}$ and take a net $\left\{x_\alpha\right\}$ in $\mathcal{D}$ such that $x_\alpha \longrightarrow x$ and $A(x_\alpha) \longrightarrow \overline{A}^{\tau_1}(x)$, both convergences being in the $\tau_1$-topology. This implies $$\text{R}_1(\lambda, \overline{A}) \left[(\lambda I - \overline{A}^{\tau_1})(x_\alpha)\right] = \text{R}(\lambda, \overline{A}) [(\lambda I - A)(x_\alpha)] = x_\alpha \longrightarrow x$$ in the $\tau_1$-topology. But $$(\lambda I - \overline{A}^{\tau_1})(x_\alpha) = (\lambda I - A)(x_\alpha) \longrightarrow (\lambda I - \overline{A}^{\tau_1})(x)$$ relatively to the $\tau_1$-topology, which implies $$\text{R}_1(\lambda, \overline{A}) \left[(\lambda I - \overline{A}^{\tau_1})(x_\alpha)\right] \longrightarrow \text{R}_1(\lambda, \overline{A}) \left[(\lambda I - \overline{A}^{\tau_1})(x)\right]$$ in the $\tau_1$-sense, by $\tau_1$-continuity of the linear operator $\text{R}_1(\lambda, \overline{A})$. By the uniqueness of the limit, $$\text{R}_1(\lambda, \overline{A}) \left[(\lambda I - \overline{A}^{\tau_1})(x)\right] = x.$$ This proves that $\text{R}_1(\lambda, \overline{A})$ is a left inverse for $(\lambda I - \overline{A}^{\tau_1})$. To prove that $\text{R}_1(\lambda, \overline{A})$ is a right inverse, fix an arbitrary $x \in \mathcal{D}_1$ and a net $\left\{x_\alpha\right\}$ in $\mathcal{D}$ such that $x_\alpha \longrightarrow x$ in the $\tau_1$-topology. Then, $$(\lambda I - \overline{A}^{\tau_1}) \text{R}_1(\lambda, \overline{A}) (x_\alpha) = (\lambda I - A) \text{R}(\lambda, \overline{A}) (x_\alpha) = x_\alpha \longrightarrow x,$$ in the $\tau_1$-topology. On the other hand, $\tau_1$-continuity of $\text{R}_1(\lambda, \overline{A})$ and $\tau_1$-closedness of $\overline{A}^{\tau_1}$ imply $$(\lambda I - \overline{A}^{\tau_1}) \text{R}_1(\lambda, \overline{A}) (x_\alpha) \longrightarrow (\lambda I - \overline{A}^{\tau_1}) \text{R}_1(\lambda, \overline{A}) (x).$$ As before, uniqueness of the limit establishes that $\text{R}_1(\lambda, \overline{A})$ is a right inverse for $\lambda I - \overline{A}^{\tau_1}$. This proves $\text{R}_1(\lambda, \overline{A})$ is a two-sided inverse for the operator $\lambda I - \overline{A}^{\tau_1}$ on $\mathcal{D}_1$ and that $\text{R}_1(\lambda, \overline{A}) = \text{R}(\lambda, \overline{A}^{\tau_1})$, as claimed. \hfill $\blacksquare$\\

The proof of the next theorem is strongly inspired in that of \cite[Theorem 5.2, page 114]{jorgensenmoore}, and gives a version of formula (1) of the reference in the context of locally convex spaces, when $A$ is assumed to be a pregenerator of an equicontinuous group.\\

\textbf{Theorem 2.3 (Series Expansion of Commutation Relations):} \textit{Let $\lambda \in \mathbb{C}$ satisfy $$|\text{Re }\lambda| \geq M_A := \max \left\{|\text{Re }\lambda_{(-, A)}|, |\text{Re }\lambda_{(+, A)}|\right\} > 0.$$ For all $x \in \mathcal{D}$ and all closable $B \in \mathcal{L}$, it is verified that $R(\lambda, \overline{A})(x) \in \text{Dom }\overline{B}$ and $$(2.3.1) \qquad \overline{B} \, \text{R}(\lambda, \overline{A})(x) = \sum_{k = 0}^{+ \infty} (-1)^k \text{R}(\lambda, \overline{A})^{k + 1} [(\text{ad }A)^k(B)](x).$$}

\textbf{Proof of Theorem 2.3:} By (2.1.1) of Theorem 2.1, $$B \, \text{R}(\lambda, \overline{A})(x) = \left[ \sum_{k = 0}^n (-1)^k \text{R}(\lambda, \overline{A})^{k + 1} (\text{ad }A)^k(B)(x) \right]$$ $$+ \, (-1)^{n + 1} \text{R}(\lambda, \overline{A})^{n + 1} (\text{ad }A)^{n + 1}(B) \, \text{R}(\lambda, \overline{A})(x),$$ for all $n \in \mathbb{N}$, $x \in \text{Ran }(\lambda I - A)$ and $B \in \mathcal{L}$. In order to prove that the remainder $$(-1)^{n + 1} \text{R}(\lambda, \overline{A})^{n + 1} (\text{ad }A)^{n + 1}(B) \, \text{R}(\lambda, \overline{A})(x)$$ goes to 0, as $n \longrightarrow + \infty$, it will be shown the existence of an $M > 0$ with the property that, for all $p \in \Gamma_A$, there exists $0 < r_p < 1$ in a way that $$p(\text{R}(\lambda, \overline{A})^k (\text{ad }A)^k (B)x) \leq M r_p^k \, \rho_{p, 1}(x), \qquad k \in \mathbb{N}, \, x \in \mathcal{D}.$$ To this purpose, it will be necessary to embed $\mathcal{L}$ into its complexification, $\mathcal{L}_\mathbb{C} := \mathcal{L} + i \mathcal{L}$, since results coming from the complex Banach algebras' realm will be used. Denote by $\text{R}_p(\lambda, \overline{A})$, as in Theorem 2.2, the unique everywhere-defined continuous extension of the densely defined $\| \, \cdot \, \|_p$-continuous linear operator $(\text{R}(\lambda, \overline{A}))_p$ induced by $\text{R}(\lambda, \overline{A})$ on the completion $\mathcal{X}_p$. Since for each fixed $1 \leq k \leq d$, $$(\text{ad }A)^k(B) = \sum_{l = 1}^d c^{(k)}_l B_l, \qquad c_l \in \mathbb{R},$$ one obtains $$p(\text{R}(\lambda, \overline{A})^k (\text{ad }A)^k(B)x) = \|\text{R}_p(\lambda, \overline{A})^k [(\text{ad }A)^k(B)x]_p\|_p$$ $$\leq \|\text{R}_p(\lambda, \overline{A})^k\|_p \, \|[(\text{ad }A)^k (B)x]_p\|_p = \|\text{R}_p(\lambda, \overline{A})^k\|_p \, p((\text{ad }A)^k(B)x)$$ $$\leq \|\text{R}_p(\lambda, \overline{A})^k\|_p \, |(\text{ad }A)^k(B)|_1 \, \rho_{p, 1}(x) \leq \|\text{R}_p(\lambda, \overline{A})^k\|_p \, \|(\tilde{\text{ad}}\text{ }A)^k\| \, |B|_1 \, \rho_{p, 1}(x),$$ for all $x \in \mathcal{D}$ - note that the majoration of $p((\text{ad }A)^k(B)x)$ is obtained in a similar way of that done in the proof of Theorem 2.2. Making $M := |B|_1$, it is sufficient to show that for some $0 < r_p < 1$ and for all sufficiently large $k$, $$\|\text{R}_p(\lambda, \overline{A})^k\|_p \, \|(\tilde{\text{ad}}\text{ }A)^k\| < r_p^k$$ or, equivalently, $$\|\text{R}_p(\lambda, \overline{A})^k\|_p^{1/k} \|(\tilde{\text{ad}}\text{ }A)^k\|^{1/k} < r_p.$$ Since $\mathcal{L}(\mathcal{X}_p)$ and $\mathcal{L}(\mathcal{L}_\mathbb{C})$ are both complex unital Banach algebras, by Gelfand's spectral radius formula it is sufficient to show that the product of $$\nu(\text{R}_p(\lambda, \overline{A})) := \lim_{k \rightarrow + \infty} \|\text{R}_p(\lambda, \overline{A})^k\|_p^{1/k}$$ and $$\nu(\tilde{\text{ad}}\text{ }A) := \lim_{k \rightarrow + \infty} \|(\tilde{\text{ad}}\text{ }A)^k\|^{1/k}$$ is strictly less than 1. By Lemma 1.6.1 and \cite[3.5 Generation Theorem (contraction case), page 73]{engel}, the linear operator $\left(\overline{A}\right)_p$ is closable and the resolvent of its closure satisfies the norm inequality $$\left\|\text{R}\left(\lambda, \overline{\left(\overline{A}\right)_p}\right)\right\|_p \leq \frac{1}{|\text{Re }\lambda|},$$ since $\overline{\left(\overline{A}\right)_p}$ generates a strongly continuous group of isometries on $\mathcal{X}_p$. It will now be proved that $$\text{R}\left(\lambda, \overline{\left(\overline{A}\right)_p}\right) = \text{R}_p(\lambda, \overline{A}).$$ So fix $x_p \in \text{Dom }\overline{\left(\overline{A}\right)_p}$ and take a net $\left\{[x_\alpha]_p\right\}$ in $\text{Dom}\left(\overline{A}\right)_p$ such that $$[x_\alpha]_p \xrightarrow{\alpha} x_p$$ and $$\left(\overline{A}\right)_p([x_\alpha]_p) \xrightarrow{\alpha} \overline{\left(\overline{A}\right)_p}(x_p).$$ This implies $$\text{R}_p(\lambda, \overline{A}) \left[\left(\lambda I - \overline{\left(\overline{A}\right)_p}\right)([x_\alpha]_p)\right] = (\text{R}(\lambda, \overline{A}))_p [(\lambda I - \left(\overline{A}\right)_p)([x_\alpha]_p)] = [x_\alpha]_p \longrightarrow x_p.$$ But $$\left(\lambda I - \overline{\left(\overline{A}\right)_p}\right)([x_\alpha]_p) = (\lambda I - \left(\overline{A}\right)_p)([x_\alpha]_p) \longrightarrow \left(\lambda I - \overline{\left(\overline{A}\right)_p}\right)(x_p),$$ which implies $$\text{R}_p(\lambda, \overline{A}) \left[\left(\lambda I - \overline{\left(\overline{A}\right)_p}\right)([x_\alpha]_p)\right] \longrightarrow \text{R}_p(\lambda, \overline{A}) \left[\left(\lambda I - \overline{\left(\overline{A}\right)_p}\right)(x_p)\right],$$ by the continuity of the operator $\text{R}_p(\lambda, \overline{A})$. By uniqueness of the limit, $$\text{R}_p(\lambda, \overline{A}) \left[\left(\lambda I - \overline{\left(\overline{A}\right)_p}\right)(x_p)\right] = x_p.$$ This proves that $\text{R}_p(\lambda, \overline{A})$ is a left inverse for $\left(\lambda I - \overline{\left(\overline{A}\right)_p}\right)$. To prove that $\text{R}_p(\lambda, \overline{A})$ is a right inverse fix an arbitrary $x_p \in \mathcal{X}_p$ and a net $\left\{[x_\alpha]_p\right\}$ in $\pi_p[\mathcal{D}]$ such that $[x_\alpha]_p \longrightarrow x_p$. Then, $$\left(\lambda I - \overline{\left(\overline{A}\right)_p}\right) \text{R}_p(\lambda, \overline{A}) ([x_\alpha]_p) = (\lambda I - \left(\overline{A}\right)_p) (\text{R}(\lambda, \overline{A}))_p ([x_\alpha]_p) = [x_\alpha]_p \longrightarrow x_p.$$ On the other hand, continuity of $\text{R}_p(\lambda, \overline{A})$ and closedness of $\overline{\left(\overline{A}\right)_p}$ imply that $$(\lambda I - \overline{\left(\overline{A}\right)_p}) \text{R}_p(\lambda, \overline{A}) ([x_\alpha]_p) \longrightarrow (\lambda I - \overline{\left(\overline{A}\right)_p}) \text{R}_p(\lambda, \overline{A}) (x_p).$$ As before, uniqueness of the limit establishes that $\text{R}_p(\lambda, \overline{A})$ is a right inverse for $\lambda I - \overline{\left(\overline{A}\right)_p}$. This proves $\text{R}_p(\lambda, \overline{A})$ is a two-sided inverse for the operator $\lambda I - \overline{\left(\overline{A}\right)_p}$ defined on $\text{Dom }\overline{\left(\overline{A}\right)_p}$ and that $\text{R}_p(\lambda, \overline{A}) = \left(\lambda I - \overline{\left(\overline{A}\right)_p}\right)^{-1}$. Since $\text{R}_p(\lambda, \overline{A})$ is a continuous operator, the equality $$\text{R}_p(\lambda, \overline{A}) = \text{R}\left(\lambda, \overline{\left(\overline{A}\right)_p}\right)$$ follows, as claimed.

Hence, combining all of these results with the hypothesis $$|\text{Re }\lambda| > \|\tilde{\text{ad}}\text{ }A\|,$$ one finally concludes that $$\|\text{R}_p(\lambda, \overline{A})^k\|_p^{1/k} \leq \|\text{R}_p(\lambda, \overline{A})\|_p = \|\text{R}(\lambda, \overline{\left(\overline{A}\right)_p})\|_p \leq \frac{1}{|\text{Re }\lambda|} < \frac{1}{\|\tilde{\text{ad}}\text{ }A\|}, \qquad k \in \mathbb{N}.$$ Noting that the sequence $\left\{\|\text{R}_p(\lambda, \overline{A})^k\|^{1/k}\right\}_{k \in \mathbb{N}}$ converges to $$\inf \left\{\|\text{R}_p(\lambda, \overline{A})^k\|_p^{1/k}: k \in \mathbb{N}\right\} = \nu(\text{R}_p(\lambda, \overline{A})),$$ one sees that $$\nu(\text{R}_p(\lambda, \overline{A})) < \frac{1}{\|\tilde{\text{ad}}\text{ }A\|}$$ and $$\nu(\text{R}_p(\lambda, \overline{A})) \, \nu(\tilde{\text{ad}}\text{ }A) < \frac{1}{\|\tilde{\text{ad}}\text{ }A\|} \, \|\tilde{\text{ad}}\text{ }A\| = 1.$$ This gives the desired result and shows that there exist $M > 0$ and $0 < r_p < 1$ such that $$p(\text{R}(\lambda, \overline{A})^k [(\text{ad }A)^k(B)x]) \leq M r_p^k \, \rho_{p, 1}(x), \qquad x \in \mathcal{D},$$ for sufficiently large $k$. Therefore, $$\sum_{k = 0}^{+ \infty} p(\text{R}(\lambda, \overline{A})^k (\text{ad }A)^k(B)x) < \infty, \qquad x \in \mathcal{D},$$ which implies $$\lim_{k \rightarrow + \infty} p(\text{R}(\lambda, \overline{A})^k (\text{ad }A)^k(B)x) = 0, \qquad x \in \mathcal{D}.$$ Since $p$ is arbitrary, $\text{R}(\lambda, \overline{A})^k (\text{ad }A)^k(B)(x) \longrightarrow 0$ in $\mathcal{X}$, for all $x \in \mathcal{D}$, and $$B \, \text{R}(\lambda, \overline{A})(x) = \sum_{k = 0}^{+ \infty} (-1)^k \text{R}(\lambda, \overline{A})^{k + 1} [(\text{ad }A)^k(B)](x), \qquad x \in \text{Ran }(\lambda I - A).$$
Now, the only step missing is to extend what was just proved to every element of $\mathcal{D}$. So let $y \in \mathcal{D}$. By the density hypothesis there exists a net $\left\{y_\alpha\right\}$ in $\text{Ran }(\lambda I - A)$ such that $\rho_{p, 1}(y_\alpha - y) \xrightarrow{\alpha} 0$, for all $p \in \Gamma_A$. Fix $p' \in \Gamma_A$. Using what was just proved one sees that, for all $\alpha$, $$\sum_{k = 0}^{+ \infty} p' \left((-1)^k \text{R}(\lambda, \overline{A})^{k + 1} [(\text{ad }A)^k (B)(y_{\alpha} - y)] \right)$$ $$\leq \rho_{p', 1}(y_{\alpha} - y) \, \|\text{R}_{p'}(\lambda, \overline{A})\|_{p'} \left[\sum_{k = 0}^{+ \infty} M r_{p'}^k \right].$$ Since $p'$ is arbitrary, taking limits on $\alpha$ on both sides gives, in particular, that $$\sum_{k = 0}^{+ \infty} (-1)^k \text{R}(\lambda, \overline{A})^{k + 1} [(\text{ad }A)^k (B)(y_{\alpha})] \longrightarrow \sum_{k = 0}^{+ \infty} (-1)^k \text{R}(\lambda, \overline{A})^{k + 1} [(\text{ad }A)^k (B)(y)].$$ Hence, $\left\{B \, \text{R}(\lambda, \overline{A})(y_{\alpha})\right\}$ converges and, by closedness of $\overline{B}$ and $\tau$-continuity of $\text{R}(\lambda, \overline{A})$, $$B \, \text{R}(\lambda, \overline{A})(y_{\alpha}) \xrightarrow{\alpha} \overline{B} \, \text{R}(\lambda, \overline{A})(y).$$ This implies $$\overline{B} \, \text{R}(\lambda, \overline{A})(y) = \sum_{k = 0}^{+ \infty} (-1)^k \text{R}(\lambda, \overline{A})^{k + 1} [(\text{ad }A)^k(B)](y), \qquad y \in \mathcal{D},$$ giving the desired result. \hfill $\blacksquare$ \\

Theorem 2.4 below is a ``locally convex version'' of \cite[Theorem 6.1, page 133]{jorgensenmoore}, when $A$ is assumed to be a pregenerator of an equicontinuous group.\\

\textbf{Theorem 2.4 (C$\bm{^1}$-Continuity of One-Parameter Groups):} \textit{\textbf{Suppose also that all the basis operators $\bm{\left(B_k\right)_{1 \leq k \leq d}}$ satisfy the (KIP) with respect to $\bm{\Gamma_A}$.} Then, $\overline{A}^{\tau_1}$ is the infinitesimal generator of (see Theorem 2.2 for the definition of this operator) a $\Gamma_{A, 1}$-group of bounded type $t \longmapsto \text{V}_1\left(t, \overline{A}^{\tau_1} \right)$ on $\mathcal{D}_1$. Moreover, $\text{V}(t, \overline{A})$ leaves $\mathcal{D}_1$ invariant, for all $t \in \mathbb{R}$, and restricts there as the $\tau_1$-continuous linear operator $\text{V}_1\left(t, \overline{A}^{\tau_1} \right)$, for all $t \in \mathbb{R}$.\\
Every $B \in \mathcal{L}$ can be seen as a continuous linear operator from $(\mathcal{D}, \tau_1)$ to $(\mathcal{X}, \tau)$, so denote by $\tilde{B}^{(1)} \colon \mathcal{D}_1 \longmapsto \mathcal{X}$ its unique continuous extension to all of $\mathcal{D}_1$. Then, for all closable $B \in \mathcal{L}$, $$(2.4.1) \qquad \tilde{B}^{(1)} \, V_1(t, A_1)(x) = V(t, \overline{A}) \, [\exp (-t \, \text{ad }A)(B)]\, \tilde{} \, ^{(1)}(x), \qquad t \in \mathbb{R}, \, x \in \mathcal{D}_1,$$ where $\exp (-t \, \text{ad }A)(B)$ is defined by the $|\, \cdot \,|_1$-convergent series $$\exp (-t \, \text{ad }A)(B) := \sum_{k = 0}^{+ \infty} \frac{(\text{ad }A)^k(B)}{k!} (-t)^k \in \mathcal{L}.$$}

\textbf{Proof of Theorem 2.4:} The first initiative of the proof is to obtain the hypotheses which are necessary to invoke the $(\Leftarrow)$ implication of \cite[Theorem 4.2]{babalola}, but for the complete Hausdorff locally convex space $(\mathcal{D}_1, \tau_1)$ (note that its topology is generated by a saturated family of seminorms). Let $\lambda \in \mathbb{C}$ satisfy $|\text{Re }\lambda| \geq M_A := \max \left\{|\text{Re }\lambda_{(-, A)}|, |\text{Re }\lambda_{(+, A)}|\right\} > 0$. The additional hypothesis that all the basis operators $\left(B_k\right)_{1 \leq k \leq d}$ satisfy the (KIP) with respect to $\Gamma_A$ will already be used in the following argument: $\overline{A}$ possesses the (KIP) with respect to $\Gamma_A$ and, since all the basis elements also possess the (KIP) with respect to $\Gamma_A$, the operator $$\overline{A}^{\tau_1} \colon \text{Dom }\overline{A}^{\tau_1} \subseteq \mathcal{D}_1 \longrightarrow \mathcal{D}_1$$ is a $\tau_1$-densely defined $\tau_1$-closed linear operator on $\mathcal{D}_1$ which possesses the (KIP) with respect to $\Gamma_{A, 1}$, since $\overline{A}^{\tau_1} \subset \overline{A}$: in fact, if $\rho_{p, 1} \in \Gamma_{A, 1}$, $N_{\rho_{p, 1}} := \left\{x \in \mathcal{D}_1: \rho_{p, 1}(x) = 0\right\}$ and $x \in \text{Dom }\overline{A}^{\tau_1} \cap N_{\rho_{p, 1}}$, then $p(x) = 0$, so $x \in N_p$, which implies $\overline{A}^{\tau_1}x = \overline{A}x \in N_p$; by the (KIP) of the basis elements, it follows that $B_k(\overline{A}^{\tau_1}x) = B_k(\overline{A}x) \in N_p$, for all $1 \leq k \leq d$, so $\overline{A}^{\tau_1}x \in N_{\rho_{p, 1}}$, proving $\overline{A}^{\tau_1}$ has the (KIP) with respect to $\Gamma_{A, 1}$.

Therefore,\footnote{Remember $(\mathcal{D}_1)_{\rho_{p, 1}} = \overline{\mathcal{D}_1 / N_{\rho_{p, 1}}}^{\| \, \cdot \, \|_{\rho_{p, 1}}}$ denotes the Banach space completion of $\mathcal{D}_1 / N_{\rho_{p, 1}}$.} $$\left(\overline{A}^{\tau_1}\right)_{\rho_{p, 1}} \colon \pi_{\rho_{p, 1}}[\text{Dom }\overline{A}^{\tau_1}] \longrightarrow (\mathcal{D}_1)_{\rho_{p, 1}}$$ is a (well-defined) densely defined linear operator, for all $\rho_{p, 1} \in \Gamma_{A, 1}$. (For the rest of the proof, the more economic notations $A_1 := \overline{A}^{\tau_1}$ and $A_{\rho_{p, 1}} := (A_1)_{\rho_{p, 1}} = \left(\overline{A}^{\tau_1}\right)_{\rho_{p, 1}}$ are going to be employed, in order to facilitate the reading process) To see that it is $\| \, \cdot \, \|_{\rho_{p, 1}}$-closable, note first that $$p((\text{ad }A)(B_k)x) \leq |(\text{ad }A)(B_j)|_1 \, \rho_{p, 1}(x) \leq \|\tilde{\text{ad}}\text{ }A\| \, \rho_{p, 1}(x),$$ for all $p \in \Gamma_A$, $1 \leq k \leq d$ and $x \in \mathcal{D}$, since $|B_k|_1 = 1$. Hence, since $A$ is $\Gamma_A$-conservative, one obtains for each $p \in \Gamma_A$ and $1 \leq k \leq d$ that $$p(B_k(\mu I - A)x) = p((\mu I - A) B_k(x) + (\text{ad }A)(B_k)(x))$$ $$\geq p((\mu I - A) B_k(x)) - p((\text{ad }A)(B_k)x)$$ $$\geq |\mu| \, p(B_k(x)) - \|\tilde{\text{ad}}\text{ }A\| \, \rho_{p, 1}(x), \qquad |\mu| > \|\tilde{\text{ad}}\text{ }A\|, \, x \in \mathcal{D}.$$ Taking the maximum over $0 \leq k \leq d$ produces $$\rho_{p, 1}((\mu I - A)x) \geq (|\mu| - \|\tilde{\text{ad}}\text{ }A\|) \, \rho_{p, 1}(x), \qquad |\mu| > \|\tilde{\text{ad}}\text{ }A\|, \, x \in \mathcal{D}.\footnote{Note that the case $k = 0$ is just a consequence of $\Gamma_A$-conservativity of $A$.}$$ Now, let $\left\{[x_\alpha]_{\rho_{p, 1}}\right\}$ be a net in $\pi_p[\text{Dom }A_1]$ and $\left\{A_{\rho_{p, 1}}([x_\alpha]_{\rho_{p, 1}})\right\} = \left\{[A_1(x_\alpha)]_{\rho_{p, 1}}\right\}$ be a net in $\mathcal{D}_1 / N_{\rho_{p, 1}}$, where $$[x_\alpha]_{\rho_{p, 1}} \xrightarrow{\alpha} 0$$ and $$[A_1(x_\alpha)]_{\rho_{p, 1}} \xrightarrow{\alpha} y \in (\mathcal{D}_1)_{\rho_{p, 1}},$$ both convergences being in the $\| \, \cdot \, \|_{\rho_{p, 1}}$-topology. By the above inequality, one has for all $\mu > \|\tilde{\text{ad}}\text{ }A\|$ and $[x']_{\rho_{p, 1}} \in \pi_{\rho_{p, 1}}[\mathcal{D}]$ that $$\left\|\left(I - \frac{1}{\mu} A_{\rho_{p, 1}}\right)\left(\left[x_\alpha + \frac{1}{\mu} x'\right]_{\rho_{p, 1}}\right)\right\|_{\rho_{p, 1}} = \left\|\left[\left(I - \frac{1}{\mu} A_1\right)\left(x_\alpha + \frac{1}{\mu} x'\right)\right]_{\rho_{p, 1}}\right\|_{\rho_{p, 1}}$$ $$= \rho_{p, 1}\left(\left(I - \frac{1}{\mu} A_1\right)\left(x_\alpha + \frac{1}{\mu} x'\right)\right) \geq \frac{\mu - \|\tilde{\text{ad}}\text{ }A\|}{\mu} \, \rho_{p, 1}\left(x_\alpha + \frac{1}{\mu} x'\right)$$ $$= \frac{\mu - \|\tilde{\text{ad}}\text{ }A\|}{\mu} \, \left\|\left[x_\alpha + \frac{1}{\mu} x'\right]_{\rho_{p, 1}}\right\|_{\rho_{p, 1}}.$$ Hence, $$\left\|[x_\alpha]_{\rho_{p, 1}} + \frac{1}{\mu} [x']_{\rho_{p, 1}} - \frac{1}{\mu} A_{\rho_{p, 1}}([x_\alpha]_{\rho_{p, 1}}) - \frac{1}{\mu^2} A_{\rho_{p, 1}}([x']_{\rho_{p, 1}})\right\|_{\rho_{p, 1}}$$ $$\geq \frac{\mu - \|\tilde{\text{ad}}\text{ }A\|}{\mu} \, \left\|[x_\alpha]_{\rho_{p, 1}} + \frac{1}{\mu} [x']_{\rho_{p, 1}}\right\|_{\rho_{p, 1}}.$$ Taking limits on $\alpha$ on both members and multiplying by $\mu$, one obtains $$\left\|[x']_{\rho_{p, 1}} - y - \frac{1}{\mu} A_{\rho_{p, 1}}([x']_{\rho_{p, 1}})\right\|_{\rho_{p, 1}} \geq \frac{\mu - \|\tilde{\text{ad}}\text{ }A\|}{\mu} \, \|[x']_{\rho_{p, 1}}\|_{\rho_{p, 1}}.$$ Sending $\mu$ to $+ \infty$ and using the density of $\pi_{\rho_{p, 1}}[\mathcal{D}]$ in $(\mathcal{D}_1)_{\rho_{p, 1}}$, it follows that $\|y\|_{\rho_{p, 1}} = 0$.  This shows $y = 0$ and proves the $\| \, \cdot \, \|_{\rho_{p, 1}}$-closability of $A_{\rho_{p, 1}}$.\\

To prove the necessary resolvent bounds, substitute $B$ by $B_j$ on equation (2.3.1) of Theorem 2.3, $1 \leq j \leq d$, where $\left\{B_j\right\}_{1 \leq j \leq d}$ is the fixed basis of $\mathcal{L}$. This gives $$\overline{B_j} \, \text{R}(\lambda, \overline{A})(x) = \sum_{k = 0}^{+ \infty} (-1)^k \text{R}(\lambda, \overline{A})^{k + 1} [(\text{ad }A)^k (B_j)](x), \qquad 1 \leq j \leq d, \, x \in \mathcal{D}.$$ Now, since Theorem 2.3 implies $\text{R}_p(\lambda, \overline{A})$ is the $\lambda$-resolvent of $\overline{\left(\overline{A}\right)_p}$, which is the infinitesimal generator of a group of isometries, $$\frac{1}{\|\tilde{\text{ad}}\text{ }A\|} > \frac{1}{|\text{Re }\lambda|} \geq \|\text{R}_p(\lambda, \overline{A})\|_p.$$ Hence, $\|\text{R}_p(\lambda, \overline{A})\|_p \, \|\tilde{\text{ad}}\text{ }A\| < 1$. This implies $$p(\overline{B_j} \, \text{R}(\lambda, \overline{A})x) \leq \|\text{R}_p(\lambda, \overline{A})\|_p \, \left[\sum_{k = 0}^{+ \infty} p(\text{R}(\lambda, \overline{A})^k [(\text{ad }A)^k (B_j)]x) \right]$$ $$\leq \|\text{R}_p(\lambda, \overline{A})\|_p \, \left[\sum_{k = 0}^{+ \infty} \|\text{R}_p(\lambda, \overline{A})\|_p^k \, \|\tilde{\text{ad}}\text{ }A\|^k \, |B_j|_1 \, \rho_{p, 1}(x) \right]$$ $$= \frac{\|\text{R}_p(\lambda, \overline{A})\|_p}{1 - \left\|\text{R}_p(\lambda, \overline{A})\right\|_p \, \|\tilde{\text{ad}}\text{ }A\|} \, \rho_{p, 1}(x) = \frac{1}{\frac{1}{\left\|\text{R}_p(\lambda, \overline{A})\right\|_p} - \|\tilde{\text{ad}}\text{ }A\|} \, \rho_{p, 1}(x)$$ $$\leq \frac{1}{|\text{Re }\lambda| - \|\tilde{\text{ad}}\text{ }A\|} \, \rho_{p, 1}(x), \qquad p \in \Gamma_A, \, 1 \leq j \leq d, \, x \in \mathcal{D},$$ since $|B_j|_1 = 1$. Also, $$p(\overline{B_0} \, \text{R}(\lambda, \overline{A})x) = p(\text{R}(\lambda, \overline{A})x) \leq \|\text{R}_p(\lambda, \overline{A})\|_p \, p(x)$$ $$\leq \frac{\|\text{R}_p(\lambda, \overline{A})\|_p}{1 - \|\text{R}_p(\lambda, \overline{A})\|_p \, \|\tilde{\text{ad}}\text{ }A\|} \, \rho_{p, 1}(x) \leq \frac{1}{|\text{Re }\lambda| - \|\tilde{\text{ad}}\text{ }A\|} \, \rho_{p, 1}(x), \qquad p \in \Gamma_A, \, x \in \mathcal{D},$$ because $0 < 1 - \|\text{R}_p(\lambda, \overline{A})\|_p \, \|\tilde{\text{ad}}\text{ }A\| < 1$. Therefore, taking the maximum over $0 \leq j \leq d$ and using the fact (proved in Theorem 2.2) that $\left. \text{R}(\lambda, \overline{A})\right|_{\mathcal{D}_1} = \text{R}(\lambda, A_1)$, a $\tau_1$-density argument provides $$\rho_{p, 1}(\text{R}(\lambda, A_1)x) \leq \frac{1}{|\text{Re }\lambda| - \|\tilde{\text{ad}}\text{ }A\|} \, \rho_{p, 1}(x), \qquad p \in \Gamma_A, \, x \in \mathcal{D}_1,$$ showing that the operator $\left. \text{R}(\lambda, \overline{A})\right|_{\mathcal{D}_1}$ possesses the (KIP) relatively to the family $\Gamma_{A, 1}$. This yields $$\left\|\left(\text{R}(\lambda, A_1)\right)_{\rho_{p, 1}} \left( [x]_{\rho_{p, 1}} \right)\right\|_{\rho_{p, 1}} \leq \frac{1}{|\text{Re }\lambda| - \|\tilde{\text{ad}}\text{ }A\|} \, \left\|[x]_{\rho_{p, 1}}\right\|_{\rho_{p, 1}}, \qquad p \in \Gamma_A, \, [x]_{\rho_{p, 1}} \in \mathcal{D}_1 / N_{\rho_{p, 1}}.$$ Since $\mathcal{D}_1 / N_{\rho_{p, 1}}$ is $\| \, \cdot \, \|_{\rho_{p, 1}}$-dense in $(\mathcal{D}_1)_{\rho_{p, 1}}$, the above estimate allows one to define by limits the unique continuous linear operator $\text{R}_{\rho_{p, 1}}(\lambda, A_1)$ on $(\mathcal{D}_1)_{\rho_{p, 1}}$ which extends $\text{R}(\lambda, A_1)$. Arguing in a similar manner of that done at the end of the first paragraph of Theorem 2.3 (exploring the closedness of the operator $\overline{A_{\rho_{p, 1}}}$), it follows that $$\text{R}_{\rho_{p, 1}}(\lambda, A_1) = \text{R}(\lambda, \overline{A_{\rho_{p, 1}}}).$$ Hence, $$\left\|\text{R}(\lambda, \overline{A_{\rho_{p, 1}}})(x_{\rho_{p, 1}})\right\|_{\rho_{p, 1}} \leq \frac{1}{|\text{Re }\lambda| - \|\tilde{\text{ad}}\text{ }A\|} \, \|x_{\rho_{p, 1}}\|_{\rho_{p, 1}}, \qquad p \in \Gamma_A, \, x_{\rho_{p, 1}} \in (\mathcal{D}_1)_{\rho_{p, 1}}.$$ In particular, $$\left\|[\text{R}(\lambda, \overline{A_{\rho_{p, 1}}})]^n(x_{\rho_{p, 1}})\right\|_{\rho_{p, 1}} \leq \frac{1}{(|\lambda| - \|\tilde{\text{ad}}\text{ }A\|)^n} \, \|x_{\rho_{p, 1}}\|_{\rho_{p, 1}}, \qquad p \in \Gamma_A, \, n \in \mathbb{N}, \, x_{\rho_{p, 1}} \in (\mathcal{D}_1)_{\rho_{p, 1}},$$ for all $\lambda \in \mathbb{R}$ satisfying $|\lambda| > \|\tilde{\text{ad}}\text{ }A\|$, which finishes the verification of the last hypothesis needed to apply Hille-Yosida-Phillips theorem for locally convex spaces (\cite[Theorem 4.2]{babalola}),\footnote{Actually, a straightforward adaptation to a group version of \cite[Theorem 4.2]{babalola} must be performed, by considering projective limits of strongly continuous groups, instead of semigroups - note that \cite[Theorem 4.2]{babalola} regards generation of semigroups.} so that $A_1$ is the generator of a $\Gamma_{A_1}$-group $t \longmapsto \text{V}_1(t, A_1)$ on $\mathcal{D}_1$.\\

Fix $\lambda \in \mathbb{C}$ satisfying $\text{Re }\lambda \geq M_A := \max \left\{|\text{Re }\lambda_{(-, A)}|, |\text{Re }\lambda_{(+, A)}|\right\} > \|\tilde{\text{ad}}\text{ }A\|$.\\

Since $$\text{R}(\lambda, A_1) = \left. \text{R}(\lambda, \overline{A}) \right|_{\mathcal{D}_1},$$ it follows for all $x \in \mathcal{D}_1$ that $$\int_0^{+ \infty} e^{- \lambda s} V(s, \overline{A})x \, ds = \text{R}(\lambda, \overline{A})x = \text{R}(\lambda, A_1)x = \int_0^{+ \infty} e^{- \lambda s} V_1(s, A_1)x \, ds,$$ where the member on the left-hand side is a strongly $\tau$-convergent integral and the member on the right-hand side is a strongly $\tau_1$-convergent integral. Therefore, for fixed $x \in \mathcal{D}_1$ and $f \in \mathcal{X}'$, $$\int_0^{+ \infty} e^{- \lambda s} \, f \left(V(s, \overline{A})x - V_1(s, A_1)x\right) \, ds = 0.$$ As the function $$f_A \colon s \longmapsto f \left(V(s, \overline{A})x - V_1(s, A_1)x\right)$$ is continuous on $[0, + \infty)$ and is such that $\left|f_A(s)\right| \leq C e^{as}$, for all $s \in [0, + \infty)$ and some $C \geq 0$, $a \in \mathbb{R}$, it follows from the injectivity of the Laplace transform that $f \left(V(s, \overline{A})x \right) = f \left(V_1(s, A_1)x\right)$, for all $s \in [0,+ \infty)$.\footnote{See, for example \cite[Corollary 8.1, page 267]{follandfourier}.} By a corollary of Hahn-Banach's theorem, it follows that $\left. V(s, \overline{A}) \right|_{\mathcal{D}_1} = V_1(s, A_1)$, for all $s \geq 0$. But if $s < 0$, then $$V(s, \overline{A})x = (V(-s, \overline{A}))^{-1}x = (V_1(-s, A_1))^{-1}x = V_1(s, A_1)x, \qquad x \in \mathcal{D}_1,$$ so $$\left. V(s, \overline{A}) \right|_{\mathcal{D}_1} = V_1(s, A_1),$$ for all $s \in \mathbb{R}$, proving that $s \longmapsto V(s, \overline{A})$ restricts to a $\Gamma_{A, 1}$-group on $\mathcal{D}_1$.\\

Before proceeding it is important to note that, for all $x \in \mathcal{D}$, the sequence $$\left\{\sum_{k = 0}^n \frac{(\text{ad }A)^k(B)}{k!} (-t)^k (x)\right\}_{n \in \mathbb{N}}$$ in $\mathcal{D}$ converges to the element $$\sum_{k = 0}^{+ \infty} \frac{(\text{ad }A)^k(B)}{k!} (-t)^k (x) \in \mathcal{D}$$ - see \cite[Theorem 3.2, page 64]{jorgensenmoore} - because the family $\left\{\mathcal{L} \ni C \longmapsto p(Cx): p \in \Gamma\right\}$ of seminorms generates a Hausdorff topology on $\mathcal{L}$ and any two Hausdorff topologies on a finite-dimensional vector space are equivalent.\\

The next step is to prove equality $(2.4.1)$. Note, first, that $t \longmapsto \text{V}_1(t, A_1)$ is of \textbf{bounded type}: if $w_{\rho_{p, 1}}$ denotes the type of the semigroup on $(\mathcal{D}_1)_{\rho_{p, 1}}$ induced (after being extended by limits)\footnote{Note that the operators $V_1(t, A_1)$, $t \geq 0$, possess the (KIP) with respect to $\Gamma_{A, 1}$, by the very definition of a $\Gamma_{A, 1}$-group.} by $[0, +\infty) \ni t \longmapsto V_1(t, A_1)$, then $w_{\rho_{p, 1}} \leq \|\tilde{\text{ad}}\text{ }A\|$, for all $\rho_{p, 1} \in \Gamma_{A, 1}$. This implies $$w := \sup \left\{w_{\rho_{p, 1}}: \rho_{p, 1} \in \Gamma_{A, 1}\right\} \leq \|\tilde{\text{ad}}\text{ }A\|,$$ so $t \longmapsto V_1(t, A_1)$ is of bounded type, following the terminology of \cite[page 170]{babalola}. Fix, as before, $\lambda \in \mathbb{C}$ satisfying $\text{Re }\lambda \geq M_A$. By what was just noted, \cite[Theorem 3.3]{babalola} is applicable and the relation $$\text{R}(\lambda, A_1)x = \int_0^{+ \infty} e^{-\lambda s} \, V_1(s, A_1)x \, ds$$ is valid, for all $x \in \mathcal{D}_1$, where the integral is $\tau_1$-convergent. Also, $(\tau_1 \times \tau)$-continuity of $\tilde{B}^{(1)}$ ensures that $$(2.4.2) \qquad \tilde{B}^{(1)} \, \text{R}(\lambda, A_1)x = \int_0^{+ \infty} e^{-\lambda s} \, \tilde{B}^{(1)} \, V_1(s, A_1)x \, ds, \qquad x \in \mathcal{D}_1,$$ the integral being $\tau$-convergent. The next step will be to prove that $$\int_0^{+ \infty} e^{-\lambda s} \, V(s, \overline{A}) \, [\exp (-s \, \text{ad }A)(B)]\, \tilde{} \, ^{(1)}(x) \, ds$$ $$= \sum_{n = 0}^{+ \infty} \left[\int_0^{+ \infty} e^{-\lambda s} \, \frac{(-s)^n}{n!} \, V(s, \overline{A}) \, (\text{ad }A)^n (B)(x) \, ds\right], \qquad x \in \mathcal{D}_1.$$ Fix $x \in \mathcal{D}$, $f \in \mathcal{X}'$ and $p' \in \Gamma_A$ such that $|f(x)| \leq M_f \, p'(x)$, for all $x \in \mathcal{X}$. Now, for all $t \in \mathbb{R}$, $$V(t, \overline{A}) \, [\exp (-t \, \text{ad }A)(B)]\, \tilde{} \, ^{(1)}(x) = \sum_{n = 0}^{+ \infty} \frac{(-t)^n}{n!} \, V(t, \overline{A}) \, (\text{ad }A)^n (B)(x),$$ by $\tau$-continuity of the operator $V(t, \overline{A})$, and the integral $$\int_0^{+ \infty} \left[\sum_{n = 0}^{+ \infty} \left|f \left(e^{-\lambda s} \, V(s, \overline{A}) \, \frac{(-s)^n}{n!} \, (\text{ad }A)^n (B)(x)\right)\right| \, ds\right]$$ $$\leq M_f \int_0^{+ \infty} e^{-\lambda s} \left[\sum_{n = 0}^{+ \infty} \frac{s^n}{n!} \, p'(V(s, \overline{A}) \, (\text{ad }A)^n (B)(x)) \, ds\right]$$ $$\leq M_f \int_0^{+ \infty} e^{-\lambda s} \left[\sum_{n = 0}^{+ \infty} \frac{s^n}{n!} \, \|(V(s, \overline{A}))_{p'}\|_{p'} \, p'((\text{ad }A)^n (B)(x)) \, ds\right]$$ $$\leq M_f \int_0^{+ \infty} e^{-\lambda s} \left[\sum_{n = 0}^{+ \infty} \frac{s^n}{n!} \, |(\text{ad }A)^n (B)|_1 \, \rho_{p', 1}(x) \, ds\right]$$ $$\leq M_f \int_0^{+ \infty} e^{-\lambda s} \left[\sum_{n = 0}^{+ \infty} \frac{(s \, \|\tilde{\text{ad}}\text{ }A\|)^n}{n!} \, |B|_1 \, \rho_{p', 1}(x) \, ds\right] = M_f \, |B|_1 \, \rho_{p', 1}(x) \int_0^{+ \infty} e^{(-\lambda + \|\tilde{\text{ad}}\text{ }A\|) s} \, ds$$ converges (since the last integral converges), where it was used that the semigroup on $\mathcal{X} / N_{p'}$ induced by $[0, + \infty) \ni t \longmapsto V(t, \overline{A})$, being a contraction semigroup, satisfies $$\|(V(s, \overline{A}))_{p'}\|_{p'} \leq 1, \qquad s \geq 0.$$ The convergence is ensured by the hypothesis $- \text{Re }\lambda + \|\tilde{\text{ad}}\text{ }A\| < 0$. This justifies an application of Fubini's Theorem, from which it follows that $$\int_0^{+ \infty} f \left(e^{-\lambda s} \, V(s, \overline{A}) \, [\exp (-s \, \text{ad }A)(B)]\, \tilde{} \, ^{(1)}(x) \, ds\right)$$ $$= \int_0^{+ \infty} \left[\sum_{n = 0}^{+ \infty} f \left(e^{-\lambda s} \, V(s, \overline{A}) \, \frac{(-s)^n}{n!} \, (\text{ad }A)^n (B)(x)\right) \, ds\right]$$ $$= \sum_{n = 0}^{+ \infty} \left[\int_0^{+ \infty} f \left(e^{-\lambda s} \, \frac{(-s)^n}{n!} \, V(s, \overline{A}) \, (\text{ad }A)^n (B)(x)\right) \, ds\right].$$ Finally, a Hahn-Banach's Theorem corollary gives the desired equality $$\int_0^{+ \infty} e^{-\lambda s} \, V(s, \overline{A}) \, [\exp (-s \, \text{ad }A)(B)]\, \tilde{} \, ^{(1)}(x) \, ds$$ $$= \sum_{n = 0}^{+ \infty} \left[\int_0^{+ \infty} e^{-\lambda s} \, \frac{(-s)^n}{n!} \, V(s, \overline{A}) \, (\text{ad }A)^n (B)(x) \, ds\right].$$ On the other hand, a resolvent formula for closed linear operators on Banach spaces\footnote{To be more precise, let $T \colon \text{Dom }T \subseteq \mathcal{Y} \longrightarrow \mathcal{Y}$ be a (not necessarily densely defined) closed linear operator on a Banach space $\mathcal{Y}$. If $\lambda$ is in the resolvent set of $T$, then $$\frac{d^n}{d \lambda^n} \text{R}(\lambda, T) = (-1)^n \, n! \, \text{R}(\lambda, T)^{n + 1}, \qquad n \in \mathbb{N}$$ - see \cite[Proposition 1.3 (ii), page 240]{engel}, for example.} implies that $$\frac{d^n}{d \lambda^n} \text{R}_p(\lambda, \overline{A})([x]_p) = \frac{d^n}{d \lambda^n} \text{R}\left(\lambda, \overline{\left(\overline{A}\right)_p}\right)([x]_p) = (-1)^n \, n! \, \text{R}\left(\lambda, \overline{\left(\overline{A}\right)_p}\right)^{n + 1}([x]_p)$$ $$= (-1)^n \, n! \, \text{R}_p(\lambda, \overline{A})^{n + 1}([x]_p), \qquad p \in \Gamma_A, \, n \in \mathbb{N}.$$ Thus $$\frac{d^n}{d \lambda^n} \text{R}\left(\lambda, \overline{A}\right)(x) = (-1)^n \, n! \, \text{R}(\lambda, \overline{A})^{n + 1}(x),$$ for all $n \in \mathbb{N}$, so $$\sum_{n = 0}^{+ \infty} \left[\int_0^{+ \infty} e^{-\lambda s} \, \frac{(-s)^n}{n!} \, V(s, \overline{A}) \, (\text{ad }A)^n (B)(x) \, ds\right]$$ $$= \sum_{n = 0}^{+ \infty} \frac{1}{n!} \left[\int_0^{+ \infty} \frac{d^n \, e^{-\lambda s}}{d \lambda^n} \, V(s, \overline{A}) \, (\text{ad }A)^n (B)(x) \, ds\right]$$ $$= \sum_{n = 0}^{+ \infty} \frac{1}{n!} \left[\frac{d^n}{d \lambda^n} \int_0^{+ \infty} e^{-\lambda s} \, V(s, \overline{A}) \, (\text{ad }A)^n (B)(x) \, ds\right] = \sum_{n = 0}^{+ \infty} \frac{1}{n!} \left[\frac{d^n}{d \lambda^n} \text{R}(\lambda, \overline{A}) \, (\text{ad }A)^n (B)(x)\right]$$ $$= \sum_{n = 0}^{+ \infty} (-1)^n \text{R}(\lambda, \overline{A})^{n + 1} \, (\text{ad }A)^n (B)(x) = \overline{B} \, \text{R}(\lambda, \overline{A})(x) = \overline{B} \, \text{R}(\lambda, A_1)(x),$$ by formula (2.3.1) and $\text{R}(\lambda, A_1) = \left. \text{R}(\lambda, \overline{A}) \right|_{\mathcal{D}_1}$. Since $\tilde{B}^{(1)} \subset \overline{B}$, the above equality together with equation (2.4.2) gives $$\int_0^{+ \infty} e^{-\lambda s} \, V(s, \overline{A}) \, [\exp (-s \, \text{ad }A)(B)]\, \tilde{} \, ^{(1)}(x) \, ds = \int_0^{+ \infty} e^{-\lambda s} \, \tilde{B}^{(1)} \, V_1(s, A_1)x \, ds.$$ Using a corollary of Hahn-Banach's Theorem and injectivity of the Laplace transform, as before, formula (2.4.1) is proved for $t \geq 0$ and $x \in \mathcal{D}$.\\
Repeating the arguments of this last paragraph with $-A_1$, instead of $A_1$, and using the identity $V_1(t, -A_1) = V(-t, A_1)$, for all $t \geq 0$, gives the equality $$\tilde{B}^{(1)} \, V_1(t, A_1)(x) = V(t, \overline{A}) \, [\exp (-t \, \text{ad }A)(B)]\, \tilde{} \, ^{(1)}(x), \qquad t \leq 0, \, x \in \mathcal{D}.$$ Now, the only thing missing is to extend formula (2.4.1) to all of $\mathcal{D}_1$, but this is easily accomplished by a $\tau_1$-density argument.\hfill $\blacksquare$ \\

With the aid of the four theorems above one can finally prove the theorem which justifies the title of this subsection. Its proof follows the lines of \cite[Theorem 7.4, page 161]{jorgensenmoore}, establishing a similar result in the more general setting of complete Hausdorff locally convex spaces:\\

\textbf{Theorem 2.5 (Construction of the Group Invariant C$\bm{^\infty}$ Domain):} \textit{Let $(\mathcal{X}, \tau)$ be a complete Hausdorff locally convex space, $\mathcal{D}$ a dense subspace of $\mathcal{X}$ and $\mathcal{L} \subseteq \text{End}(\mathcal{D})$ a finite dimensional \textbf{real} Lie algebra of linear operators. Fix an ordered basis $\left(B_k\right)_{1 \leq k \leq d}$ of $\mathcal{L}$ and define a basis-dependent norm on $\mathcal{L}$ by $$\left|\sum_{1 \leq k \leq d} b_k B_k\right|_1 := \sum_{1 \leq k \leq d} |b_k|.$$ Suppose $\mathcal{L}$ is generated, as a Lie algebra, by a finite set $\mathcal{S}$ of infinitesimal pregenerators of \textbf{equicontinuous} groups and denote by $t \longmapsto V(t, \overline{A})$ the respective one-parameter group generated by $\overline{A}$, for all $A \in \mathcal{S}$. For each $A \in \mathcal{S}$ choose a fundamental system of seminorms $\Gamma_A$ for $\mathcal{X}$ with respect to which the operator $\overline{A}$ has the (KIP), is $\Gamma_A$-conservative and $V( \, \cdot \, , \overline{A})$ is $\Gamma_A$-isometrically equicontinuous. Define $\Gamma_{A, 1} := \left\{\rho_{p, 1}: p \in \Gamma_A\right\}$. Now, assume that the following hypotheses are verified:}

\begin{enumerate}

\item \textit{the basis $\left(B_k\right)_{1 \leq k \leq d}$ is formed by closable elements \textbf{having the (KIP) with respect to} $$\Gamma_\mathcal{S} := \bigcup_{A \in \mathcal{S}} \Gamma_A;$$}
\item \textit{for each $A \in \mathcal{S}$ there exist two complex numbers $\lambda_{(-, A)}$, $\, \lambda_{(+, A)}$ satisfying $\text{Re }\lambda_{(-, A)} < - \|\tilde{\text{ad}}\text{ }A\|$ and $\text{Re }\lambda_{(+, A)} > \|\tilde{\text{ad}}\text{ }A\|$, where $\|\tilde{\text{ad}}\text{ }A\|$ denotes the usual operator norm of the extension $\tilde{\text{ad}}\text{ }A \colon \mathcal{L}_\mathbb{C} \longrightarrow \mathcal{L}_\mathbb{C}$ such that, for all $\lambda \in \mathbb{C}$ satisfying $$|\text{Re }\lambda| \geq M_A := \max \left\{|\text{Re }\lambda_{(-, A)}|, |\text{Re }\lambda_{(+, A)}|\right\},$$ the subspace $\text{Ran }(\lambda I - A)$ is dense in $\mathcal{D}$ with respect to the topology $\tau_1$ induced by the family $\Gamma_{A, 1}$.}

\end{enumerate}

\textit{Define $\mathcal{D}_\mathcal{S}$ as the smallest subspace of $\mathcal{X}$ (with respect to set inclusion) that contains $\mathcal{D}$ and is left invariant by the operators $\left\{V(t, \overline{A}): A \in \mathcal{S}, t \in \mathbb{R}\right\}$ - such subspace is well-defined since, by Theorem 2.4, $\mathcal{D}_1$ is left invariant by $V(t, \overline{A})$, for all $t \in \mathbb{R}$ and $A \in \mathcal{S}$; also, $\mathcal{D}_\mathcal{S} \subseteq \mathcal{D}_1$. Then:}

\begin{enumerate}

\item \textit{$\mathcal{D}_\mathcal{S}$ is a C$^\infty$ domain for $\overline{\mathcal{S}} := \left\{\overline{A}: A \in \mathcal{S}\right\}$, so $\mathcal{D}_\mathcal{S} \subseteq \mathcal{X}_\infty(\overline{\mathcal{S}})$ - for the definition of $\mathcal{X}_\infty(\overline{\mathcal{S}})$, see the beginning of this subsection.}
\item \textit{The estimates $$(2.5.1) \qquad \rho_{p, n}(V(t, \overline{A})x) \leq \exp (n \|\text{ad }A\| \, |t|) \, \rho_{p, n}(x)$$ hold for every $A \in \mathcal{S}$, $p \in \Gamma_A$, $n \in \mathbb{N}$, $t \in \mathbb{R}$ and $x \in \mathcal{D}_\mathcal{S}$, where $\rho_{p, n}$ was defined at the beginning of this subsection. Define $\mathcal{D}_n$ as the $\tau_n$-closure of $\mathcal{D}$ inside the complete Hausdorff locally convex space $(\mathcal{X}_n, \tau_n)$, for each $n \geq 1$. Then, $\mathcal{D}_n$ is left invariant by the operators $\left\{V(t, \overline{A}): A \in \mathcal{S}, t \in \mathbb{R}\right\}$ so, in particular, the one-parameter groups $t \longmapsto V(t, \overline{A})$ restrict to exponentially equicontinuous groups on $\mathcal{D}_n$ of type $\left\{w_p\right\}_{p \in \Gamma_A}$, where $w_p \leq n \|\text{ad }A\|$, for each $n \geq 1$ and $p \in \Gamma_A$. By formula (2.5.1), they are all locally equicontinuous one-parameter groups.}
\item \textit{The complete Hausdorff locally convex space $$\mathcal{D}_\infty := \bigcap_{n = 1}^{+ \infty} \mathcal{D}_n,$$ equipped with the $\tau_\infty$-topology generated by the family $$\Gamma_{A, \infty} := \bigcup_{n = 0}^{+ \infty} \Gamma_{A, n}$$ of seminorms, for some fixed $A \in \mathcal{S}$ (note that this topology does not depend on the choice of $A \in \mathcal{S}$), is left invariant by the operators $\left\{V(t, \overline{A}): A \in \mathcal{S}, t \in \mathbb{R}\right\}$, and the restriction of the one-parameter groups $t \longmapsto V(t, \overline{A})$ to $\mathcal{D}_\infty$ are strongly continuous one-parameter groups of $\tau_\infty$-continuous operators. By (2.5.1), these restricted groups are locally equicontinuous, so each $t \longmapsto V(t, \overline{A})$ is a $\Gamma_{A, \infty}$-group on $\mathcal{D}_\infty$, where $A \in \mathcal{S}$. The same conclusions are true for the space $\overline{\mathcal{D}_\mathcal{S}}^\infty$, defined as the closure of $\mathcal{D}_\mathcal{S}$ inside the complete Hausdorff locally convex space $(\mathcal{X}_\infty(\overline{\mathcal{S}}), \tau_\infty)$. Even more, $\overline{\mathcal{D}_\mathcal{S}}^\infty$ is a C$^\infty$-domain for $\overline{\mathcal{S}}$.}
\item \textit{The Lie algebra $\mathcal{L}_\infty := \left\{\left. \overline{B} \right|_{\overline{\mathcal{D}_\mathcal{S}}^\infty}: B \in \mathcal{L}\right\}$ is isomorphic to $\mathcal{L}$.}

\end{enumerate}

\textbf{Proof of Theorem 2.5:} \textbf{Proof of (1):} To see that $\mathcal{D}_\mathcal{S}$ is a C$^\infty$-domain, define inductively a non-decreasing sequence of subspaces by $$D_\mathcal{S}^0 := \mathcal{D}, \qquad \mathcal{D}_\mathcal{S}^{n + 1} := \text{span}_\mathbb{C} \left\{V(t, \overline{A})x, \, t \in \mathbb{R}, \, A \in \mathcal{S}, \, x \in \mathcal{D}_\mathcal{S}^n\right\}, \text{ for } n \geq 0,$$ and $$\mathcal{D}_\mathcal{S}^\infty := \bigcup_{n = 0}^{+ \infty} \mathcal{D}_\mathcal{S}^n.$$ Since each $V(t, \overline{A})$ sends $\mathcal{D}_\mathcal{S}^n$ into $\mathcal{D}_\mathcal{S}^{n + 1}$, for all $n \in \mathbb{N}$, it follows that each $V(t, \overline{A})$ sends $\mathcal{D}_\mathcal{S}^\infty$ into itself. Besides that, $\mathcal{D} \subseteq \mathcal{D}_\mathcal{S}^\infty$, so $\mathcal{D}_\mathcal{S} \subseteq \mathcal{D}_\mathcal{S}^\infty$. Now, $\mathcal{D} \subseteq \mathcal{D}_\mathcal{S}$ implies $$\left[\sum_{j = 1}^l c_j \, V(t_j, \overline{A_j})x_j\right] \in \mathcal{D}_\mathcal{S},$$ for every $l \in \mathbb{N}$, $c_j \in \mathbb{C}$, $t_j \in \mathbb{R}$, $A_j \in \mathcal{S}$ and $x_j \in \mathcal{D}$, proving that $\mathcal{D}_\mathcal{S}^1 \subseteq \mathcal{D}_\mathcal{S}$. Applying this argument inductively one sees that $\mathcal{D}_\mathcal{S}^n \subseteq \mathcal{D}_\mathcal{S}$, for all $n \in \mathbb{N}$, proving that $\mathcal{D}_\mathcal{S}^\infty \subseteq \mathcal{D}_\mathcal{S}$. Hence, $\mathcal{D}_\mathcal{S} = \mathcal{D}_\mathcal{S}^\infty$.

By hypothesis, $\mathcal{D}_\mathcal{S}^0 = \mathcal{D}$ is a C$^\infty$ domain for $\overline{\mathcal{S}}$. Suppose, for some $n \geq 0$, that $\mathcal{D}_\mathcal{S}^n$ is a C$^\infty$ domain for $\overline{\mathcal{S}}$. Since $\mathcal{D}_\mathcal{S} \subseteq \mathcal{D}_1$, equation (2.4.1), proved in Theorem 2.4, shows that $$\tilde{B}^{(1)} \, V_1(t, A_1)(x) = V(t, \overline{A}) \, [\exp (-t \, \text{ad }A)(B)]\, \tilde{} \, ^{(1)}(x),$$ for all closable $B \in \mathcal{L}$, $A \in \mathcal{S}$, $t \in \mathbb{R}$ and $x \in \mathcal{D}_\mathcal{S}$. Let $y = V(t, \overline{A})x$ be an element in the spanning set of $\mathcal{D}_\mathcal{S}^{n + 1}$, so that $t \in \mathbb{R}$, $A \in \mathcal{S}$ and $x \in \mathcal{D}_\mathcal{S}^n$. Then, since $\tilde{B}^{(1)} = \overline{B}|_{\mathcal{D}_1}$ and $V_1(t, A_1) = V(t, \overline{A})|_{\mathcal{D}_1}$, for all $t \in \mathbb{R}$, $$\overline{B}(y) = \overline{B} \, V(t, \overline{A})(x) = \tilde{B}^{(1)} \, V_1(t, A_1)(x) = V(t, \overline{A}) \, [\exp (-t \, \text{ad }A)(B)]\, \tilde{} \, ^{(1)}(x).$$ To see that $[\exp (-t \, \text{ad }A)(B)]\, \tilde{} \, ^{(1)}(x) \in \mathcal{D}_\mathcal{S}^n$, note that every $C \in \mathcal{L}$, with $C = \sum_{k = 1}^d c_k \, B_k$, satisfies $$\tilde{C}^{(1)} \subset \sum_{k = 1}^d c_k \, \overline{B_k}.$$ Indeed, if $z \in \mathcal{D}_1$ there exists a net $\left\{z_\alpha\right\}$ in $\mathcal{D}$ which converges to $z$ with respect to the $\tau_1$-topology, so $C(z_\alpha) \xrightarrow{\alpha} \tilde{C}^{(1)}(z)$ in the $\tau$-topology. But, by the definition of $\tau_1$, $B_k(z_\alpha) \xrightarrow{\alpha} \overline{B_k}(z)$, for all $1 \leq k \leq d$. Consequently, $$C(z_\alpha) = \sum_{k = 1}^d c_k \, B_k(z_\alpha) \xrightarrow{\alpha} \sum_{k = 1}^d c_k \, \overline{B_k}(z),$$ proving that $\tilde{C}^{(1)} = \sum_{k = 1}^d c_k \, \overline{B_k}|_{\mathcal{D}_1}$, as desired. Applying this reasoning to the operator $C = \exp (-t \, \text{ad }A)(B) \in \mathcal{L}$ and taking into account the hypothesis of $\mathcal{D}_\mathcal{S}^n$ being a C$^\infty$ domain for $\overline{\mathcal{S}}$, the conclusion $[\exp (-t \, \text{ad }A)(B)]\, \tilde{} \, ^{(1)}(x) \in \mathcal{D}_\mathcal{S}^n$ follows. Hence, $\overline{B}(y) \in \mathcal{D}_\mathcal{S}^{n + 1}$ which, by linearity of $\overline{B}$, proves that $$\overline{B}[\mathcal{D}_\mathcal{S}^{n + 1}] \subseteq \mathcal{D}_\mathcal{S}^{n + 1}.$$ Since $B$ was an arbitrary closable operator of $\mathcal{L}$, this establishes that each $\mathcal{D}_\mathcal{S}^n$ is a C$^\infty$ domain for $\overline{\mathcal{S}}$. Therefore, $\mathcal{D}_\mathcal{S}$, being a union of C$^\infty$ domains, is itself a C$^\infty$ domain. Since $\mathcal{X}_\infty(\overline{\mathcal{S}})$ is the largest C$^\infty$ domain for $\overline{\mathcal{S}}$, the inclusion $$\mathcal{D}_\mathcal{S} \subseteq \mathcal{X}_\infty(\overline{\mathcal{S}})$$ must hold, finishing the proof of (1).\\

\textbf{Proofs of (2) and (3):} To prove that $$(2.5.1) \qquad \rho_{p, n}(V(t, \overline{A})x) \leq \exp (n \|\text{ad }A\| \, |t|) \, \rho_{p, n}(x)$$ holds for every $A \in \mathcal{S}$, $p \in \Gamma_A$, $n \in \mathbb{N}$, $t \in \mathbb{R}$ and $x \in \mathcal{D}_\mathcal{S}$, the idea is to proceed by induction. Note that the evaluation of the seminorms $\rho_{p, n}$ is legitimate, for all $n \in \mathbb{N}$, because $\mathcal{D}_\mathcal{S}$ is contained in the complete Hausdorff vector space $$\bigcap_{n = 1}^{+ \infty} \mathcal{X}_n =: \mathcal{X}_\infty = \mathcal{X}_\infty(\overline{\mathcal{S}})$$ introduced at the beginning of this subsection. The case $n = 0$ is immediate since, for each $A \in \mathcal{S}$, $t \longmapsto V(t, \overline{A})$ is a $\Gamma_A$-isometrically equicontinuous one-parameter group. Suppose that $$\rho_{p, n}(V(t, \overline{A})x) \leq \exp (n \|\text{ad }A\| \, |t|) \, \rho_{p, n}(x)$$ for a fixed $n \in \mathbb{N}$ and all $A \in \mathcal{S}$, $p \in \Gamma_A$, $x \in \mathcal{D}_\mathcal{S}$ and $t \in \mathbb{R}$. First, note that for a fixed $C \in \mathcal{L}$, with $C = \sum_{k = 1}^d c_k \, B_k$, one has that $$\rho_{p, n}(\tilde{C}^{(1)}(x)) \leq |C|_1 \, \rho_{p, n + 1}(x),$$ for all $x \in \mathcal{D}$: $\rho_{p, n}(\tilde{C}^{(1)}(x)) = \rho_{p, n}(C(x)) = p(B^\mathsf{u}(C(x)))$, where $B^\mathsf{u}$ is a certain monomial of size $n$ in the elements $\left(B_k\right)_{1 \leq k \leq d}$, by the definition of $\rho_{p, n}$. Then, $$\rho_{p, n}(\tilde{C}^{(1)}(x)) = p \left(B^\mathsf{u} \left( \sum_{k = 1}^d c_k \, B_k \right)(x) \right) \leq \left( \sum_{k = 1}^d |b_k| \right) \rho_{p, n + 1}(x) = |C|_1 \, \rho_{p, n + 1}(x),$$ establishing the desired inequality. Hence, for $A \in \mathcal{S}$, $p \in \Gamma_A$, $x \in \mathcal{D}$, $t \in \mathbb{R}$ and $1 \leq k \leq d$, $$\rho_{p, n}(\overline{B_k} \, V(t, \overline{A})x) = \rho_{p, n}(\tilde{B_k} \, V(t, \overline{A})x) = \rho_{p, n}(V(t, \overline{A}) \, [\exp (-t \, \text{ad }A)(B_k)]\, \tilde{} \, ^{(1)}x)$$ $$\leq \exp (n \|\text{ad }A\| \, |t|) \, \|\exp (-t \, \text{ad }A)\| \, |B_k|_1 \, \rho_{p, n + 1}(x)$$ $$\leq \exp (n \|\text{ad }A\| \, |t|) \, \exp (|t| \, \|\text{ad }A\|) \, \rho_{p, n + 1}(x) = \exp ((n + 1) \|\text{ad }A\| \, |t|) \, \rho_{p, n + 1}(x).$$ Upon taking the maximum over $0 \leq k \leq d$ - the estimate for $k = 0$ is trivial - the inequality $$(2.5.1)' \qquad \rho_{p, n + 1}(V(t, \overline{A})x) \leq \exp ((n + 1) \|\text{ad }A\| |t|) \, \rho_{p, n + 1}(x)$$ follows. Now, note that the inequality just obtained implies that $V(t, \overline{A})$ leaves $\mathcal{D}_{n + 1}$ invariant: indeed, fix $x \in \mathcal{D}_{n + 1}$ and let $\left\{x_\alpha\right\}_{\alpha \in \mathcal{A}}$ be a net in $\mathcal{D}$ which $\tau_{n + 1}$-converges to $x \in \mathcal{D}_{n + 1}$. Then, by the inequality above, $$\rho_{p, n + 1}(V(t, \overline{A})(x_\alpha - x_{\alpha'})) \leq \exp ((n + 1) \|\text{ad }A\| |t|) \, \rho_{p, n + 1}(x_\alpha - x_{\alpha'}),$$ for all $\alpha, \alpha' \in \mathcal{A}$. This implies $\left\{V(t, \overline{A})(x_\alpha)\right\}_{\alpha \in \mathcal{A}}$ is a $\tau_{n + 1}$-Cauchy net and, since $\mathcal{D}_{n + 1}$ is complete, it follows this net is $\tau_{n + 1}$-convergent to an element in $\mathcal{D}_{n + 1}$. Since the $\tau_{n + 1}$-topology is stronger that the $\tau$-topology, this element must be $V(t, \overline{A})x$. Therefore, $V(t, \overline{A})x \in \mathcal{D}_{n + 1}$, for all $A \in \mathcal{S}$. By the definition of $\mathcal{D}_\mathcal{S}$, $$\mathcal{D}_\mathcal{S} \subseteq \mathcal{D}_{n + 1}.$$ The argument just exposed proves (2.5.1)' not only for $x \in \mathcal{D}$, but establishes it for all $x \in \mathcal{D}_{n + 1}$. In particular, it is valid for all $x \in \mathcal{D}_\mathcal{S}$, finishing the induction proof. Hence, (2.5.1) is true and, by the proof just developed, it follows that $V(t, \overline{A})$ leaves $\mathcal{D}_n$ invariant, for all $A \in \mathcal{S}$, $t \in \mathbb{R}$, $n \geq 1$. Moreover, (2.5.1) and a $\tau_n$-density argument of $\mathcal{D}_\mathcal{S}$ in $\mathcal{D}_n$ proves (just as it was already noted during the induction proof) that each $V(t, \overline{A})$ is a $\tau_n$-continuous operator on $\mathcal{D}_n$ and that (2.5.1) extends to $x \in \mathcal{D}_n$, for each fixed $n \geq 1$. This is a fact which will be used in the next argument, to prove that the restrictions of the one-parameter groups $t \longmapsto V(t, \overline{A})$ to $\mathcal{D}_n$ are strongly continuous with respect to the topology $\tau_n$, for every $n \geq 1$.\\

The estimates $$\sum_{k = 0}^{+ \infty} \rho_{p, n} \left( \frac{(\text{ad }A)^k(B_j)}{k!} \, (-t)^k(x) - B_j(x) \right) \leq \sum_{k = 1}^{+ \infty} \frac{|(\text{ad })^k(B_j)|_1}{k!} \, |t|^k \, \rho_{p, n + 1}(x)$$ $$\leq \rho_{p, n + 1}(x) \, \sum_{k = 1}^{+ \infty} \frac{\|\text{ad }A\|^k}{k!} \, |t|^k, \qquad t \in \mathbb{R}, \, p \in \Gamma_A, \, x \in \mathcal{D}$$ show each $\exp (-t \, \text{ad }A)(B_j)(x)$ $\tau_n$-converges to $B_j(x)$, when $t \longrightarrow 0$, for all $n \in \mathbb{N}$ and $1 \leq j \leq d$. Hence, using formula (2.5.1), one sees that $V(t, \overline{A}) \, \exp (-t \, \text{ad }A)(B_j)(x)$ $\tau_n$-converges to $B_j(x)$, for all $A \in \mathcal{S}$, $n \in \mathbb{N}$ and $1 \leq j \leq d$. Therefore, an application of formula (2.4.1) shows that if $t \longmapsto V(t, \overline{A})x$ is $\tau_n$-continuous for some fixed $n \in \mathbb{N}$ and all $x \in \mathcal{D}$, then $t \longmapsto B_j \, V(t, \overline{A})x$ is also $\tau_n$-continuous for all $x \in \mathcal{D}$, where $A \in \mathcal{S}$ and $1 \leq j \leq d$. But $t \longmapsto V(t, \overline{A})x$ is $\tau_0$-continuous, for every $A \in \mathcal{S}$ and $x \in \mathcal{D}$, so an inductive argument on $n$, followed by a $\tau_n$-density argument of $\mathcal{D}$ in $\mathcal{D}_n$ combined with an $\epsilon/3$-argument, shows that the group $t \longmapsto V(t, \overline{A})$ restricts to an exponentially equicontinuous group on $\mathcal{D}_n$ with type $w_p \leq n \|\text{ad }A\|$, where $p \in \Gamma_A$ and $n \geq 1$, which ends the proof of (2). Also, this proves that the restriction of the groups $t \longmapsto V(t, \overline{A})$ to the complete Hausdorff locally convex space $\mathcal{D}_\infty$, equipped with the $\tau_\infty$-topology, are strongly continuous one-parameter groups of $\tau_\infty$-continuous operators. Consequently, for each fixed $A \in \mathcal{S}$, $t \longmapsto V(t, \overline{A})$ restricts to a $\Gamma_{A, \infty}$-group on $\mathcal{D}_\infty$, since local equicontinuity follows at once from estimates (2.5.1), after extending them to all $x \in \mathcal{D}_n$, for each fixed $n \geq 1$.\\

By the group invariance just mentioned, together with completeness of $(\mathcal{D}_\infty, \tau_\infty)$, the chain of inclusions $$\mathcal{D}_\mathcal{S} \subseteq \overline{\mathcal{D}_\mathcal{S}}^\infty \subseteq \mathcal{D}_\infty \subseteq \mathcal{X}_\infty(\overline{\mathcal{S}})$$ is true. But even more is true: estimates (2.5.1) show, together with a $\tau_\infty$-density argument, that the $\tau_\infty$-continuous operators $V(t, \overline{A})$, $A \in \mathcal{S}$, $t \in \mathbb{R}$, leave $\overline{\mathcal{D}_\mathcal{S}}^\infty$ invariant. This observation, combined with the strong continuity claim of the previous paragraph and the inclusion $\overline{\mathcal{D}_\mathcal{S}}^\infty \subseteq \mathcal{D}_\infty$, above, shows that each $t \longmapsto V(t, \overline{A})$ is a $\Gamma_{A, \infty}$-group on $\overline{\mathcal{D}_\mathcal{S}}^\infty$, where $A$ belongs to $\mathcal{S}$. The fact that $\overline{\mathcal{D}_\mathcal{S}}^\infty$ is a C$^\infty$ domain for $\overline{\mathcal{S}}$ is also a consequence of a $\tau_\infty$-density argument on the C$^\infty$ domain $\mathcal{D}_\mathcal{S}$, and establishes (3).\\

\textbf{Proof of (4):} To prove the isomorphism claim first note that, since $\overline{\mathcal{D}_\mathcal{S}}^\infty$ is a C$^\infty$ domain for $\overline{\mathcal{S}}$, compositions of operators in $\mathcal{L}_\infty$ are still operators in $\mathcal{L}_\infty$, so that it makes sense to consider a Lie algebra structure inside $\mathcal{L}_\infty$. To prove that the map $B \longmapsto \overline{B}|_{\overline{\mathcal{D}_\mathcal{S}}^\infty}$ from $\mathcal{L}$ to $\mathcal{L}_\infty$ is a Lie algebra isomorphism it will first be proved that it is linear and that it preserves commutators. To these ends first note that, for every closable $B \in \mathcal{L}$, $\overline{B}|_{\overline{\mathcal{D}_\mathcal{S}}^\infty} \subset \tilde{B}^{(1)} = \overline{B}|_{\mathcal{D}_1}$ so, if $B = \sum_{k = 1}^d b_k \, B_k$, then $$\overline{B}|_{\overline{\mathcal{D}_\mathcal{S}}^\infty} = \sum_{k = 1}^d b_k \, \overline{B_k}|_{\overline{\mathcal{D}_\mathcal{S}}^\infty},$$ by what was already noted at the very beginning of the proof. Therefore, if $C, D \in \mathcal{L}$ and $\mu \in \mathbb{C}$ with $C = \sum_{k = 1}^d c_k \, B_k$ and $D = \sum_{k = 1}^d d_k \, B_k$, then $$\overline{C + \mu \, D}|_{\overline{\mathcal{D}_\mathcal{S}}^\infty} = \sum_{k = 1}^d (c_k + \mu \, d_k) \, \overline{B_k}|_{\overline{\mathcal{D}_\mathcal{S}}^\infty}$$ and, since $\overline{C}|_{\overline{\mathcal{D}_\mathcal{S}}^\infty} = \sum_{k = 1}^d c_k \, \overline{B_k}|_{\overline{\mathcal{D}_\mathcal{S}}^\infty}$ and $\overline{D}|_{\overline{\mathcal{D}_\mathcal{S}}^\infty} = \sum_{k = 1}^d d_k \, \overline{B_k}|_{\overline{\mathcal{D}_\mathcal{S}}^\infty}$, linearity follows. Preservation of commutators follows from a similar argument. That the map is a bijection is immediate, so $\mathcal{L}$ and $\mathcal{L}_\infty$ are indeed isomorphic. \hfill $\blacksquare$

\subsection*{Exponentiation Theorems}
\addcontentsline{toc}{subsection}{Exponentiation Theorems}

\subsubsection*{$\bullet$ The First Exponentiation Theorems}
\addcontentsline{toc}{subsubsection}{The First Exponentiation Theorems}

\indent

\textbf{Definition 2.6 (Exponentiable Lie Algebras of Linear Operators):} \textit{Let $\mathcal{X}$ be a Hausdorff locally convex space, $\mathcal{D}$ a dense subspace of $\mathcal{X}$ and $\mathcal{L} \subseteq \text{End}(\mathcal{D})$ a real finite-dimensional Lie algebra of linear operators. Then, $\mathcal{L}$ is said to be \textbf{exponentiable}, or to \textbf{exponentiate}, if every operator in $\mathcal{L}$ is closable and if there exists a simply connected Lie group $G$ having a Lie algebra $\mathfrak{g}$ isomorphic to $\mathcal{L}$, via an isomorphism $\eta \colon \mathfrak{g} \ni X \longmapsto \eta(X) \in \mathcal{L}$, and a strongly continuous locally equicontinuous representation $V \colon G \longrightarrow \mathcal{L}(\mathcal{X})$ such that $\mathcal{D} \subseteq C^\infty(V)$ and $dV(X) = \overline{\eta(X)}$, for all $X \in \mathfrak{g}$. If, moreover, every $dV(X) = \overline{\eta(X)}$, $X \in \mathfrak{g}$, is the generator of an equicontinuous one-parameter group (respectively, of a $\Gamma$-isometrically equicontinuous group, where $\Gamma$ is a fundamental system of seminorms for $\mathcal{X}$), it will be said that $\mathcal{L}$ exponentiates to a representation by equicontinuous one-parameter groups (respectively, exponentiates to a representation by $\Gamma$-isometrically equicontinuous one-parameter groups). If $\mathcal{X}$ is a Banach space, then an analogous terminology will be used: if $\mathcal{L}$ exponentiates and every $dV(X) = \overline{\eta(X)}$, $X \in \mathfrak{g}$, is the generator of a group of isometries, then it will be said that $\mathcal{L}$ exponentiates to a representation by isometries. Sometimes, $\mathfrak{g}$ and $\mathcal{L}$ will be identified, and the $\eta$ symbol will be omitted.}\\

A few words about Definition 2.6 need to be said. In finite-dimensional representation theory - meaning $\mathcal{X}$ is finite-dimensional - every real finite-dimensional Lie algebra of linear operators is exponentiable. One way to arrive at this conclusion is to first consider a simply connected Lie group $G$ having a Lie algebra isomorphic to $\mathcal{L}$ ($G$ is unique, up to isomorphism), which always exists, by a theorem - this follows, for example, from Lie's Third Theorem - see \cite[Theorem 9.4.11, page 334]{neeb} or \cite[Remark 9.5.12, page 341]{neeb} - combined with \cite[Proposition 9.5.1, page 335]{neeb}. Then, using the Baker-Campbell-Hausdorff Formula, one manages to construct a local homomorphism\footnote{$\phi$ is said to be a local homomorphism from $G$ to $\mathcal{L}(\mathcal{X})$ if there exists an open symmetric connected neighborhood of the identity $\Omega$ such that $\phi \colon \Omega \longrightarrow \mathcal{L}(\mathcal{X})$ is a map satisfying $\phi(gh) = \phi(g) \phi(h)$, for all $g, h$ in $\Omega$ for which $gh$ also belongs to $\Omega$.} from $G$ to $\mathcal{L}(\mathcal{X})$ which, by simple connectedness of $G$, is always extendable to a global homomorphism $\tilde{\phi} \colon G \longrightarrow \mathcal{L}(\mathcal{X})$, by the so-called Monodromy Principle (see \cite[Proposition 9.5.8, page 339]{neeb} and \cite[Theorem 1.4.5, page 50]{hossein}). Proceeding in this way, a group homomorphism $\tilde{\phi} \colon G \longrightarrow \mathcal{L}(\mathcal{X})$ with the desired properties may be finally obtained \cite[Theorem 9.5.9, page 340]{neeb}. For a brief discussion of exponentiability on finite-dimensional spaces, see also \cite[page 271]{jorgensen}.

To obtain the first exponentiation theorem of this subsection (Theorem 2.7), which is an ``equicontinuous locally convex version'' of \cite[Theorem 9.2, page 197]{jorgensenmoore}, another theorem of \cite{jorgensenmoore} is going to be needed:\\

\textbf{\cite[Theorem 9.1, page 196]{jorgensenmoore}:} \textit{Let $\mathcal{X}$ be a Hausdorff locally convex space,\footnote{Locally convex spaces are automatically assumed to be Hausdorff in \cite{jorgensenmoore} (see the second paragraph of Appendix A, of the reference).} $\mathcal{D}$ a dense subspace of $\mathcal{X}$ and $\mathcal{L} \subseteq \text{End}(\mathcal{D})$ a finite dimensional real Lie algebra of linear operators. Suppose $\mathcal{L}$ is generated, as a Lie algebra, by a finite set $\mathcal{S}$ of infinitesimal pregenerators of strongly continuous locally equicontinuous groups (those are abbreviated as cle groups, in \cite{jorgensenmoore} - see page 62; see also page 178 of this same reference), where $t \longmapsto V(t, \overline{A})$ denotes the group generated by $\overline{A}$, for all $A \in \mathcal{S}$. If the following two conditions are satisfied then, $\mathcal{L}$ exponentiates to a strongly continuous representation of a Lie group:}

\begin{enumerate}

\item \textit{the domain $\mathcal{D}$ is left invariant by the operators $\left\{V(t, \overline{A}): A \in \mathcal{S}, t \in \mathbb{R}\right\}$;}
\item \textit{for each $x \in \mathcal{D}$ and each pair of elements $A, B \in \mathcal{S}$ there is an open interval (which may depend on $x$, $A$ and $B$) such that the function $t \longmapsto B \, V(t, \overline{A})$ is bounded in $I$.}

\end{enumerate}

\textbf{Theorem 2.7:} \textit{Let $(\mathcal{X}, \tau)$ be a complete Hausdorff locally convex space, $\mathcal{D}$ a dense subspace of $\mathcal{X}$ and $\mathcal{L} \subseteq \text{End}(\mathcal{D})$ a finite dimensional \textbf{real} Lie algebra of linear operators. Assume the same hypotheses of Theorem 2.5. Then, $\mathcal{L}$ exponentiates.}\\

\textbf{Proof of Theorem 2.7:} Using the notations and the results of the previous subsection, consider the complete Hausdorff locally convex space $\overline{\mathcal{D}_\mathcal{S}}^\infty$, constructed in Theorem 2.5, which is a C$^\infty$ domain for $\mathcal{\overline{S}}$ that contains $\mathcal{D}$ and is left invariant by the operators $\left\{V(t, \overline{A}): A \in \mathcal{S}, t \in \mathbb{R}\right\}$. The estimates $$p(B \, V(t, \overline{A})x) \leq |B|_1 \, \rho_{p, 1}(V(t, \overline{A})x),$$ for all $p \in \Gamma_A$, $x \in \mathcal{D}$ and $A, B \in \mathcal{S}$, follow from the usual arguments and extend by $\tau_1$-density to $$p(\overline{B}|_{\overline{\mathcal{D}_\mathcal{S}}^\infty} \, V(t, \overline{A})x) \leq |B|_1 \, \rho_{p, 1}(V(t, \overline{A})x),$$ for all $x \in \overline{\mathcal{D}_\mathcal{S}}^\infty$, by $\tau_1$-continuity of $\left. V(t, \overline{A}) \right|_{\mathcal{D}_1}$, $(\tau_1 \times \tau_1)$-closedness of $\overline{B}^{\tau_1}$ and $\overline{B}|_{\overline{\mathcal{D}_\mathcal{S}}^\infty} \subset \overline{B}^{\tau_1} = \overline{B}|_{\mathcal{D}_1}$. This proves the existence of a non-empty open interval on which the function $$I \ni t \longmapsto \overline{B}|_{\overline{\mathcal{D}_\mathcal{S}}^\infty} \, V(t, \overline{A})x$$ is bounded. Also, $$\overline{A} = \overline{\overline{A}|_{\overline{\mathcal{D}_\mathcal{S}}^\infty}}, \qquad A \in \mathcal{S},$$ so each $\overline{A}|_{\overline{\mathcal{D}_\mathcal{S}}^\infty}$ is a pregenerator of an equicontinuous group on $\mathcal{X}$. Applying \cite[Theorem 9.1, page 196]{jorgensenmoore}, above, with $\mathcal{D}$ replaced by $\overline{\mathcal{D}_\mathcal{S}}^\infty$, yields that $\mathcal{L}_\infty$ is exponentiable. Since item 4 of Theorem 2.5 established that $\mathcal{L}_\infty$ and $\mathcal{L}$ are isomorphic, the result follows. \hfill $\blacksquare$ \\

Since every basis of a finite-dimensional Lie algebra $\mathcal{L}$ is a set of generators for $\mathcal{L}$, the following corollary is also true:\\

\textbf{Corollary 2.8:} \textit{Assume the same hypotheses of Theorem 2.7 are valid, with $\mathcal{S}$ replaced by $\left(B_k\right)_{1 \leq k \leq d}$. Then, $\mathcal{L}$ exponentiates. If the stronger hypothesis that there exists a fundamental system of seminorms $\Gamma$ for $\mathcal{X}$ such that the operators $B_k$ are all $\Gamma$-conservative (that is, they are conservative with respect to the \textbf{same} $\Gamma$), then $\mathcal{L}$ exponentiates to a representation by $\Gamma$-isometrically equicontinuous one-parameter groups.}\\

\textbf{Proof of Corollary 2.8:} Let $G \ni g \longmapsto V(g) \in \mathcal{L}(\mathcal{X})$ be the underlying Lie group representation, $\mathfrak{g}$ the Lie algebra of $G$ and $\eta \colon \mathfrak{g} \longrightarrow \mathcal{L} \subseteq \text{End}(\mathcal{D})$ be the Lie algebra representation as in Definition 2.6. Suppose that the stronger hypotheses of Corollary 2.8 are valid. Then, formula (7) at \cite[page 248]{yosida} says that $$V(\exp t \, \eta^{-1}(B_k))x = \lim_{n \rightarrow + \infty} \text{exp}\left(t \, \overline{B_k} \left(I - \frac{1}{n} \, \overline{B_k}\right)^{-1}\right)x, \qquad t \in [0, + \infty), \, x \in \mathcal{X}.$$ Then, a repetition of the argument in the last paragraph of the proof of Theorem 1.6.3 shows that each $t \longmapsto V(\exp t \, \eta^{-1}(B_k))$ is a $\Gamma$-isometrically equicontinuous group. Now, there exist $d$ real-valued analytic functions $\left\{t_k\right\}_{1 \leq k \leq d}$ defined on a neighborhood $\Omega$ of the identity of $G$ such that $g \longmapsto (t_k(g))_{1 \leq k \leq d}$ maps $\Omega$ diffeomorphically onto a neighborhood of the origin of $\mathbb{R}^d$, with $$g = \exp (t_1(g) \, \eta^{-1}(B_1)) \ldots \exp (t_k(g) \, \eta^{-1}(B_k)) \ldots \exp (t_d(g) \, \eta^{-1}(B_d)), \qquad g \in \Omega.$$ Therefore, connectedness of $G$ guarantees that the group generated by $\Omega$ is all of $G$, so every $g \in G$ can be written as $g = \prod_{j = 1}^{n(g)} \exp s_j \, \eta^{-1}(B_{k_j})$, where $s_j \in \mathbb{R}$ and $n(g) \in \mathbb{N}$. Therefore, the group property of $V$ gives $p(V(g)x) = p(x)$, for all $p \in \Gamma$, $g \in G$ and $x \in \mathcal{X}$ so, in particular, $t \longmapsto V(\exp t X)$ is a $\Gamma$-isometrically equicontinuous group, for all $X \in \mathfrak{g}$. This shows that every $dV(X) = \overline{\eta(X)}$, $X \in \mathfrak{g}$, is the generator of a $\Gamma$-isometrically equicontinuous group, so $\mathcal{L}$ exponentiates to a representation by $\Gamma$-isometrically equicontinuous one-parameter groups. \hfill $\blacksquare$\\

The next theorem is a version of \cite[Theorem 3.1]{jorgensengoodman} for pregenerators of equicontinuous groups on complete Hausdorff locally convex spaces. It is a $d$-dimensional noncommutative version of Lemma 1.6.3, in the sense that it guarantees exponentiability of a $d$-dimensional and (possibly) noncommutative Lie algebra of linear operators, instead of a 1-dimensional one:\\

\textbf{Theorem 2.9:} \textit{Let $(\mathcal{X}, \tau)$ be a complete Hausdorff locally convex space, $\mathcal{D}$ a dense subspace of $\mathcal{X}$ and $\mathcal{L} \subseteq \text{End}(\mathcal{D})$ a finite dimensional real Lie algebra of linear operators defined on $\mathcal{D}$. Assume that the following hypotheses are valid:}

\begin{enumerate}

\item \textit{$\mathcal{L}$ is generated, as a Lie algebra, by a finite set $\mathcal{S}$ of infinitesimal pregenerators of equicontinuous groups. For each $A \in \mathcal{S}$, let $\Gamma_A$ be the fundamental system of seminorms for $\mathcal{X}$ constructed in Subsection 1.5, with respect to which the operator $\overline{A}$ has the (KIP), is $\Gamma_A$-conservative and $V( \, \cdot \, , \overline{A})$ is $\Gamma_A$-isometrically equicontinuous;}
\item \textit{$\mathcal{L}$ possesses a basis $\left(B_k\right)_{1 \leq k \leq d}$ formed by closable elements \textbf{having the (KIP) with respect to} $$\Gamma_\mathcal{S} := \bigcup_{A \in \mathcal{S}} \Gamma_A;$$}
\item \textit{for each fixed $A \in \mathcal{S}$, every element of $\mathcal{D}$ is a $\tau$-projective analytic vector for $A$.}

\end{enumerate}

\textit{Then, $\mathcal{L}$ exponentiates.}\\

\textbf{Proof of Theorem 2.9:} Fix $A \in \mathcal{S}$. The proof will be divided in steps:

\begin{enumerate}[label=\textbf{\arabic*}]

\item \textbf{Each $\bm{x \in \mathcal{D}}$ is a $\bm{\tau_1}$-projective analytic vector for $\bm{A}$:} fix $p \in \Gamma_A$ and $x \in \mathcal{D}$. Since $B_k(x)$ is a $\tau$-projective analytic vector for $A$, for all $1 \leq k \leq d$ (remember $B_0 := I$), there exist $C_{x, p}, s_{x, p} > 0$ such that $$p(A^n \, B_k (x)) \leq C_{x, p} \, s_{x, p}^n \, n!, \qquad n \in \mathbb{N}, \, 1 \leq k \leq d.$$ In fact, if for each fixed $0 \leq k \leq d$, a certain $r_{x, p}(k) > 0$ satisfying $$\sum_{j = 0}^{+ \infty} \frac{p(A^j \, B_k (x))}{j!} \, r_{x, p}(k)^j < \infty$$ is chosen, it suffices to make $$C_{x, p} := \max \left\{\sum_{j = 0}^{+ \infty} \frac{p(A^j \, B_k (x))}{j!} \, r_{x, p}(k)^j, 0 \leq k \leq d\right\}$$ and $$s_{x, p} := \max \left\{\frac{1}{r_{x, p}(k)}: 0 \leq k \leq d\right\}.$$ Hence, if $C = \sum_{j = 1}^d c_j \, B_j$, then $$p(A^n(C(x))) \leq \sum_{k = 1}^d |c_k| \, p(A^n \, B_k(x)) \leq |C|_1 \, C_{x, p} \, s_{x, p}^n \, n!,$$ for all $n \in \mathbb{N}$. This implies $$p(A^n \, (\text{ad }A)^m(B)(x)) \leq |(\text{ad }A)^m(B)|_1 \, C_{x, p} \, s_{x, p}^n \, n! \leq \|\text{ad }A\|^m \, |B|_1 \, C_{x, p} \, s_{x, p}^n \, n!,$$ for all $m, n \in \mathbb{N}$ and $B \in \mathcal{L}$. One can prove by induction that $$B \, A^n(x) = \sum_{j = 0}^n \binom{n}{j} \, (-1)^{n - j} \, A^j \, (\text{ad }A)^{n - j}(B)(x), \qquad n \in \mathbb{N}, \, B \in \mathcal{L},$$ so that $$p(B \, A^n (x)) \leq \sum_{j = 0}^n \binom{n}{j} \, p(A^j \, (\text{ad }A)^{n - j}(B) (x))$$ $$\leq \sum_{j = 0}^n \binom{n}{j} \, \|\text{ad }A\|^{n - j} \, |B|_1 \, C_{x, p} \, s_{x, p}^j \, j!, \qquad p \in \Gamma_A, \, x \in \mathcal{D}.$$ Applying the inequality just obtained, it follows that for all $z \in \mathbb{C}$, $$\sum_{n = 0}^{+ \infty} \frac{p(B \, A^n(x))}{n!} \, |z|^n  \leq |B|_1 \, C_{x, p} \, \sum_{n = 0}^{+ \infty} |z|^n \sum_{k = 0}^{n} \binom{n}{k} \, \|\text{ad }A\|^{n - k} \, s_{x, p}^k \, \frac{k!}{n!}$$ $$\leq |B|_1 \, C_{x, p} \, \sum_{n = 0}^{+ \infty} |z|^n \sum_{k = 0}^{n} \binom{n}{k} \, \|\text{ad }A\|^{n - k} \, s_{x, p}^k = |B|_1 \, C_{x, p} \, \sum_{n = 0}^{+ \infty} (\|\text{ad }A\| + s_{x, p})^n \, |z|^n,$$ and the latter series converges for all $z \in \mathbb{C}$ satisfying $(\|\text{ad }A\| + s_{x, p})|z| < 1$. In particular, substituting $B$ by $B_k$ and taking the maximum over $0 \leq k \leq d$ proves the desired analyticity result.

\item \textbf{For every $\bm{p \in \Gamma_A}$ and every $\bm{\lambda \in \mathbb{\textbf{C}}}$ satisfying $\bm{|\text{\textbf{Re} }\lambda| > \|\tilde{\text{\textbf{ad}}}\text{ }A\|}$, $$\bm{\rho_{p, 1}((\lambda I - A)(x)) \geq (|\text{\textbf{Re} }\lambda| - \|\tilde{\text{\textbf{ad}}}\text{ }A\|) \, \rho_{p, 1}(x), \qquad x \in \mathcal{D}:}$$}By the choice of $\Gamma_A$, one sees that, for all $\lambda \in \mathbb{C} \backslash i \, \mathbb{R}$, the resolvent operator $\text{R}(\lambda, \overline{A})$ satisfies the bounds $$(2.9.1) \qquad p(\text{R}(\lambda, \overline{A})(x)) \leq \frac{1}{|\text{Re }\lambda|} \, p(x), \qquad p \in \Gamma_A, \, x \in \mathcal{X}.$$ Then, a repetition of the argument at the beginning of the proof of Theorem 2.4 gives the desired estimates.

\end{enumerate}

Now, since the basis elements possess the (KIP) with respect to $\Gamma_\mathcal{S}$, the operator $$A \colon \mathcal{D} \subseteq \mathcal{D}_1 \longrightarrow \mathcal{D}_1$$ is a $\tau_1$-densely defined $(\tau_1 \times \tau_1)$-closable linear operator on $\mathcal{D}_1$ which possesses the (KIP) with respect to $\Gamma_{A, 1} \subseteq \Gamma_\mathcal{S}$, by the same argument explained at the beginning of Theorem 2.4. Therefore, the induced linear operator $$A_{\rho_{p, 1}} \colon \pi_{\rho_{p, 1}}[\mathcal{D}] \longrightarrow (\mathcal{D}_1)_{\rho_{p, 1}}, \quad A_{\rho_{p, 1}}([x]_{\rho_{p, 1}}) = [\overline{A}^{\tau_1}(x)]_{\rho_{p, 1}},$$ is well-defined for all $\rho_{p, 1} \in \Gamma_{A, 1}$, and $\pi_{\rho_{p, 1}}[\mathcal{D}]$ is a $\|\, \cdot \,\|_{\rho_{p, 1}}$-dense set of analytic vectors for $A_{\rho_{p, 1}}$, by item 1, above. Furthermore, the inequality proved in item 2 implies $$\rho_{p, 1}((\lambda I - A)^n(x)) \geq (|\text{Re }\lambda| - \|\tilde{\text{ad}}\text{ }A\|)^n \, \rho_{p, 1}(x), \qquad |\text{Re }\lambda| > \|\tilde{\text{ad}}\text{ }A\|,$$ for all $\rho_{p, 1} \in \Gamma_{A, 1}$, $n \in \mathbb{N}$, $x \in \mathcal{D}$. Therefore, \cite[Theorem 1]{rusinek} becomes applicable, and the operators $A_{\rho_{p, 1}}$ are infinitesimal pregenerators of strongly continuous groups $t \longmapsto V(t, \overline{A_{\rho_{p, 1}}})$ on $(\mathcal{D}_1)_{\rho_{p, 1}}$ which satisfy $$\|V(t, \overline{A_{\rho_{p, 1}}})x_{\rho_{p, 1}}\|_{\rho_{p, 1}} \leq \text{exp }(\|\tilde{\text{ad}}\text{ }A\| \, |t|) \, \|x_{\rho_{p, 1}}\|_{\rho_{p, 1}}, \qquad x_{\rho_{p, 1}} \in (\mathcal{D}_1)_{\rho_{p, 1}}.$$ Hence, by a classical generation-type theorem for groups (see the Generation Theorem for Groups in \cite[page 79]{engel}, it follows that, for each $p \in \Gamma$ and every $\lambda \in \mathbb{C}$ satisfying $|\text{Re }\lambda| > \|\tilde{\text{ad}}\text{ }A\|$, the resolvent operators $\text{R}(\lambda, \overline{A_{\rho_{p, 1}}})$ are well-defined (in other words, $\lambda \in \rho(\overline{A_{\rho_{p, 1}}})$, for all $\rho_{p, 1} \in \Gamma_{A, 1}$ and all $\lambda \in \mathbb{C}$ satisfying $|\text{Re }\lambda| > \|\tilde{\text{ad}}\text{ }A\|$) and the ranges $$(\lambda I - A_{\rho_{p, 1}})[\pi_{\rho_{p, 1}}[\mathcal{D}]]$$ are dense in their respective spaces $(\mathcal{D}_1)_{\rho_{p, 1}}$. In particular, they are dense in $\pi_{\rho_{p, 1}}[\mathcal{D}]$, for each $\rho_{p, 1} \in \Gamma_{A, 1}$.\footnote{This follows from the fact that $\lambda I - \overline{A_{\rho_{p, 1}}}$ is surjective and from $\text{Ran}(\lambda I - \overline{A_{\rho_{p, 1}}}) = \overline{\text{Ran}(\lambda I - A_{\rho_{p, 1}})}$.} Just as it was argued at the beginning of the proof of Lemma 1.6.3 (remember $\Gamma_{A, 1}$ is a \underline{saturated} family of seminorms), one concludes $\text{Ran }(\lambda I - A)$ is $\tau_1$-dense in $\mathcal{D}$ if $|\text{Re }\lambda| > \|\tilde{\text{ad}}\text{ }A\|$, for all $A \in \mathcal{S}$. Therefore, since $A \in \mathcal{S}$ is arbitrary, Theorem 2.7 guarantees $\mathcal{L}$ is exponentiable. \hfill $\blacksquare$\\

Just as in the case of Corollary 2.8, the following is also true:\\

\textbf{Corollary 2.10:} \textit{Assume the same hypotheses of Theorem 2.9 are valid, with $\mathcal{S}$ replaced by $\left(B_k\right)_{1 \leq k \leq d}$. Then, $\mathcal{L}$ exponentiates. If the stronger hypothesis that there exists a fundamental system of seminorms $\Gamma$ for $\mathcal{X}$ such that the operators $B_k$ are all $\Gamma$-conservative (that is, they are conservative with respect to the \textbf{same} $\Gamma$), then $\mathcal{L}$ exponentiates to a representation by $\Gamma$-isometrically equicontinuous one-parameter groups.}

\subsubsection*{$\bullet$ Strongly Elliptic Operators - Sufficient Conditions for Exponentiation}
\addcontentsline{toc}{subsubsection}{Strongly Elliptic Operators - Sufficient Conditions for Exponentiation}

An important concept which will permeate the next theorems is that of a strongly elliptic operator.\\

\textbf{Definition (Strongly Elliptic Operators):} \textit{Let $\mathcal{X}$ be a complete Hausdorff locally convex space, $G$ a Lie group with Lie algebra $\mathfrak{g}$, with $\left(X_k\right)_{1 \leq k \leq d}$ being an ordered basis for $\mathfrak{g}$, and $V \colon G \longrightarrow \mathcal{L}(\mathcal{X})$ a strongly continuous locally equicontinuous representation, where $$\partial V \colon \mathfrak{g} \ni X \longmapsto \partial V(X) \in \text{End}(C^\infty(V))$$ is the induced Lie algebra representation. Assume, without loss of generality, that $\partial V$ is \textbf{faithful} (in other words, that it is an injective map). It is said that an operator $$H_m := \sum_{|\alpha| \leq m} c_\alpha \, \partial V(X^\alpha)$$ of order\footnote{Remember the notation introduced in Subsection 1.4.} $m$ belonging to the complexification of the enveloping algebra of $\mathfrak{g}$, $(\mathfrak{U}[\partial V[\mathfrak{g}]])_\mathbb{C} := \mathfrak{U}[\partial V[\mathfrak{g}]] + i \, \mathfrak{U}[\partial V[\mathfrak{g}]]$, is \textbf{elliptic} if the polynomial function $$P_m(\xi) := \sum_{|\alpha| = m} c_\alpha \, \xi^\alpha$$ on $\mathbb{R}^d$, called the principal part of $H_m$, satisfies $P_m(\xi) \neq 0$, if $\xi \neq 0$. This definition is basis independent - see \cite[page 21]{robinson}. Finally, an elliptic operator is said to be \textbf{strongly elliptic} if $$\text{Re }(-1)^{m/2} \sum_{|\alpha| = m} c_\alpha \, \xi^\alpha > 0,$$ for all $\xi \in \mathbb{R}^d \, \backslash \left\{0\right\}$ - see \cite[page 28]{robinson}. The prototypical example of a strongly elliptic operator is the ``minus Laplacian'' operator $$-\sum_{k = 1}^d [\partial V(X_k)]^2.$$ For a matter of terminological convenience, in the future (Theorem 2.14), the same definition which was just made for elliptic and strongly elliptic operators will be extended to (complexifications of) universal enveloping algebras of operators which are \textbf{not} necessarily associated with a strongly continuous Lie group representation.}\footnote{For an exposition on the usefulness of the concept of strong ellipticity in the realm of PDE's, see \cite{nirenberg}.}\\

Motivated by \cite{bratteliheat}, the following definition will be employed:\\

\textbf{Definition (Representation by Closed Linear Operators):} In the next theorems the following notations will be used: let $\mathcal{X}$ be a complete Hausdorff locally convex space, $\mathcal{D} \subseteq \mathcal{X}$ a dense subspace, $\mathfrak{g}$ a real finite-dimensional Lie algebra and $\eta$ a function defined on $\mathfrak{g}$ assuming values on (not necessarily continuous) closed linear operators on $\mathcal{X}$ - in other words, $\eta(X)$ is a closed linear operator on $\mathcal{X}$, for all $X \in \mathfrak{g}$. Assume, also, that $\mathcal{D}$ is a core for all the operators in $\eta[\mathfrak{g}]$, so that $\mathcal{D}$ will be called a \textbf{core domain} for $\eta[\mathfrak{g}]$, and that $\eta_\mathcal{D} \colon \mathfrak{g} \longrightarrow \mathcal{L} \subseteq \text{End}(\mathcal{D})$, defined by $\eta_\mathcal{D}(X) := \eta(X)|_\mathcal{D}$, is a faithful Lie algebra representation onto $\mathcal{L}$. Then, the triple $(\mathcal{D}, \mathfrak{g}, \eta)$ is said to be a \textbf{representation of $\bm{\mathfrak{g}}$ by closed linear operators on $\bm{\mathcal{X}}$}. It is said to \textbf{exponentiate}, or to be \textbf{exponentiable}, if there exists a simply connected Lie group $G$ having $\mathfrak{g}$ as its Lie algebra and a strongly continuous locally equicontinuous representation $V \colon G \longrightarrow \mathcal{L}(\mathcal{X})$ such that $\mathcal{D} \subseteq C^\infty(V)$ and $dV(X) = \overline{\eta(X)|_\mathcal{D}} = \eta(X)$, for all $X \in \mathfrak{g}$. Also, identical terminologies of Definition 2.6 will be used to indicate if the operators $dV(X)$, $X \in \mathfrak{g}$, are generators of equicontinuous, or even $\Gamma$-isometrically equicontinuous groups.\\

Theorem 2.11, below, generalizes \cite[Theorem 2.8]{bratteliheat} to higher-order strongly elliptic operators, but still in the Banach space context. Differently of what is done in \cite{bratteliheat}, this theorem will be proved for an arbitrary dense core domain $\mathcal{D}$ which is left invariant by the operators in $\mathcal{L}$ - in \cite{bratteliheat}, the claim is proved for the particular case in which the domain is $$\mathcal{X}_\infty := \bigcap_{n = 1}^{+ \infty} \bigcap \left\{\text{Dom }[\eta(X_{i_1}) \ldots \eta(X_{i_k}) \ldots \eta(X_{i_n})]: X_{i_k} \in \mathcal{B}, 1 \leq k \leq n \right\},$$ where $\mathcal{B}$ is a basis for $\mathfrak{g}$ (by what was proved at the beginning of this section, $\mathcal{X}_\infty$ is the maximal C$^\infty$ domain for the set $\left\{\eta(X_k)\right\}_{1 \leq k \leq d}$). The references \cite[Theorem 2.1]{bratteliheat} and \cite[Theorem 2.2, page 80; Lemma 2.3, page 82]{robinson} were vital sources of inspiration for the proof of Theorem 2.11.\footnote{See also \cite{langlandsthesis} and \cite{langlands}.} The theorems of the latter reference obtain some estimates which are fulfilled by strongly elliptic linear operators coming from a strongly continuous representation of a given Lie group and, together with Theorem 5.1 of this same book (page 30), prove (in particular) the following:\\

\textit{Let $\mathcal{X}$ be a Banach space, $G$ a Lie group with Lie algebra $\mathfrak{g}$, $V \colon G \longrightarrow \mathcal{L}(\mathcal{X})$ a strongly continuous representation of $G$ by bounded operators and $H_m \in \mathfrak{U}(\partial V[\mathfrak{g}])_\mathbb{C}$ an element of order $m$ which is strongly elliptic. Then, $-H_m$ is an infinitesimal pregenerator of a strongly continuous semigroup $t \longmapsto S(t)$ satisfying $S(t)[\mathcal{X}] \subseteq C^\infty(V)$ for $t \in (0, 1]$ and, for each $n \in \mathbb{N}$, there exists a constant $C_n > 0$ such that $$\rho_n(S(t)y) \leq C_n \, t^{-\frac{n}{m}} \|y\|, \qquad t \in (0, 1], \, y \in \mathcal{X}.$$ More precisely, such constants are of the form $C_n = K \, L^n n!$, for some $K, L \geq 1$.\footnote{In order to avoid unnecessary confusions, it should be noted that references \cite{bratteliheat}, \cite{brattelidissipative} and \cite{robinson} (see page 30) use a slightly different convention from the one employed in this manuscript to define infinitesimal generators: if an operator $H_m$ is strongly elliptic then, according to their definitions, $\overline{H_m}$ will be the generator of a one-parameter semigroup. However, using the definition of infinitesimal generators of this paper, if $H_m$ is strongly elliptic, then $- \overline{H_m}$ will be the generator of a one-parameter group. Therefore, for example, $-\sum_{k = 1}^d [\partial V(X_k)]^2$ is a strongly elliptic operator and $\overline{\sum_{k = 1}^d [\partial V(X_k)]^2}$ generates a one-parameter semigroup, with the definitions employed here.}}\\

Regarding the reciprocal question or, in other words, the issue of exponentiability, one has the following result:\\
 
\textbf{Theorem 2.11 (Exponentiation - Banach Space Version, Isometric Case):} \textit{Let $\mathcal{X}$ be a Banach space, $\mathfrak{g}$ a real finite-dimensional Lie algebra with an ordered basis $\left(X_k\right)_{1 \leq k \leq d}$, $(\mathcal{D}, \mathfrak{g}, \eta)$ a representation of $\mathfrak{g}$ by closed linear operators on $\mathcal{X}$, where $B_k := \eta(X_k)$, $1 \leq k \leq d$, and $H_m$ an element of order $m \geq 2$ of $\mathfrak{U}(\eta_\mathcal{D}[\mathfrak{g}])_\mathbb{C}$, where $\mathfrak{U}(\eta_\mathcal{D}[\mathfrak{g}])_\mathbb{C} := \mathfrak{U}(\eta_\mathcal{D}[\mathfrak{g}]) + i \, \mathfrak{U}(\eta_\mathcal{D}[\mathfrak{g}])$ is the complexification of $\mathfrak{U}(\eta_\mathcal{D}[\mathfrak{g}])$. Suppose that the following hypotheses are valid:}

\begin{enumerate}

\item \textit{$\eta(X)$ is a conservative operator, for every $X \in \mathfrak{g}$;}

\item \textit{$-H_m$ is an infinitesimal pregenerator of a strongly continuous semigroup $t \longmapsto S(t)$ satisfying $S(t)[\mathcal{X}] \subseteq \mathcal{D}$, for each $t > 0$ and, for each $n \in \mathbb{N}$ satisfying $0 < n \leq m - 1$, there exists $C_n > 0$ for which the estimates $$(2.11.1) \qquad \rho_n(S(t)y) \leq C_n \, t^{-\frac{n}{m}} \|y\|$$ are verified for all $t \in (0, 1]$ and all $y \in \mathcal{X}$, where $B_0 := I$ and $$\rho_n(x) := \max \left\{\|B_{i_1} \ldots B_{i_k} \ldots B_{i_n}x\|: 1 \leq k \leq n, 0 \leq i_k \leq d\right\},$$ for all $x \in \mathcal{D}$.}

\end{enumerate}

\textit{Then, $(\mathcal{D}, \mathfrak{g}, \eta)$ exponentiates to a representation by isometries: in other words, each $\eta(X)$, $X \in \mathfrak{g}$, is the generator of a group of isometries.}\\

\textbf{Proof of Theorem 2.11:} The first task of this proof is to obtain the estimates $$(2.11.2) \qquad \rho_n(y) \leq \epsilon^{m - n} \|H_m(y)\| + \frac{E_n}{\epsilon^n}\|y\|, \qquad y \in \mathcal{D},$$ valid for all $0 < \epsilon \leq 1$, $0 < n \leq m - 1$ and some $E_n > 0$.\\

The semigroup $t \longmapsto S(t)$ is strongly continuous, so there exist constants $M \geq 1$ and $w \geq 0$ such that $$\|S(t)y\| \leq M \exp (w t) \, \|y\|,$$ for all $y \in \mathcal{X}$ and $t \geq 0$. Now, if $0 < \delta < 1$ satisfies $0 \leq w < \frac{1}{\delta}$, the resolvent $\text{R}(1, - \delta \, \overline{H_m}) = (1 + \delta \, \overline{H_m})^{-1}$ is well-defined and $$\text{R}(1, - \delta \, \overline{H_m}) = \int_0^{+ \infty} \exp (-t) \, S_{\delta t}(y) \, dt, \qquad y \in \mathcal{X},$$ with $$\int_0^{+ \infty} \exp (-t) \, \|S_{\delta t}(y)\| \, dt < \infty.$$ Fix a monomial $B^\mathsf{u}$ of size $0 < n \leq m - 1$ in the operators $\left(B_k\right)_{1 \leq k \leq d}$ and $t > 0$. If $\delta t \leq 1$, then $$\|B^\mathsf{u} S_{\delta t}(y)\| \leq C_n \, (\delta t)^{-\frac{n}{m}} \|y\|, \qquad y \in \mathcal{X},$$ by (2.11.1); if $\delta t > 1$, then $$\|B^\mathsf{u} S_{\delta t}(y)\| = \|B^\mathsf{u} S_{1 + (\delta t - 1)}(y)\| \leq C_n \, M \exp (w (\delta t - 1)) \, \|y\|, \qquad y \in \mathcal{X}.$$ These inequalities and $$\frac{1}{\exp (t)} \leq \sqrt[m]{\frac{1}{\exp (t)}} \leq t^{-\frac{k}{m}} \sqrt[m]{k!},$$ with $k = m - n$ and $k = n$, respectively, allow one to obtain, for every $0 < \delta < 1/(2w)$, the estimates $$\int_0^{+ \infty} \exp (-t) \, \|B^\mathsf{u} S_{\delta t}(y)\| \, dt$$ $$\leq C_n \, \|y\| \left(\int_0^{1 / \delta} \frac{\exp (-t)}{(\delta t)^{\frac{n}{m}}} \, dt + M \exp (-w) \int_{1 / \delta}^{+ \infty} \exp ((w \delta - 1)t) \, dt\right)$$ $$= C_n \, \|y\| \left(\delta^{-\frac{n}{m}} \left. \left[\exp (-t) \, \frac{m}{m - n} \, t^\frac{m - n}{m}\right] \right|_{t = 0}^{t = 1 / \delta} + \delta^{-\frac{n}{m}} \, \frac{m}{m - n} \int_0^{1 / \delta} \exp (-t) \, t^\frac{m - n}{m} \, dt \right.$$ $$\left. + \, M \exp (-w) \, \frac{\exp (w - 1 / \delta)}{1 - w \delta} \right) \leq C_n \, \|y\| \left( \delta^{-\frac{n}{m}} \, \delta^{\frac{m - n}{m}} \, \sqrt[m]{(m - n)!} \, \frac{m}{m - n} \, \delta^\frac{n - m}{m} \right.$$ $$\left. + \, \delta^{-\frac{n}{m}} \, \frac{m}{m - n} \int_0^{+ \infty} \exp (-t) \, t^\frac{m - n}{m} \, dt + 2 \, M \, \delta^{\frac{n}{m}} \sqrt[m]{n!} \right)$$ $$\leq C_n \left( \sqrt[m]{(m - n)!} \, \frac{m}{m - n} + I_0 \, \frac{m}{m - n} + 2 \, M \sqrt[m]{n!} \right) \delta^{-\frac{n}{m}} \, \|y\|,$$ for all $y \in \mathcal{X}$, where $$I_0 := \int_0^{+ \infty} \exp (-t) \, t^\frac{m - n}{m} \, dt$$ (in the last inequality it was used that $\delta^{-\frac{n}{m}} > \delta^{\frac{n}{m}}$). Therefore, if $y$ is any vector in $\mathcal{X}$, then $$\int_0^{+ \infty} \exp (-t) \, \|B^\mathsf{u} S_{\delta t}(y)\| \, dt < \infty$$ and, in particular, the integral $\int_0^{+ \infty} \exp (-t) \, B^\mathsf{u} S_{\delta t}(y) \, dt$ defines an element in $\mathcal{X}$. By closedness of each operator $B_k$ it follows that $(1 + \delta \, \overline{H_m})^{-1}(y) \in \text{Dom }B^\mathsf{u}$ and $$B^\mathsf{u} \, (1 + \delta \, \overline{H_m})^{-1}(y) = \int_0^{+ \infty} \exp (-t) \, B^\mathsf{u} \, S_{\delta t}(y) \, dt.$$ Also, the estimates above show that there exists $E_n' > 0$ independent of $\delta$ such that $$\|B^\mathsf{u} \, (1 + \delta \, \overline{H_m})^{-1}(y)\| \leq E_n' \, \delta^{- \frac{n}{m}} \|y\|, \qquad y \in \mathcal{X}.$$ Hence, if $y$ is of the form $y = (1 + \delta \, \overline{H_m})(y') = (1 + \delta \, H_m)(y')$, $y' \in \mathcal{D}$, one obtains $$\|B^\mathsf{u}(y')\| = \|B^\mathsf{u} \, (1 + \delta \, \overline{H_m})^{-1}[(1 + \delta \, \overline{H_m})(y')]\| \leq E_n' \, \delta^{- \frac{n}{m}} \|(1 + \delta \, \overline{H_m})(y')\|$$ $$\leq E_n' \, \delta^{\frac{m - n}{m}} \|H_m(y')\| + E_n' \, \delta^{- \frac{n}{m}} \|y'\|.$$ Now, define $\epsilon_n(\delta) := (E_n')^{\frac{1}{m - n}} \delta^{\frac{1}{m}}$. If $\epsilon_n(\delta_0) > 1$, for some $\delta_0$, then (2.11.2) is proved since, by continuity, the range of the function $\epsilon_n$ must contain the interval $(0, 1]$ when $\delta$ runs through $(0, 1 / (2w))$, as a consequence of Weierstrass Intermediate Value Theorem. Then, defining $E_n := (E_n')^{\frac{m}{m - n}}$, one obtains $$E_n \, \epsilon_n(\delta)^{-n} = (E_n')^{\frac{m}{m - n}} \epsilon_n(\delta)^{-n} = E_n' \, \delta^{-\frac{n}{m}},$$ so $$\rho_n(y') \leq \epsilon^{m - n} \|H_m(y')\| + \frac{E_n}{\epsilon^n} \|y'\|, \qquad y' \in \mathcal{D}, \, 0 < \epsilon \leq 1,$$ as desired. However, if $\epsilon_n(\delta) < 1$, for all $\delta \in (0, 1 / (2w))$, then multiply the function $\epsilon_n$ by a real constant $a$ greater than 1 so that the range of $\epsilon_n' := a \, \epsilon_n$ contains 1. This reduces the problem to the previous case and finishes the first step of the proof.\\

The idea, now, is to divide the rest of the proof into three items:

\begin{itemize}

\item Given $q \in \mathbb{N}$, $q \geq 1$, there exists a strictly positive constant $K_q$ such that $$\rho_n(S(t)y) \leq K_q \sup_{0 \leq j \leq q} \|H_m^j(y)\|,$$ for all $y \in \mathcal{D}$, $(q - 1)(m - 1) < n \leq q (m - 1)$ and $t \in (0, 1]$. This implies, in particular, that given $q \in \mathbb{N}$, $q \geq 1$, and $y \in \mathcal{D}$, there exists $K_{y, q} > 0$ satisfying $\rho_n(S(t)y) \leq K_{y, q}$, for all $t \in (0, 1]$ and $(q - 1)(m - 1) < n \leq q (m - 1)$ - actually, this is the result that will be invoked, later:\\

The procedure is by induction on $q$. By what was proved earlier, for each fixed $0 < n \leq m - 1$ the estimates $$\rho_n(y) \leq \epsilon^{m - n} \|H_m(y)\| + \frac{E_n}{\epsilon^n} \|y\|, \qquad y \in \mathcal{D},$$ are valid for some $E_n > 0$ and all $0 < \epsilon \leq 1$. This shows the existence of a $K_1 > 0$ satisfying $$\rho_n(S(t)y) \leq K_1 \sup_{0 \leq j \leq 1} \|H_m^j(y)\|, \qquad y \in \mathcal{D}, \, 0 < n \leq m - 1, \, t \in (0, 1],$$ since $H_m S(t) = S(t) H_m$ on $\text{Dom }H_m$ and $\|S(t)\| \leq M e^{wt}$, for some $M \geq 1$, $w \geq 0$ and all $t \geq 0$. Now, suppose that, for some fixed $q \in \mathbb{N}$, $q \geq 1$, the inequality $$\rho_n(S(t)y) \leq K_q \sup_{0 \leq j \leq q} \|H_m^j(y)\|$$ is true, for all $y \in \mathcal{D}$, $(q - 1)(m - 1) < n \leq q (m - 1)$, $t \in (0, 1]$. Fix $y \in \mathcal{D}$. A monomial $B^\mathsf{u}$ of size $q (m - 1) < n \leq (q + 1)(m - 1)$ in the operators $\left(B_k\right)_{1 \leq k \leq d}$ may be decomposed as $B^\mathsf{u} = B^{\mathsf{u}_0} B^{\mathsf{u}'}$, with $|\mathsf{u}_0| = m - 1$ and $|\mathsf{u}'| = n - (m - 1)$. Using $(2.11.2)$, one obtains $$\|B^\mathsf{u} S(t)y\| = \|B^{\mathsf{u}_0}B^{\mathsf{u}'}S(t)y\| \leq \epsilon^{m - (m - 1)} \|H_m \, B^{\mathsf{u}'} S(t)y\| + \frac{E_{m - 1}}{\epsilon^{m - 1}} \|B^{\mathsf{u}'}S(t)y\|$$ $$= \epsilon \, \|H_m \, B^{\mathsf{u}'} S(t)y\| + \frac{E_{m - 1}}{\epsilon^{m - 1}} \|B^{\mathsf{u}'}S(t)y\|, \qquad 0 < t \leq 1, \, 0 < \epsilon \leq 1.$$ On the other hand, by the induction hypothesis, $$\|H_m \, B^{\mathsf{u}'} S(t)y\| \leq \|[\text{ad }(H_m)(B^{\mathsf{u}'})] \, S(t)y\| + \|B^{\mathsf{u}'} H_m \, S(t)y\|$$ $$\leq k \, (n - (m - 1)) \, m \, \rho_n(S(t)y) + K_q \sup_{0 \leq j \leq q} \|H_m^j H_m(y)\|, \qquad 0 < t \leq 1,$$ $k$ being a constant depending only on $d$, $H_m$ and on the numbers $c_{ij}$ which are defined by $[B_i, B_j] = \sum_{k = 1}^d c_{ij}^{(k)} B_k$, for all $1 \leq i \neq j \leq d$ (see Subsection 1.7). Hence, using the induction hypothesis again, one obtains $$\|B^\mathsf{u} S(t)y\| \leq \epsilon \, [k \, (n - (m - 1)) \, m \, \rho_n(S(t)y) + K_q \sup_{0 \leq j \leq q + 1} \|H_m^j(y)\|]$$ $$+ \, \frac{E_{m - 1}}{\epsilon^{m - 1}} \, K_q \sup_{0 \leq j \leq q} \|H_m^j(y)\|, \qquad 0 < t \leq 1, \, 0 < \epsilon \leq 1.$$ Choosing an $1 \geq \epsilon_0 > 0$ such that $\epsilon_0 \, k \, (n - (m - 1)) \, m < 1$, and taking the maximum over all monomials $B^\mathsf{u}$ of size $n$ gives $$\rho_n(S(t)y) \leq \frac{\epsilon_0 \, K_q \, \sup_{0 \leq j \leq q + 1} \|H_m^j(y)\| + \frac{E_{m - 1}}{\epsilon_0^{m - 1}} \, K_q \, \sup_{0 \leq j \leq q} \|H_m^j(y)\|}{1 - \epsilon_0 \, k \, (n - (m - 1))}$$ $$\leq \left[\frac{K_q \left(\epsilon_0 + \frac{E_{m - 1}}{\epsilon_0^{m - 1}}\right)}{1 - \epsilon_0 \, k \, (n - (m - 1))}\right] \sup_{0 \leq j \leq q + 1} \|H_m^j(y)\|, \qquad 0 < t \leq 1.$$ This concludes the induction proof.\\

\item There exist $K, L \geq 1$ such that, given $q \in \mathbb{N}$, $q \geq 1$ and $(q - 1)(m - 1) < n \leq q (m - 1)$, there exists a strictly positive constant $C_n$ defined by $C_n := K \, L^n n!$ satisfying $$\rho_n(S(t)y) \leq C_n \, t^{-\frac{n}{m}} \|y\|, \qquad y \in \mathcal{X}, \, t \in (0, 1].$$ More precisely, there exist $K, L \geq 1$ such that $$\rho_n(S(t)y) \leq K \, L^n n! \, t^{- \frac{n}{m}} \|y\|,$$ for all $n \in \mathbb{N}$, $n \geq 1$, $y \in \mathcal{X}$ and $t \in (0, 1]$:\\

The idea is, again, to proceed by induction. The case $q = 1$ follows from (2.11.1), with $K := D_1 := \max \left\{1, \max \left\{C_j: 0 < j \leq m - 1\right\} \right\}$ and $L := 1$, so it is already known to be true. Fix $q \geq 1$ and suppose that there exists, for every $p \leq q$ and $l \in \mathbb{N}$ satisfying $(p - 1)(m - 1) < l \leq p (m - 1)$, a strictly positive constant $C_l$ (not necessarily of the form $K \, L^l l!$, $K, L \geq 1$) such that $$\rho_l(S(t)y) \leq C_l \, t^{-\frac{l}{m}} \|y\|, \qquad y \in \mathcal{X}, \, t \in (0, 1].$$ Fix $n \in \mathbb{N}$ satisfying $q (m - 1) < n \leq (q + 1)(m - 1)$ (in particular, $n > m - 1$), $0 < t \leq 1$ and $\bm{y \in \mathcal{D}}$. As in the previous item, take a monomial $B^\mathsf{u}$ of size $q (m - 1) < n \leq (q + 1)(m - 1)$ in the operators $\left(B_k\right)_{1 \leq k \leq d}$ and decompose it as $B^\mathsf{u} = B^{\mathsf{u}_0} B^{\mathsf{u}'}$, with $|\mathsf{u}_0| = m - 1$ and $|\mathsf{u}'| = n - (m - 1)$. Then, by Lemma 1.7.1, $$(2.11.3) \qquad B^\mathsf{u} S(t)y = B^{\mathsf{u}_0} S_s \, B^{\mathsf{u}'} S_{t - s}(y) + B^{\mathsf{u}_0} [(\text{ad }B^{\mathsf{u}'})(S_s)] \, S_{t - s}(y)$$ $$= B^{\mathsf{u}_0} S_s \, B^{\mathsf{u}'} S_{t - s}(y) - \int_0^s B^{\mathsf{u}_0} S_r \, [(\text{ad }B^{\mathsf{u}'})(H_m)] \, S_{t - r}(y) \, dr,$$ for all $0 < s < t$. Like before, one has the inequality $$\|[(\text{ad }B^{\mathsf{u}'})(H_m)] \, S_{t - r}(y)\| \leq k \, (n - (m - 1)) \, m \, \rho_n(S_{t - r}(y))$$ where $k$, again, depends only on $d$, $H_m$ and on the numbers $c_{ij}$ defined by $[B_i, B_j] = \sum_{j = 1}^d c_{ij} \, B_j$, for all $1 \leq i \neq j \leq d$. Applying this inequality, together with the induction hypothesis, one obtains from $(2.11.3)$ the estimates $$\|B^\mathsf{u} \, S(t)y\| \leq C_{m - 1} \, s^{- \frac{m - 1}{m}} \|B^{\mathsf{u}'} S_{t - s}(y)\|$$ $$+ \int_0^s C_{m - 1} \, r^{- \frac{m - 1}{m}} \|[\text{ad }(B^{\mathsf{u}'})(H_m)] \, S_{t - r}(y)\| \, dr$$ $$\leq C_{m - 1} \, s^{- \frac{m - 1}{m}} C_{n - (m - 1)} \, (t - s)^{- \frac{n - (m - 1)}{m}} \|y\|$$ $$+ \, C_{m - 1} \, k \, (n - (m - 1)) \, m \int_0^s r^{- \frac{m - 1}{m}} \rho_n(S_{t - r}(y)) \, dr.$$ Making the changes of variables $s = \lambda t$, where $\lambda \in (0, 1)$, and $r = ut$, inside the integral, such inequality becomes $$\|B^\mathsf{u} S(t)y\| \leq C_{m - 1} \, C_{n - (m - 1)} \, t^{- \frac{n}{m}} \lambda^{- \frac{m - 1}{m}} (1 - \lambda)^{- \frac{n - (m - 1)}{m}} \|y\|$$ $$+ \, C_{m - 1} \, k \, (n - (m - 1)) \, m \, t^{\frac{1}{m}} \int_0^\lambda u^{- \frac{m - 1}{m}} \rho_n(S_{t (1 - u)}(y)) \, du.$$ Putting $\lambda = n^{-m}$, one has that $\lambda^{- \frac{m - 1}{m}} = n^{m - 1}$ and $$(1 - \lambda)^{- \frac{n - (m - 1)}{m}} = (1 - \lambda)^{- \frac{n}{m}} (1 - \lambda)^{\frac{m - 1}{m}} \leq (1 - \lambda)^{- \frac{n}{m}} \leq c_m,$$ where $$c_m := \sup_{n \geq m} \left\{(1 - n^{-m})^{- \frac{n}{m}}\right\} > 1.$$ Hence, by taking the maximum over all the monomials of size $n$, the above inequality becomes $$\rho_n(S(t)y) \leq C_{m - 1} \, C_{n - (m - 1)} \, t^{- \frac{n}{m}} \lambda^{- \frac{m - 1}{m}} (1 - \lambda)^{- \frac{n - (m - 1)}{m}} \|y\|$$ $$+ \, C_{m - 1} \, k \, (n - (m - 1)) \, m \, t^{\frac{1}{m}} \int_0^{\lambda} u^{- \frac{m - 1}{m}} \rho_n(S_{t (1 - u)}(y)) \, du$$ $$\leq C_{m - 1} \, C_{n - (m - 1)} \, t^{- \frac{n}{m}} n^{m - 1} c_m \, \|y\|$$ $$+ \, C_{m - 1} \, k \, (n - (m - 1)) \, m \, t^{\frac{1}{m}} \int_0^{n^{-m}} u^{- \frac{m - 1}{m}} \rho_n(S_{t (1 - u)}(y)) \, du.$$ But, since the $0 < t \leq 1$ fixed was arbitrary, one may repeat this procedure and iterate $j - 1$ times the last inequality to obtain $$(2.11.4) \qquad \rho_n(S(t)y) \leq C_{m - 1} \, C_{n - (m - 1)} \, t^{- \frac{n}{m}} n^{m - 1} c_m \, \|y\| \sum_{i = 0}^{j - 1} a_i + R_j,$$ where $a_0 := 1$, $$a_i := (C_{m - 1} \, k \, (n - (m - 1)) \, m \, t^{\frac{1}{m}})^i \int_0^{n^{-m}} u_1^{- \frac{m - 1}{m}} (1 - u_1)^{- \frac{n}{m}} \, du_1 \ldots$$ $$\ldots \int_0^{n^{-m}} u_i^{- \frac{m - 1}{m}} (1 - u_i)^{- \frac{n}{m}} \, du_i$$ $$\leq (C_{m - 1} \, k \, n \, m \, t^{\frac{1}{m}})^i \left[\sup_{n \geq m} (1 - n^{-m})^{- \frac{n}{m}} \right]^i \left(\int_0^{n^{-m}} u^{- \frac{m - 1}{m}} \, du \right)^i$$ $$= (C_{m - 1} \, k \, n \, m \, c_m \, t^{\frac{1}{m}})^i (m \, n^{-1})^i = (C_{m - 1} \, k \, m^2 \, c_m \, t^{\frac{1}{m}})^i,$$ if $i \geq 1$, and $$R_j := (C_{m - 1} \, k (n - (m - 1)) m t^{\frac{1}{m}})^j \int_0^{n^{-m}} u_1^{- \frac{m - 1}{m}} \ldots$$ $$\ldots \int_0^{n^{-m}} u_j^{- \frac{m - 1}{m}} \, \rho_n(S_{t (1 - u_1) \ldots (1 - u_j)}(y)) \, du_1 \ldots du_j,$$ for all $j \geq 1$. Therefore, by the previous item, $$R_j \leq (C_{m - 1} \, k \, m^2 \, t^{\frac{1}{m}})^j K_{y, q + 1}.$$ If $t_m > 0$ is chosen so that $C_{m - 1} \, k \, m^2 \, c_m \, (t_m)^{\frac{1}{m}} = \frac{1}{2}$, that is, if $$t_m = (2 \, C_{m - 1} \, k \, m^2 \, c_m)^{-m},$$ then for all $0 < t \leq t_m$, $$\sum_{k = 0}^{+ \infty} a_k \leq \sum_{k = 0}^{+ \infty} (C_{m - 1} \, k \, m^2 \, c_m \, t^{\frac{1}{m}})^k \leq \sum_{k = 0}^{+ \infty} (C_{m - 1} \, k \, m^2 \, c_m \, (t_m)^{\frac{1}{m}})^k = \sum_{k = 0}^{+ \infty} \left( \frac{1}{2} \right)^k = 2$$ and $$\lim_{k \rightarrow + \infty} R_k \leq\footnote{Here, it was used that $c_m > 1$.} \, K_{y, q + 1} \left[\lim_{k \rightarrow + \infty}(C_{m - 1} \, k \, m^2 \, c_m \, t^{\frac{1}{m}})^k\right]$$ $$\leq K_{y, q + 1} \left[\lim_{k \rightarrow + \infty}(C_{m - 1} \, k \, m^2 \, c_m \, (t_m)^{\frac{1}{m}})^k\right] = K_{y, q + 1} \left[\lim_{k \rightarrow + \infty} \left(\frac{1}{2}\right)^k\right] = 0,$$ showing that $\lim_{k \rightarrow + \infty} R_k = 0$. Hence, taking the limit $j \rightarrow + \infty$ in $(2.11.4)$, it follows that $$\rho_n(S(t)y) \leq 2 \, C_{m - 1} \, C_{n - (m - 1)} \, n^{m - 1} c_m \, t^{- \frac{n}{m}} \|y\|,$$ if $t \in (0, t_m]$ and $y \in \mathcal{D}$. A quick induction on this formula\footnote{More precisely, an induction proof on $q \geq 1$, on the formula $$\rho_n(S(t)y) \leq \rho_n(S(t)y) \leq D_1 \, [2 \, D_1 \, c_m]^{q - 1} \, n^{(q - 1)(m - 1)} t^{- \frac{n}{m}} \|y\|, \qquad (q - 1)(m - 1) < n \leq q(m - 1)$$ subjected to the restrictions $y \in \mathcal{D}$ and $t \in (0, t_m]$. The base case follows easily from (2.11.1) and from $D_1 \geq 1$. The inductive step follows from a repetition of the argument done so far, in the present item, with $C_l$ substituted by $D_1 \, [2 \, D_1 \, c_m]^{p - 1} \, l^{(p - 1)(m - 1)}$.} gives $$\rho_n(S(t)y) \leq D_1 \, [2 \, D_1 \, c_m]^q \, n^{q(m - 1)} t^{- \frac{n}{m}} \|y\|, \qquad t \in (0, t_m], \, y \in \mathcal{D}$$ where, again, $D_1 := \max \left\{1, \max \left\{C_j: 0 < j \leq m - 1\right\} \right\}$. Consequently, $$\rho_n(S(t)y) \leq D_1 \, [2 e \, D_1 \, c_m]^n n! \, t^{- \frac{n}{m}} \|y\|, \qquad t \in (0, t_m], \, y \in \mathcal{D},$$ where it was used the estimate $n^{q(m - 1)} \leq e^n [q(m - 1)]!$. Therefore, defining $K_0 := D_1$ and $L_0 := 2 e \, D_1 \, c_m$, the constants $C_l$, with $p (m - 1) < l \leq (p + 1)(m - 1)$ and $0 < p \leq q$, may be chosen as $C_l = K_0 \, L_0^l \, l!$.\\

The last task is to extend the inequality just obtained to all $y \in \mathcal{X}$ and to $t \in (t_m, 1]$. Since $S(t)[\mathcal{X}] \subseteq \mathcal{D}$, for all $t > 0$, one has that, for fixed $t_m \geq t > 0$ and $y \in \mathcal{X}$, $$\rho_n(S(t)y) = \inf_{0 < \epsilon < t} \rho_n(S_{t - \epsilon} \, S_\epsilon(y)) \leq K_0 \, L_0^n \, n! \left[\inf_{0 < \epsilon < t} [(t - \epsilon)^{- \frac{n}{m}} \|S_\epsilon(y)\|]\right]$$ $$\leq K_0 \, L_0^n \, n! \left[\inf_{0 < \epsilon < t} [(t - \epsilon)^{- \frac{n}{m}} M e^{w \epsilon}]\right] \|y\| = (M \, K_0) \, L_0^n \, n! \, t^{- \frac{n}{m}} \|y\|.$$ If $t_m \geq 1$, the induction proof is complete. So suppose $0 < t_m < 1$. Then, if $1 \geq t > t_m$ and $y \in \mathcal{X}$, $$\rho_n(S(t)y) = \rho_n(S_{t_m} \, S_{t - t_m}(y)) \leq K_0 \, L_0^n \, n! \, (t_m)^{- \frac{n}{m}} \|S_{t - t_m}(y)\|$$ $$\leq K_0 \, L_0^n \, n! \, (t_m)^{- \frac{n}{m}} M \, e^{w(t - t_m)} \|y\| \leq (M \, e^w \, K_0) \, (t_m^{- \frac{1}{m}} L_0)^n n! \, t^{- \frac{n}{m}} \|y\|,$$ where in the last inequality it was used that $t^{- \frac{n}{m}} \geq 1$ and that $t - t_m < 1$. Hence, with the definitions $K := M \, e^w \, K_0$ and $L := t_m^{- \frac{1}{m}} L_0$ the desired inequality $$\rho_n(S(t)y) \leq K \, L^n n! \, t^{- \frac{n}{m}} \|y\|,$$ for all $y \in \mathcal{X}$ and $t \in (0, 1]$ is finally obtained, concluding the induction step.

It is therefore proved that, with the definitions $K := M \, e^w \, D_1$ and $L := 2 e \, M \, D_1 \, c_m \, t_m^{- \frac{1}{m}}$, it is true that $$\rho_n(S(t)y) \leq K \, L^n n! \, t^{- \frac{n}{m}} \|y\|,$$ for all $n \in \mathbb{N}$, $n \geq 1$, $x \in \mathcal{X}$ and $t \in (0, 1]$.\\

\item Define $\mathcal{S}_0 := \bigcup_{0 < t \leq 1} S(t)[\mathcal{X}]$ and $$C^\omega(\eta) := \left\{x \in \mathcal{X}_\infty: \sum_{n = 0}^{+ \infty} \frac{\rho_n(x)}{n!} s^n < \infty, \text{ for some } s > 0 \right\}.$$ Since $S$ is strongly continuous, $\mathcal{S}_0$ is dense in $\mathcal{X}$. Besides, the estimates from the previous item show that $$\mathcal{S}_0 \subseteq C^\omega(\eta) \cap \mathcal{D} =: \mathcal{D}^\omega,$$ in view of the following argument: fix $0 < t \leq 1$ and $y \in \mathcal{X}$; if one chooses $r > 0$ so that $$r < \frac{1}{L \, t^{- \frac{1}{m}}},$$ then the inequality $$\sum_{n = 0}^{+ \infty} \frac{\rho_n(S(t)y)}{n!} \, r^n \leq K \, \|y\| \sum_{n = 0}^{+ \infty} L^n t^{- \frac{n}{m}} r^n < \infty$$ proves the inclusion. Therefore, $\mathcal{D}^\omega$ is dense in $\mathcal{X}$. Now, because the inclusion $C^\omega(\eta) \subseteq C^\omega(B_k)$ is valid for every $1 \leq k \leq d$ and each $B_k$ is conservative, it follows that every $B_k|_{\mathcal{D}^\omega}$ is an infinitesimal pregenerator of a strongly continuous one parameter group of isometries, by \cite[Theorem 1]{rusinek}. Hence, since $\mathcal{D}^\omega$ is also left invariant by the operators in $\left\{\eta(X): X \in \mathfrak{g}\right\}$,\footnote{To see $C^\omega(\eta)$ is left invariant by the $B_k$'s an analogous argument used to prove that the series obtained by term-by-term differentiation of a complex absolutely convergent power series is still absolutely convergent, and with the same radius of convergence, can be performed.} it follows from \cite[Theorem 3.1]{jorgensengoodman} that the finite-dimensional real Lie algebra $$\left\{\eta(X)|_{\mathcal{D}^\omega}: X \in \mathfrak{g}\right\} \subseteq \text{End}(\mathcal{D}^\omega)$$ is exponentiable. The inclusion $\overline{\eta(X)|_{\mathcal{D}^\omega}} \subset \overline{\eta(X)|_\mathcal{D}} = \eta(X)$ together with dissipativity of $\eta(X)$ yields the equality $\overline{\eta(X)|_{\mathcal{D}^\omega}} = \eta(X) = \overline{\eta(X)|_\mathcal{D}}$, by \cite[Theorem 3.1.15, (3), page 177]{bratteli1}. Therefore, $(\mathcal{D}, \mathfrak{g}, \eta)$ is also exponentiable. But then, again, since each generator $\eta(X)$ is a conservative operator having a dense set of analytic vectors contained in its domain, they actually generate groups of isometries.\hfill $\blacksquare$

\end{itemize}

\textbf{Remark 2.11.1:} Note that the conservativity hypothesis was just invoked in the last item of the proof. In order to obtain a more general exponentiation result (beyond the isometric case) one could, instead of imposing the conservativity hypothesis on the operators $\eta(X)$, $X \in \mathfrak{g}$, replace it by the following (more general) condition: for each $X \in \mathfrak{g}$, there exist $\sigma_x \geq 0$ and $M_x > 0$ such that, for all $n \in \mathbb{N}$, $|\mu| > \sigma_x$ and $y \in \mathcal{D}$, $$\|(\mu I - \eta(X))^n y\| \geq M_x^{-1} \, (|\mu| - \sigma_x)^n \, \|y\|.$$ Then, invoking \cite[Theorem 1]{rusinek}, \cite[Theorem 3.1]{jorgensengoodman} and an adaptation of the proof of \cite[Theorem 3.1.15, (3), page 177]{bratteli1},\footnote{Namely, the following theorem is also true: \textit{Let $S$ be a densely defined closable linear operator on the Banach space $\mathcal{Y}$. If $\overline{S}$ generates a strongly continuous semigroup $t \longmapsto S(t)$ satisfying $$\|S(t)y\| \leq M \exp(wt) \, \|y\|, \qquad y \in \mathcal{Y}, \, t > 0,$$ for constants $M \geq 1$ and $w \geq 0$, and $T$ is an extension of $\overline{S}$ such that $\mu_0 I - T$ is an injective operator, for some $\mu_0 > w$, then $\overline{S} = T$.}} just like it was done in the proof of the last item of Theorem 2.11, one obtains sufficient conditions for the exponentiation of $(\mathcal{D}, \mathfrak{g}, \eta)$ to a strongly continuous representation of a simply connected Lie group (which is not necessarily implemented by isometries).\\

Using Corollary 2.10, it is possible to obtain a ``locally convex space version'' of Theorem 2.11, by making just a few adjustments on its proof. Hence, the following generalization of Theorem 2.11 is true:\\

\textbf{Theorem 2.12 (Exponentiation - Locally Convex Space Version, I):} \textit{Let $\mathcal{X}$ be a complete Hausdorff locally convex space, $\mathfrak{g}$ a real finite-dimensional Lie algebra with an ordered basis $\mathcal{B} := \left(X_k\right)_{1 \leq k \leq d}$, $(\mathcal{D}, \mathfrak{g}, \eta)$ a representation of $\mathfrak{g}$ by closed linear operators on $\mathcal{X}$ and $H_m$ an element of order $m \geq 2$ of $\mathfrak{U}(\eta_\mathcal{D}[\mathfrak{g}])_\mathbb{C}$. Define the subspace of projective analytic vectors for $\eta$ as $$C^\omega_\leftarrow(\eta) := \left\{x \in \mathcal{X}_\infty: \sum_{n = 0}^{+ \infty} \frac{\rho_{p, n}(x)}{n!} s_{x, p}^n < \infty, \text{ for some } s_{x, p} > 0, \text{ for all } p \in \Gamma_m \right\},$$ $\mathcal{D}^\omega := \mathcal{D} \cap C^\omega_\leftarrow(\eta)$ and assume the following hypotheses:}

\begin{enumerate}

\item \textit{each $B_k := \eta(X_k)$ is a $\Gamma_k$-conservative operator, where $\Gamma_k$ is a fundamental system of seminorms for $\mathcal{X}$;}

\item \textit{$-H_m$ is an infinitesimal pregenerator of an equicontinuous semigroup $t \longmapsto S(t)$ satisfying $S(t)[\mathcal{X}] \subseteq \mathcal{D}$, for all $t \in (0, 1]$ and, being a pregenerator of an equicontinuous semigroup, let $\Gamma_m$ be a fundamental system of seminorms for $\mathcal{X}$ with respect to which the operator $-H_m$ has the (KIP), is $\Gamma_m$-dissipative and $S$ is $\Gamma_m$-contractively equicontinuous;}

\item \textit{define $$\Gamma_\mathcal{B} := \bigcup_{1 \leq k \leq d} \Gamma_k,$$ suppose $$\Gamma_\mathcal{B} \subseteq \Gamma_m$$ and that, for each $p \in \Gamma_m$ and $n \in \mathbb{N}$ satisfying $0 < n \leq m - 1$, there exists $C_{p, n} > 0$ for which the estimates $$(2.12.1) \qquad \rho_{p, n}(S(t)y) \leq C_{p, n} \, t^{-\frac{n}{m}} p(y)$$ are verified for all $t \in (0, 1]$ and $y \in \mathcal{X}$, where $B_0 := I$ and $$\rho_{p, n}(x) := \max \left\{p(B_{i_1} \ldots B_{i_k} \ldots B_{i_n}x): 1 \leq k \leq n, 0 \leq i_k \leq d\right\},$$ for all $x \in \mathcal{D}$.}

\end{enumerate}

\textit{Then: (a) The Lie algebra $\left\{\eta(X)|_{\mathcal{D}^\omega}: X \in \mathfrak{g}\right\}$ has a dense subspace of projective analytic vectors and is exponentiable to a strongly continuous locally equicontinuous representation of a Lie group. (b) If, instead, one has the stronger hypotheses that each $\eta(X)$ is $\Gamma_X$-conservative, where $\Gamma_X$ is a fundamental system of seminorms for $\mathcal{X}$, and $$\Gamma := \bigcup_{X \in \mathfrak{g}} \Gamma_X \subseteq \Gamma_m,$$ then $(\mathcal{D}, \mathfrak{g}, \eta)$ exponentiates to a representation by equicontinuous one-parameter groups, and the Lie algebra $\eta_\mathcal{D}[\mathfrak{g}]$ has a dense set of projective analytic vectors.}\\

\textit{(c) If, instead, it happens that all of the operators $\eta(X_k)$, $1 \leq k \leq d$, are $\Gamma_m$-conservative, then $(\mathcal{D}, \mathfrak{g}, \eta)$ exponentiates to a representation by $\Gamma_m$-isometrically equicontinuous one-parameter groups.}\\

\textbf{Proof of Theorem 2.12:} Since $-H_m$ is a pregenerator of an equicontinuous group, the resolvent $\text{R}(1, - \delta \, \overline{H_m})$ is well-defined for all $0 < \delta < 1$ and the identity $$\text{R}(1, - \delta \, \overline{H_m}) = \int_0^{+ \infty} \exp (-t) \, S_{\delta t}(y) \, dt, \qquad y \in \mathcal{X},$$ holds - see \cite[Theorem 1, page 240]{yosida} or \cite[Theorem 3.3, page 172]{babalola} - with $$\int_0^{+ \infty} \exp (-t) \, p(S_{\delta t}(y)) \, dt < \infty,$$ for all $p \in \Gamma_m$. Then, by the estimates of the hypotheses, one finds analogously that $$\int_0^{+ \infty} \exp (-t) \, p(B^\mathsf{u} S_{\delta t}(y)) \, dt$$ $$\leq C_{p, n} \left( \sqrt[m]{(m - n)!} \, \frac{m}{m - n} + I_0 \, \frac{m}{m - n} + 2 \, \sqrt[m]{n!} \right) \delta^{-\frac{n}{m}} \, p(y),$$ for all $y \in \mathcal{X}$. As before, closedness of each operator $B_k$ shows that $(1 + \delta \, \overline{H_m})^{-1}(y) \in \text{Dom }B^\mathsf{u}$ and $$B^\mathsf{u} \, (1 + \delta \, \overline{H_m})^{-1}(y) = \int_0^{+ \infty} \exp (-t) \, B^\mathsf{u} \, S_{\delta t}(y) \, dt.$$ Hence, the estimates obtained show that the inequality $$(2.12.2) \qquad \rho_{p, n}(y) \leq \epsilon^{m - n} p(H_m(y)) + \frac{E_{p, n}}{\epsilon^n}p(y)$$ is valid for every $p \in \Gamma_m$, $0 < \epsilon \leq 1$, $0 < n \leq m - 1$ and $y \in \mathcal{D}$, for some $E_{p, n} > 0$. A very important observation is that (2.12.2) proves, in particular, that the operators $\left(B_k\right)_{1 \leq k \leq d}$ possess the (KIP) with respect to $\Gamma_m$, so the extra hypothesis $\Gamma_\mathcal{B} \subseteq \Gamma_m$ guarantees that these basis elements possess the (KIP) relatively to $\Gamma_\mathcal{B}$, a property which must be verified if one wants to apply Corollary 2.10.\\

The rest of the adaptation of Theorem 2.11 goes as follows:

\begin{itemize}

\item Given $p \in \Gamma_m$ and $q \in \mathbb{N}$, $q \geq 1$, there exists a strictly positive constant $K_{p, q}$ such that $$\rho_{p, n}(S(t)y) \leq K_{p, q} \sup_{0 \leq j \leq q} p(H_m^j(y)),$$ for all $y \in \mathcal{D}$, $(q - 1)(m - 1) < n \leq q \, (m - 1)$ and $t \in (0, 1]$. This implies, in particular, that given $q \in \mathbb{N}$, $q \geq 1$, and $y \in \mathcal{D}$, there exists $K_{y, p, q} > 0$ satisfying $\rho_{p, n}(S(t)y) \leq K_{y, p, q}$, for all $t \in (0, 1]$ and $(q - 1)(m - 1) < n \leq q \, (m - 1)$:\\

The proof is done by induction on $q$, using (2.12.2), in perfect analogy with what is done in Theorem 2.11.\\

\item For each $p \in \Gamma_m$ there exist $K_p, L_p \geq 1$ such that, given $q \in \mathbb{N}$, $q \geq 1$ and $(q - 1)(m - 1) < n \leq q \, (m - 1)$, there exists a strictly positive constant $C_{p, n}$ defined by $C_{p, n} := K_p \, L_p^n \, n!$ satisfying $$\rho_{p, n}(S(t)y) \leq C_{p, n} \, t^{-\frac{n}{m}} p(y), \qquad y \in \mathcal{X}, \, t \in (0, 1].$$ More precisely, there exist $K_p, L_p \geq 1$ such that $$\rho_{p, n}(S(t)y) \leq K_p \, L_p^n \, n! \, t^{- \frac{n}{m}} p(y),$$ for all $n \in \mathbb{N}$, $n \geq 1$, $y \in \mathcal{X}$ and $t \in (0, 1]$:\\

The idea is, again, to proceed by induction on $q$, and the proof goes practically unchanged, evaluating seminorms on the vectors, instead of norms.\\

\item It is now that the context of locally convex spaces really makes a difference, and where Corollary 2.10 will be needed. As before, define $\mathcal{S}_0 := \bigcup_{0 < t \leq 1} S(t)[\mathcal{X}]$. Since $S$ is strongly continuous, $\mathcal{S}_0$ is dense in $\mathcal{X}$. Using the estimates from the previous item shows that $$\mathcal{S}_0 \subseteq C^\omega_\leftarrow(\eta) \cap \mathcal{D} = \mathcal{D}^\omega:$$ fix $p \in \Gamma_m$, $0 < t \leq 1$ and $y \in \mathcal{X}$; since $$\sum_{n = 0}^{+ \infty} \frac{\rho_{p, n}(S(t)y)}{n!} \, r^n \leq K_p \, \|y\| \sum_{n = 0}^{+ \infty} L_p^n t^{- \frac{n}{m}} r^n, \qquad r > 0,$$ if one chooses $r$ so that $$r < \frac{1}{L_p \, t^{- \frac{1}{m}}},$$ then the latter series converges. This proves the inclusion and, therefore, $\mathcal{D}^\omega$ is dense in $\mathcal{X}$. Each $B_k|_{\mathcal{D}^\omega}$ is a $\Gamma_k$-conservative operator having the (KIP) with respect to $\Gamma_k$, so by Lemma 1.6.3 they are all infinitesimal pregenerators of equicontinuous groups. Just as in the Banach space setting the inclusion $C^\omega_\leftarrow(\eta) \subseteq C^\omega_\leftarrow(\eta(X))$ is valid, for every $X \in \mathfrak{g}$. Also, $\mathcal{D}^\omega$ is left invariant by the operators in $\left\{\eta(X): X \in \mathfrak{g}\right\}$. Therefore, since all the elements $\left\{B_k|_{\mathcal{D}^\omega}\right\}_{1 \leq k \leq d}$ possess the (KIP) with respect to $\Gamma_\mathcal{B}$, Corollary 2.10 becomes applicable, so the \textbf{finite-dimensional} real Lie algebra $\left\{\eta(X)|_{\mathcal{D}^\omega}: X \in \mathfrak{g}\right\}$ exponentiates to a strongly continuous locally equicontinuous representation of a Lie group.\\

Now, suppose that every $\eta(X)$, $X \in \mathfrak{g}$, is a $\Gamma_X$-conservative operator and that the inclusion $$\Gamma := \bigcup_{X \in \mathfrak{g}} \Gamma_X \subseteq \Gamma_m$$ holds. Each operator $\overline{\eta(X)|_{\mathcal{D}^\omega}}$ is already known to generate a strongly continuous locally equicontinuous group, but the stronger hypotheses just assumed also imply that $\overline{\eta(X)|_{\mathcal{D}^\omega}}$ is a $\Gamma_X$-conservative operator having the (KIP) with respect to $\Gamma_X$, for all $X \in \mathfrak{g}$. Hence, since all of them have a dense set of projective analytic vectors (by the argument in the preceding paragraph), the stronger conclusion that each $\overline{\eta(X)|_{\mathcal{D}^\omega}}$ is the generator of a $\Gamma_X$-isometrically equicontinuous group follows, by Lemma 1.6.3. The inclusion $\overline{\eta(X)|_{\mathcal{D}^\omega}} \subset \overline{\eta(X)} = \eta(X)$ holds, for every $X \in \mathfrak{g}$, so all that remains to be shown is that the generator of an equicontinuous semigroup has no proper dissipative extensions, an objective which may be accomplished by an adaptation of the argument in \cite[Theorem 3.1.15, (3), page 177]{bratteli1}: let $T_X$ be a proper extension of $\overline{\eta(X)|_{\mathcal{D}^\omega}}$, and let $x \in \text{Dom }T_X \backslash \text{Dom }\overline{\eta(X)|_{\mathcal{D}^\omega}}$. Since $\overline{\eta(X)|_{\mathcal{D}^\omega}}$ is a $\Gamma_X$-dissipative operator that generates an equicontinuous semigroup (actually, it is a $\Gamma_X$-conservative operator that generates an equicontinuous group, but this additional information will not be needed), by \cite[Proposition 3.13]{albanese} and \cite[Theorem 3.14]{albanese} the equality $$\mathcal{X} = \overline{\text{Ran }(I - \eta(X)|_{\mathcal{D}^\omega})} = \text{Ran }(I - \overline{\eta(X)|_{\mathcal{D}^\omega}})$$ must hold. Therefore, there must be $y \in \text{Dom }(I - \overline{\eta(X)|_{\mathcal{D}^\omega}})$ such that $$(I - T_X)x = (I - \overline{\eta(X)|_{\mathcal{D}^\omega}})y = (I - T_X)y.$$ But, then $(I - T_X)(x - y) = 0$, so $\Gamma_X$-dissipativity of $T_X$ together with the fact that $\mathcal{X}$ is Hausdorff implies $x = y \in \text{Dom }\overline{\eta(X)|_{\mathcal{D}^\omega}}$, contradicting the hypothesis on $x$. This proves that $\overline{\eta(X)|_{\mathcal{D}^\omega}}$ does not have any proper $\Gamma_X$-dissipative extensions and that the equality $\overline{\eta(X)|_{\mathcal{D}^\omega}} = \eta(X) = \overline{\eta(X)|_\mathcal{D}}$ holds, for every $X \in \mathfrak{g}$, so $(\mathcal{D}, \mathfrak{g}, \eta)$ is exponentiable to a representation by equicontinuous one-parameter groups (and, for each $X \in \mathfrak{g}$, $\eta(X)$ is the generator of a $\Gamma_X$-isometrically equicontinuous group).\\

Now, suppose, instead of hypothesis 1., that the stronger assumption that all of the operators $\eta(X_k)$, $1 \leq k \leq d$, are $\Gamma_m$-conservative, holds. Using the argumentation of the last paragraph yields $$\overline{B_k|_{\mathcal{D}^\omega}} = B_k, \qquad 1 \leq k \leq d,$$ so that each $B_k$ is the generator of a $\Gamma_m$-isometrically equicontinuous group, and a repetition of the argument in Corollary 2.8 proves that all of the operators $\eta(X)$, $X \in \mathfrak{g}$, are generators of $\Gamma_m$-isometrically equicontinuous one-parameter groups. \hfill $\blacksquare$

\end{itemize}

Before moving on to the next subsection, a small discussion on projective limits of locally convex spaces will be done.\\

Let $\left\{\mathcal{X}_i\right\}_{i \in I}$ be a family of Hausdorff locally convex spaces, where $I$ is a directed set under the partial order $\preceq$ and, for all $i, j \in I$ satisfying $i \preceq j$, suppose that there exists a continuous linear map $\mu_{ij} \colon \mathcal{X}_j \longrightarrow \mathcal{X}_i$ satisfying $\mu_{ij} \circ \mu_{jk} = \mu_{ik}$, whenever $i \preceq j \preceq k$ ($\mu_{ii}$ is, by definition, the identity map, for all $i \in I$). Then, $(\mathcal{X}_i, \mu_{ij}, I)$ is called a \textbf{projective system (or an inverse system)} of Hausdorff locally convex spaces. Consider the vector space $\mathcal{X}$ defined as $$\mathcal{X} := \left\{(x_i)_{i \in I} \in \prod_{i \in I} \mathcal{X}_i: \mu_{ij}(x_j) = x_i, \text{ for all }i, j \in I \text{ with }i \preceq j\right\}$$ and equipped with the relative Tychonoff's product topology or, equivalently, the coarsest topology for which every canonical projection $\pi_j \colon (x_i)_{i \in I} \longmapsto x_j$ is continuous, relativized to $\mathcal{X}$. Then, $\mathcal{X}$ is a Hausdorff locally convex space, and is called the \textbf{projective limit (or inverse limit)} of the family $\left\{\mathcal{X}_i\right\}_{i \in I}$. In this case, the notation $\varprojlim \mathcal{X}_i := \mathcal{X}$ is employed. If $\pi_i[\mathcal{X}]$ is dense in each $\mathcal{X}_i$, then the projective limit is said to be \textbf{reduced}. There is no loss of generality in assuming a projective limit to be reduced: indeed, the projective limit of the family $$\left\{\overline{\pi_i[\mathcal{X}]}\right\}_{i \in I}$$ is equal to the projective limit of $\left\{\mathcal{X}_i\right\}_{i \in I}$, which is $\mathcal{X}$ (see \cite[page 139]{schaefer}).\\

A complete Hausdorff locally convex space $(\mathcal{Y}, \Gamma)$ is isomorphic to the projective limit $\varprojlim \mathcal{Y}_p$, where $\mathcal{Y}_p$ is the Banach completion of $\mathcal{Y} / N_p$, for every $p \in \Gamma$:\footnote{By an isomorphism between two locally convex spaces is meant a continuous bijective linear map with a continuous inverse.} since $\Gamma$ is saturated, it is a directed set. Moreover, if $p, q \in \Gamma$ satisfy $p \preceq q$, define the continuous linear map $\mu_{pq} \colon \mathcal{Y}_q \longrightarrow \mathcal{Y}_p$ as the unique bounded linear extension of the map $[x]_q \longmapsto [x]_p$. Then, the relations $\mu_{pr} = \mu_{pq} \circ \mu_{qr}$ hold, whenever $p, q, r \in \Gamma$ satisfy $p \preceq q \preceq r$, so the projective limit $\varprojlim \mathcal{Y}_p$ is well-defined. The map $$\Phi \colon \mathcal{Y} \longrightarrow \prod_{p \in \Gamma} \mathcal{Y}_p, \qquad \Phi(x) := ([x]_p)_{p \in \Gamma}$$ is linear and injective, because $\mathcal{Y}$ is Hausdorff. Also, by the definition of $\Phi$, it is clear that $\Phi[\mathcal{Y}] \subseteq \varprojlim \mathcal{Y}_p$. To show $\Phi[\mathcal{Y}] = \varprojlim \mathcal{Y}_p$, two auxiliary steps are needed:

\begin{itemize}

\item $\Phi[\mathcal{Y}]$ is a closed subspace of $\varprojlim \mathcal{Y}_p$: consider $x = (x_p)_{p \in \Gamma} \in \overline{\Phi[\mathcal{Y}]} \subseteq \varprojlim \mathcal{Y}_p$ and a net $$\left\{\tilde{x}_\alpha = ([x_\alpha]_p)_{p \in \Gamma}\right\}_{\alpha \in \mathcal{A}}$$ in $\Phi[\mathcal{Y}]$ which converges to $x$. For each fixed $p \in \Gamma$, $\left\{[x_\alpha]_p\right\}_{\alpha \in \mathcal{A}}$ is a Cauchy net in $\mathcal{Y} / N_p$, which is equivalent to the fact that $\left\{x_\alpha\right\}_{\alpha \in \mathcal{A}}$ is a Cauchy net in $\mathcal{Y}$. Since $\mathcal{Y}$ is complete, there exists $y \in \mathcal{Y}$ for which $\lim_\alpha x_\alpha = y$. This implies $\lim_\alpha [x_\alpha]_p = [y]_p$, for every $p \in \Gamma$, so $x = ([y]_p)_{p \in \Gamma} \in \Phi[\mathcal{Y}]$. This ends the proof of this first step.

\item $\Phi[\mathcal{Y}]$ is a dense subspace of $\varprojlim \mathcal{Y}_p$: for every fixed $p \in \Gamma$ the inclusions $$\mathcal{Y} / N_p = \pi_p[\Phi[\mathcal{Y}]] \subseteq \pi_p[\varprojlim \mathcal{Y}_p] \subseteq \mathcal{Y}_p$$ hold, so the fact that $\mathcal{Y} / N_p$ is dense in $\mathcal{Y}_p$ together with the definition of the Tychonoff topology gives the desired result.

\end{itemize}

This establishes the equality $\Phi[\mathcal{Y}] = \varprojlim \mathcal{Y}_p$, and a small adaptation of the argument presented in the first item also shows that $\Phi^{-1}$ is continuous. This ends the desired proof.\\

Now, suppose there are linear operators $T_i \colon \text{Dom }T_i \subseteq \mathcal{X}_i \longrightarrow \mathcal{X}_i$ which are connected by the relations $\mu_{ij} [\text{Dom }T_j] \subseteq \text{Dom }T_i$ and $T_i \circ \mu_{ij} = \mu_{ij} \circ T_j$, whenever $i \preceq j$. Then, the family $\left\{T_i\right\}_{i \in I}$ is said to be a \textbf{projective family of linear operators}. The latter relation ensures that the linear transformation $T$ defined on $\text{Dom }T := \varprojlim \text{Dom }T_i \subseteq \varprojlim \mathcal{X}_i$ by $T(x_i)_{i \in I} := (T_i(x_i))_{i \in I}$ has its range inside $\varprojlim \mathcal{X}_i$, thus defining a linear operator on $\varprojlim \mathcal{X}_i$. This operator is called the \textbf{projective limit of $\bm{\left\{T_i\right\}_{\bm{i \in I}}}$}, as in \cite[page 167]{babalola}.\\

The next objective of this section is to give necessary conditions for the exponentiability problem treated so far.

\subsubsection*{$\bullet$ Strongly Elliptic Operators - Necessary Conditions for Exponentiation}
\addcontentsline{toc}{subsubsection}{Strongly Elliptic Operators - Necessary Conditions for Exponentiation}

\indent

\textbf{Theorem 2.13 (Strongly Elliptic Operators on Locally Convex Spaces):} \textit{Let $(\mathcal{X}, \Gamma)$ be a complete Hausdorff locally convex space and $(\mathcal{X}_\infty, \mathfrak{g}, \eta)$ an \textbf{exponentiable} representation by closed operators, where $\mathcal{B} := \left(X_k\right)_{1 \leq k \leq d}$ is an ordered basis of the finite-dimensional real Lie algebra $\mathfrak{g}$, $\eta(X_k) := B_k$ and $$\mathcal{X}_\infty := \bigcap_{n = 1}^{+ \infty} \bigcap \left\{\text{Dom }[B_{i_1} \ldots B_{i_k} \ldots B_{i_n}]: 1 \leq k \leq n, 1 \leq i_k \leq d \right\}.$$ Let $G$ be the underlying simply connected Lie group such that $V \colon G \longrightarrow \mathcal{L}(\mathcal{X})$ is the strongly continuous locally equicontinuous representation satisfying the property that $dV(X) = \overline{\eta(X)|_{\mathcal{X}_\infty}} = \eta(X)$, for all $X \in \mathfrak{g}$ - in particular, $dV(X_k) = \eta(X_k) = B_k$, for all $1 \leq k \leq d$, and $\mathcal{X}_\infty = C^\infty(V)$ - see Subsection 1.2. Suppose that:}

\begin{enumerate}

\item \textit{the operators $\left(B_k\right)_{1 \leq k \leq d}$ are all $\Gamma$-conservative;}
\item \textit{each $dV(X_k)$, $1 \leq k \leq d$, is the infinitesimal generator of an \textbf{equicontinuous} group $V_k(t) := V(\exp t X_k)$.}

\end{enumerate}

\textit{Define $H_m \in \mathfrak{U}(\partial V[\mathfrak{g}])_\mathbb{C}$ by $$H_m := \sum_{\alpha; |\alpha| \leq m} c_\alpha \, \partial V(X_1)^{\alpha_1} \ldots \partial V(X_k)^{\alpha_k} \ldots \partial V(X_d)^{\alpha_d},$$ so $\text{Dom }H_m = \mathcal{X}_\infty$. Suppose $H_m$ is a \textbf{strongly elliptic operator} on $\mathcal{X}$ and that $-H_m$ is a $\bm{\Gamma}$\textbf{-dissipative} operator. Then, $-H_m$ is an infinitesimal pregenerator of a $\Gamma$-contractively equicontinuous semigroup $S \colon [0, + \infty) \longrightarrow \mathcal{L}(\mathcal{X})$ such that $S(t)[\mathcal{X}] \subseteq \mathcal{X}_\infty$, for all $t > 0$ and, for each $p \in \Gamma$ and $n \in \mathbb{N}$ satisfying $0 < n \leq m - 1$, there exists $C_{p, n} > 0$ validating the estimates $$(2.13.1) \qquad \rho_{p, n}(S(t)y) \leq C_{p, n} \, t^{-\frac{n}{m}} p(y),$$ for every $t \in (0, 1]$ and $y \in \mathcal{X}$. Moreover, if $H_m = -\sum_{k = 1}^d \partial V(X_k)^2$, then the hypothesis of $\Gamma$-dissipativity is superfluous.}\\

\textbf{Proof of Theorem 2.13:} The first task of the proof is to obtain the global (KIP) for $V$ with respect to $\Gamma$ or, in other words, to obtain the inclusions $$V[G][N_p] \subseteq N_p, \qquad p \in \Gamma.$$ To this end, it is sufficient to repeat the argument inside of the proof of Corollary 2.8 to obtain $p(V(g)x) = p(x)$, for all $p \in \Gamma$, $g \in G$ and $x \in \mathcal{X}$, so that the global (KIP) $V[G][N_p] \subseteq N_p$ holds, for all $p \in \Gamma$. This shows that the induced strongly continuous locally equicontinuous group representation $V_p \colon G \longrightarrow \mathcal{X}_p$, defined at each fixed $g \in G$ as the extension by density of the bounded isometric linear operator $V_p(g) \colon [y]_p \longmapsto [V(g)y]_p$, is well-defined, for all $p \in \Gamma$, and that it is a representation by isometries. Because of the arguments just exposed, it follows that if $H_2 := -\sum_{k = 1}^d \partial V(X_k)^2$, then the dissipativity hypothesis on $-H_2 =: \Delta$ is superfluous: fix $p \in \Gamma$. By what was just proved in this paragraph, the operators $\partial V(X)$, $X \in \mathfrak{g}$, have the (KIP) with respect to $\Gamma$. Note first that $$\left\|\left(I - \lambda^2 (\partial V(X_k))_p^2\right)[x]_p\right\|_p = \left\|\left(I - \lambda (\partial V(X_k))_p\right)\left(I + \lambda (\partial V(X_k))_p\right)[x]_p\right\|_p$$ $$\geq \|[x]_p\|_p, \qquad 1 \leq k \leq d, \, \lambda \in \mathbb{R}, \, x \in \mathcal{X}_\infty,$$ by conservativity of the operators $(\partial V(X_k))_p$. This proves $(\partial V(X_k))_p^2$ is a dissipative operator, for all $1 \leq k \leq d$. Fix $[x]_p \in \text{Dom }\Delta_p$. By a characterization of dissipativity on Banach spaces (see pages 174 and 175 of \cite{bratteli1}), if $f \in (\mathcal{X}_p)'$ is a tangent functional at $[x]_p \in \text{Dom }\Delta_p$ (in other words, a bounded linear functional such that $|f([x]_p)| = \|f\| \, \|[x]_p\|_p$, where $\|f\|$ denotes the operator norm of $f$), then $$\text{Re }f\left((\partial V(X_k))_p^2([x]_p)\right) \leq 0, \qquad 1 \leq k \leq d.$$ Therefore, by linearity, $\text{Re }f(\Delta_p([x]_p)) \leq 0$, proving that $\Delta_p$ is a dissipative operator. But $p \in \Gamma$ is arbitrary, so this shows $\Delta$ is $\Gamma$-dissipative.\\

Since all of the operators $\partial V(X_k)$, $1 \leq k \leq d$, possess the (KIP) with respect to $\Gamma$, it follows that $-H_m$ also has this property. Fix $p \in \Gamma$ and $1 \leq k \leq d$. The inclusion $(\partial V(X_k))_p \subset \partial V_p(X_k)$ holds, because $\pi_p[C^\infty(V)] \subseteq C^\infty(V_p)$ and, if $y \in C^\infty(V)$, then $$\left\|\frac{V_p(\exp t X_k)([y]_p) - [y]_p}{t} - (\partial V(X_k))_p([y]_p)\right\|_p = p\left(\frac{V(\exp t X_k)y - y}{t} - \partial V(X_k)(y)\right)$$ goes to 0, as $t \longrightarrow 0$. Hence, defining $$H_p := \sum_{\alpha; |\alpha| \leq m} c_\alpha \, \partial V_p(X_1)^{\alpha_1} \ldots \partial V_p(X_k)^{\alpha_k} \ldots \partial V_p(X_d)^{\alpha_d},$$ the inclusion $$(-H_m)_p = -\sum_{\alpha; |\alpha| \leq m} c_\alpha \, (\partial V(X_1))_p^{\alpha_1} \ldots (\partial V(X_k))_p^{\alpha_k} \ldots (\partial V(X_d))_p^{\alpha_d} \subset -H_p$$ also holds. Now, \cite[Theorem 5.1, page 30]{robinson} says that $-H_p$ is a pregenerator of a strongly continuous semigroup $S_p$ on $\mathcal{X}_p$. Hence, by \cite[Theorem 2.5]{babalola}, the projective limit of the family $\left\{\overline{H_p}\right\}_{p \in \Gamma}$, which will be denoted by $\tilde{H}_m$, is such that $-\tilde{H}_m$ is the generator of a $\Gamma$-semigroup $S \colon [0, +\infty) \longrightarrow \mathcal{L}_\Gamma(\mathcal{X})$ on $\varprojlim \mathcal{X}_p$, which is the projective limit of the semigroups $\left\{S_p\right\}_{p \in \Gamma}$. Since each $-\overline{H_p}$ is dissipative, by hypothesis, it follows by the Feller-Miyadera-Phillips Theorem together with \cite[Proposition 3.11 (iii)]{albanese} (applied to Banach spaces) that $-\overline{H_p}$ is actually the generator of a \textbf{contraction semigroup}. Therefore, $-\tilde{H}_m$ must be the generator of a $\Gamma$-contractively equicontinuous semigroup on $\varprojlim \mathcal{X}_p$. It should be noted, at this point, that all of the results proved on $\varprojlim \mathcal{X}_p$ carry over to $\mathcal{X}$ (and vice-versa), since $\mathcal{X}$ and $\varprojlim \mathcal{X}_p$ are isomorphic locally convex spaces. Therefore, in what follows, $\mathcal{X}$ and $\varprojlim \mathcal{X}_p$ will be treated as if they were the \textbf{same} space. Hence, the relation $$-H_m = -\varprojlim (H_m)_p \subset -\varprojlim \overline{H_p} = -\tilde{H}_m$$ holds.\\

To conclude that $-\overline{H_m}$ is the generator of a $\Gamma$-contractively equicontinuous semigroup, it suffices to show that $-\overline{H_m} = -\tilde{H}_m$, an equality which will be obtained after the next paragraph.\\

To prove the inclusion $S(t)[\mathcal{X}] \subseteq \mathcal{X}_\infty$, for all $t > 0$, note first that, if $1 \leq k \leq d$ and $y \in \mathcal{X}$, then the function $t \longmapsto V(\exp t X_k)y$ is infinitely differentiable on $[0, +\infty)$ if, and only if, $t \longmapsto V_p(\exp t X_k)([y]_p)$ is infinitely differentiable on $[0, +\infty)$, for all $p \in \Gamma$. Hence, using the alternative description of smooth vectors for a strongly continuous locally equicontinuous representation proved in Subsection 1.2, one concludes that $y$ belongs to $C^\infty(V)$ if, and only if, $[y]_p$ belongs to $C^\infty(V_p)$, for every $p \in \Gamma$. But \cite[Theorem 5.1]{robinson} says, in particular, that $S_p(t)[\mathcal{X}_p] \subseteq C^\infty(V_p)$, for all $t > 0$. Therefore, $[S(t)y]_p = S_p([y]_p)$ belongs to $C^\infty(V_p)$, for all $p \in \Gamma$ and $t > 0$, so $S(t)y$ must belong to $\mathcal{X}_\infty$, for each fixed $t > 0$.\\

The argument in the previous paragraph shows, in particular, that $S(t)[\mathcal{X}_\infty] \subseteq \mathcal{X}_\infty$, for all $t > 0$, so Lemma 1.4.2 says that $\mathcal{X}_\infty$ is a core for $-\tilde{H}_m$. Hence, $$-\tilde{H}_m = -\overline{\tilde{H}_m|_{\mathcal{X}_\infty}} = -\overline{H_m},$$ as claimed.\\

Finally, validity of 2.13.1 is a consequence of \cite[Corollary 5.6, page 44]{robinson}, for if $p \in \Gamma$, $B^\mathsf{u}$ is a monomial of size $0 < |\mathsf{u}| := n \leq m - 1$ in the generators $\left\{dV(X_k)\right\}_{1 \leq k \leq d}$ and $(B^\mathsf{u})_p$ is the induced monomial in the quotient $\mathcal{X} / N_p$, then $$p(B^\mathsf{u} S(t)y) = \|(B^\mathsf{u})_p S_p(t)[y]_p\|_p \leq C_{p, n} \, t^{-\frac{n}{m}} \|[y]_p\|_p$$ $$= C_{p, n} \, t^{-\frac{n}{m}} p(y), \qquad t \in (0, 1], \, y \in \mathcal{X},$$ for some constants $C_{p, n} > 0$, since $(dV(X))_p \subset dV_p(X)$, for every $p \in \Gamma$ and $X \in \mathfrak{g}$. \hfill $\blacksquare$\\

\textbf{Observation 2.13.1:} In the same context of Theorem 2.13, it is possible to show that if the strongly continuous locally equicontinuous representation $V \colon G \longrightarrow \mathcal{L}(\mathcal{X})$ satisfies the condition that, for each $p \in \Gamma$, there exists $M_p > 0$ and another $q \in \Gamma$ satisfying $$\sup_{g \in G} p(V(g)x) \leq M_p \, q(x), \qquad x \in \mathcal{X}$$ - this is a condition which is automatically satisfied if $G$ is \textbf{compact} - then hypothesis 1 of the theorem which asks for conservativity of the operators $\left(B_k\right)_{1 \leq k \leq d}$ with respect to the \textbf{same} $\Gamma$ is automatically fulfilled. To see this, fix $p \in \Gamma$. In a similar fashion of that done in Subsection 1.5, one may define the seminorm $$x \longmapsto \tilde{p}(x) := \sup_{g \in G} p(V(g)x) < \infty$$ and see that $$p(x) \leq \tilde{p}(x) \leq M_p \, q(x) \leq M_p \, \tilde{q}(x), \qquad x \in \mathcal{X},$$ proving that the families $\Gamma$ and $\tilde{\Gamma} := \left\{\tilde{p}: p \in \Gamma\right\}$ are equivalent. But even more, $\tilde{\Gamma}$ has the very useful property that $$\tilde{p}(V(h)x) = \sup_{g \in G} p(V(gh)x) = \sup_{g \in G} p(V(g)x) = \tilde{p}(x),$$ for all $h \in G$ and $x \in \mathcal{X}$, so $$\tilde{p}(V(h)x) = \tilde{p}(x), \qquad h \in G, \, x \in \mathcal{X}.$$ Repeating the argument in Subsection 1.5 yields conservativity of the operators $\left\{dV(X_k)\right\}_{1 \leq k \leq d}$ with respect to $\tilde{\Gamma}$. Hence, the assertion follows.

\subsubsection*{$\bullet$ Exponentiation in Locally Convex Spaces - Characterization}
\addcontentsline{toc}{subsubsection}{Exponentiation in Locally Convex Spaces - Characterization}

\indent

Finally, it is possible to combine Theorems 2.12 and 2.13 (see also Subsection 1.5) to obtain a characterization of Lie algebras of linear operators which exponentiate to representations by $\Gamma$-isometrically equicontinuous one-parameter groups, generalizing \cite[Theorem 3.9]{bratteliheat}:\footnote{Note that the statement of Theorem 2.14, below, is written in a slightly different way than that of \cite[Theorem 2.14]{rodrigotese}. The statement given here is more elegant.}\\

\textbf{Theorem 2.14 (Exponentiation - Locally Convex Space Version, II):} \textit{Let $(\mathcal{X}, \Gamma)$ be a complete Hausdorff locally convex space and $(\mathcal{X}_\infty, \mathfrak{g}, \eta)$ a representation by closed operators, where $\mathcal{B} := \left(X_k\right)_{1 \leq k \leq d}$ is an ordered basis of the finite-dimensional real Lie algebra $\mathfrak{g}$, $\eta(X_k) := B_k$ and $$\mathcal{X}_\infty := \bigcap_{n = 1}^{+ \infty} \bigcap \left\{\text{Dom }[B_{i_1} \ldots B_{i_k} \ldots B_{i_n}]: 1 \leq k \leq n, 1 \leq i_k \leq d \right\}.$$ Suppose that $$H_m := \sum_{\alpha; |\alpha| \leq m} c_\alpha \, \eta_{\mathcal{X}_\infty}(X_1)^{\alpha_1} \ldots \eta_{\mathcal{X}_\infty}(X_k)^{\alpha_k} \ldots \eta_{\mathcal{X}_\infty}(X_d)^{\alpha_d}$$ is an element of $\mathfrak{U}(\eta_{\mathcal{X}_\infty}[\mathfrak{g}])_\mathbb{C}$ of order $m \geq 2$ which is strongly elliptic and that $-H_m$ is $\Gamma$-dissipative.\footnote{If $H_m = -\sum_{k = 1}^d \eta_{\mathcal{X}_\infty}(X_k)^2$, then the dissipativity hypothesis on $-H_m$ will be superfluous, by the usual arguments.} Then, the following are equivalent:}

\begin{enumerate}

\item \textit{$(\mathcal{X}_\infty, \mathfrak{g}, \eta)$ exponentiates to a strongly continuous locally equicontinuous representation \linebreak $V \colon G \longrightarrow \mathcal{L}(\mathcal{X})$ of a simply connected Lie group $G$, having $\mathfrak{g}$ as its Lie algebra, such that $dV(X) = \overline{\eta(X)|_{\mathcal{X}_\infty}}$ is the generator of a $\Gamma$-isometrically equicontinuous one-parameter group, for all $X \in \mathfrak{g}$.}
\item \textit{The operators $\left(B_k\right)_{1 \leq k \leq d}$ are $\Gamma$-conservative, $-H_m$ is an infinitesimal pregenerator of a $\Gamma$-contractively equicontinuous semigroup $S \colon [0, + \infty) \longrightarrow \mathcal{L}(\mathcal{X})$ such that $S(t)[\mathcal{X}] \subseteq \mathcal{X}_\infty$, for all $t > 0$ and, for each $p \in \Gamma$ and $n \in \mathbb{N}$ satisfying $0 < n \leq m - 1$, there exists $C_{p, n} > 0$ for which $$(2.14.1) \qquad \rho_{p, n}(S(t)y) \leq C_{p, n} \, t^{-\frac{n}{m}} p(y),$$ where $t \in (0, 1]$ and $y \in \mathcal{X}$.}

\end{enumerate}

\textit{Also, if one of the above conditions is fulfilled, then $\eta$ has a dense set of projective analytic vectors.}

\section{Some Applications to Locally Convex Algebras}

\indent

\textbf{Definition 3.1:} \textit{Let $\mathcal{A}$ be a Hausdorff locally convex space (over $\mathbb{C}$) equipped with an associative product which is compatible with the sum via the distributive rules $(a + b)c = ac + bc$, $a(b + c) = ab + ac$, for all $a, b, c \in \mathcal{A}$. Also, suppose that the multiplication is \textbf{separately continuous}: for each fixed $a \in \mathcal{A}$, the maps $b \longmapsto ab$ and $b \longmapsto ba$ are continuous. Then, $\mathcal{A}$ is called a \textbf{locally convex algebra}. If the underlying vector space of $\mathcal{A}$ is a complete (respectively, sequentially complete) Hausdorff locally convex space, then it said to be a \textbf{complete} (respectively, \textbf{sequentially complete}) \textbf{locally convex algebra}. It may also happen that the product satisfies a stronger continuity assumption: if the map $(a, b) \longmapsto ab$ from $\mathcal{A} \times \mathcal{A}$ to $\mathcal{A}$ is continuous, then the product is said to be \textbf{jointly continuous}. If $p$ is a seminorm on $\mathcal{A}$ satisfying $p(ab) \leq p(a) p(b)$, for all $a, b \in \mathcal{A}$, then it is said to be \textbf{submultiplicative}, or an \textbf{m-seminorm}. Therefore, if $\mathcal{A}$ possesses a fundamental system of submultiplicative seminorms, then its product is jointly continuous.}\\

\textit{Suppose $\mathcal{A}$ is a locally convex algebra equipped with an involution $*$, which is, by definition, an antilinear map on $\mathcal{A}$ satisfying $(ab)^* = b^*a^*$ and $(a^*)^* = a$, for all $a, b \in \mathcal{A}$. Then, $\mathcal{A}$ is called a \textbf{locally convex algebra with involution}. If $* \colon a \longmapsto a^*$ is a continuous map, $\mathcal{A}$ is called a \textbf{locally convex $\bm{*}$-algebra}.}\\

Interesting examples of locally convex ($*$-)algebras are the Arens-Michael ($*$-)algebras (which are defined as complete m($^*$)-convex algebras - see \cite[Chapter I]{fragoulopoulou}), the von Neumann algebras, the locally C$^*$-algebras (which are going to be properly defined, soon), the C$^*$-like locally convex $*$-algebras (\cite{inouekursten}, \cite{arens}) and the GB$^*$-algebras (\cite{allan2}, \cite{allan}, \cite{bhattgb}, \cite{bhattgb2}, \cite{dixon1}, \cite{dixon2}, \cite{fraginouekursten}). In particular, all locally C$^*$-algebras are C$^*$-like locally convex $*$-algebras, and the latter are special types of GB$^*$-algebras. If the topology of an Arens-Michael ($*$-)algebra can be generated by a countable family of seminorms, then it is called a \textbf{Fr\'echet ($\bm{*}$-)algebra}.

See also \cite{fragoulopoulou2}, \cite{weigtzarakas}, \cite{blackadarcuntz} and \cite{bhatt2}.

\subsection*{$\bullet$ Exponentiation in Complete Locally Convex Algebras}
\addcontentsline{toc}{subsection}{Exponentiation in Complete Locally Convex Algebras}

\indent

Now, it will be shown that, once a Lie algebra of operators on a locally convex algebra $\mathcal{A}$ (or on a locally convex $*$-algebra $\mathcal{A}$) is exponentiated, there is an automatic compatibility of the group representation $V \colon G \longrightarrow \mathcal{L}(\mathcal{A})$ obtained with the additional algebraic structure, and not only with the vector space structure. To this end, an idea which appears in the proof of \cite[Theorem 1]{brattelirob} was borrowed and adapted to the present context.\\

\textbf{Lemma 3.2 (Automatic Compatibility with the Additional Algebraic Operations):} \textit{Let $(\mathcal{A}, \Gamma)$ be a locally convex algebra (respectively, locally convex $*$-algebra), $\mathcal{D} \subseteq \mathcal{A}$ a dense subalgebra (respectively, $*$-subalgebra) and $\mathcal{L} \subseteq \text{End}(\mathcal{D})$ be an exponentiable Lie algebra of derivations (respectively, $*$-derivations) to a strongly continuous locally equicontinuous representation $g \longmapsto V(g)$ of a Lie group. Then, each one-parameter group $V(t, \overline{\delta}) \colon t \longmapsto V(\exp t \delta)$, $\delta \in \mathcal{L}$, of continuous linear operators is actually a representation by automorphisms (respectively, $*$-automorphisms).}\\

\textbf{Proof of Lemma 3.2:} First, assume that $\mathcal{A}$ is only a locally convex algebra. Fix $\delta \in \mathcal{L}$ and $t \in \mathbb{R}$. The operator $V(t, \overline{\delta})$ will be abbreviated to $V_t$, to facilitate the reading process. To prove $V_t$ is an algebra automorphism, the only thing that remains to be shown is that it preserves the multiplication of $\mathcal{A}$. Define $$F \colon s \longmapsto V_{-s}(V_s(c) \, V_s(d)),$$ where $c, d$ are fixed elements of $\text{Dom }\overline{\delta}$. Since the multiplication operation on $\mathcal{A}$ is separately continuous, the linear map from $\text{Dom }\overline{\delta}$ to $\mathcal{A}$ defined by $$L_a \colon b \longmapsto ab, \qquad a \in \mathcal{A},$$ is linear and continuous, when $\text{Dom }\overline{\delta}$ is equipped with the topology defined by the family of seminorms $\left\{a \longmapsto p(a) + p(\overline{\delta}(a)): p \in \Gamma\right\}$. Moreover, since $$\frac{V_{s + h}(a) \, b - V_s(a) \, b}{h} = V_s \left(\frac{V_h(a) - (a)}{h}\right) b, \qquad s \in \mathbb{R}, \, h \in \mathbb{R} \backslash \left\{0\right\}, \, a, b \in \text{Dom }\overline{\delta},$$ the strong derivative of $s \longmapsto L_{V(s)(a)}$, for each fixed $a \in \text{Dom }\overline{\delta}$, is $$s \longmapsto L_{V(s)(\overline{\delta}(a))},$$ which maps $\mathbb{R}$ to the space of continuous linear maps from $\text{Dom }\overline{\delta}$ to $\mathcal{A}$.

The function $s \longmapsto V_s$ also maps $\mathbb{R}$ to the space of continuous linear maps from $\text{Dom }\overline{\delta}$ to $\mathcal{A}$ in a strongly differentiable way, with $V_t \circ \overline{\delta}$ being its strong derivative at $s = t$. Therefore, applying \cite[Theorem A.1, page 440]{jorgensenmoore} to $F$, twice, gives $$F'(t) = - V_{-t} \, \overline{\delta}(V_t(c) \, V_t(d)) + V_{-t}(V_t(\overline{\delta}(c)) \, V_t(d) + V_t(c) \, V_t(\overline{\delta}(d)))$$ $$= - V_{-t} \left(\overline{\delta}(V_t(c)) \, V_t(d) + V_t(c) \, \overline{\delta}(V_t(d))\right) + \, V_{-t} \left(\overline{\delta}(V_t(c)) \, V_t(d) + V_t(c) \, \overline{\delta}(V_t(d))\right) = 0.$$ Since $F(0) = cd$, a quick application of the Hahn-Banach Theorem gives $$V_{-t}(V_t(c) \, V_t(d)) = F(t) = F(0) = cd,$$ so the arbitrariness of $t$, $c$ and $d$, together with denseness of $\text{Dom }\overline{\delta}$ in $\mathcal{A}$, prove that $V_t$ is an algebra homomorphism of $\mathcal{A}$.\\

Now, suppose $\mathcal{A}$ is a locally convex $*$-algebra. Defining the function $$G \colon s \longmapsto V_{-s}(V_s(c)^*),$$ where $c$ is a fixed element of $\text{Dom }\overline{\delta}$, and applying a similar reasoning, one obtains that $V_t$ preserves the (continuous) involution operation, for each $t \in \mathbb{R}$. \hfill $\blacksquare$\\

\textbf{Remark 3.2.1:} It should be mentioned that the proof above, although simple, has some subtleties. The precise specification of the locally convex spaces involved in order to apply \cite[Theorem A.1, page 440]{jorgensenmoore} is one of them. The other one is that the usual formula for the derivative of the product of two functions defined on $\mathbb{R}$ had to be derived via an alternative approach: since the multiplication of $\mathcal{A}$ is \textbf{not} necessarily jointly continuous, the usual proof of this formula (as it is done in many Elementary Calculus books - probably in most of them) may not be adapted to this context.\\

Finally, some applications of the exponentiation theorems of Section 2 may be given.\\

\textbf{Definition 3.3:} \textit{Let $\mathcal{A}$ be a locally convex $*$-algebra and $\mathcal{D} \subseteq \mathcal{A}$ a $*$-subalgebra. Suppose $(\mathcal{D}, \mathfrak{g}, \eta)$ is a representation of $\mathfrak{g}$ by closed linear operators on $\mathcal{A}$, in the sense of the definition made in Section 2. Suppose, also, that every element of $\text{Ran }\eta$ is a $*$-derivation - a $*$-derivation $\delta$ is a linear, $*$-preserving map defined on a $*$-subalgebra of $\mathcal{A}$ which satisfies the Leibniz rule: $\delta(ab) = \delta(a)b + a \delta(b)$, $a, b \in \text{Dom }\delta$. Then, $(\mathcal{D}, \mathfrak{g}, \eta)$ is called a \textbf{representation of $\bm{\mathfrak{g}}$ by closed $\bm{*}$-derivations on $\bm{\mathcal{A}}$}. Now, whenever $(\mathcal{D}, \mathfrak{g}, \eta)$ exponentiates to a Lie group strongly continuous locally equicontinuous representation $\alpha \colon G \longrightarrow \mathcal{L}(\mathcal{A})$ (in the context of algebras the letter ``$\alpha$'', instead of ``$V$'', will be used), it is a consequence of Lemma 3.2 that each one-parameter group $t \longmapsto \alpha(\exp tX)$, $X \in \mathfrak{g}$, is actually implemented by $*$-automorphisms. In other words, every operator $\alpha(\exp tX)$ is a $*$-automorphism on $\mathcal{A}$.\footnote{A $*$-automorphism is a continuous $*$-preserving algebraic isomorphism of $\mathcal{A}$ onto $\mathcal{A}$ which possesses a continuous inverse.} An analogous definition can be done if $\mathcal{A}$ is just a locally convex algebra, substituting the words ``$*$-subalgebra'' by ``algebra'', ``$*$-derivation'' by ``derivation'' and ``$*$-automorphism'' by ``automorphism''.}\\

Combining Theorem 2.14 with Lemma 3.2 one can obtain the following exponentiation theorem for complete locally convex ($*$-)algebras:\footnote{Just as in Theorem 2.14, the statement of Theorem 3.4, below, is written in a slightly more elegant way than \cite[Theorem 3.7]{rodrigotese}.}\\

\textbf{Theorem 3.4 (Exponentiation - Complete Locally Convex Algebras):} \textit{Let $\mathcal{A}$ be a \textbf{complete} locally convex algebra, $\mathfrak{g}$ a real finite-dimensional Lie algebra with an ordered basis $\mathcal{B} := \left(X_k\right)_{1 \leq k \leq d}$ and $(\mathcal{A}_\infty, \mathfrak{g}, \eta)$ a representation of $\mathfrak{g}$ by closed derivations, with $$\mathcal{\mathcal{A}}_\infty := \bigcap_{n = 1}^{+ \infty} \bigcap \left\{\text{Dom }[\delta_{i_1} \ldots \delta_{i_k} \ldots \delta_{i_n}]: 1 \leq k \leq n, 1 \leq i_k \leq d \right\},$$ where $\delta_X := \eta(X)$, $\delta_k := \eta(X_k)$, $1 \leq k \leq d$. Suppose that $$H_m := \sum_{\alpha; |\alpha| \leq m} c_\alpha \, \eta_{\mathcal{A}_\infty}(X_1)^{\alpha_1} \ldots \eta_{\mathcal{A}_\infty}(X_k)^{\alpha_k} \ldots \eta_{\mathcal{A}_\infty}(X_d)^{\alpha_d}$$ is an element of $\mathfrak{U}(\eta_{\mathcal{A}_\infty}[\mathfrak{g}])_\mathbb{C}$ of order $m \geq 2$ which is strongly elliptic and that $-H_m$ is $\Gamma$-dissipative. Then, the following are equivalent:}

\begin{enumerate}

\item \textit{$(\mathcal{A}_\infty, \mathfrak{g}, \eta)$ exponentiates to a strongly continuous locally equicontinuous representation \linebreak $\alpha \colon G \longrightarrow \mathcal{L}(\mathcal{A})$ \textbf{by automorphisms} of a simply connected Lie group $G$, having $\mathfrak{g}$ as its Lie algebra, such that $d\alpha(X) = \overline{\eta(X)|_{\mathcal{A}_\infty}}$ is the generator of a $\Gamma$-isometrically equicontinuous one-parameter group, for all $X \in \mathfrak{g}$.}
\item \textit{The operators $\left\{\delta_k\right\}_{1 \leq k \leq d}$ are $\Gamma$-conservative, $-H_m$ is an infinitesimal pregenerator of a $\Gamma$-contractively equicontinuous semigroup $t \longmapsto S(t)$ satisfying $S(t)[\mathcal{A}] \subseteq \mathcal{A}_\infty$, for all $t \in (0, 1]$; moreover, for each $p \in \Gamma$ and $n \in \mathbb{N}$ satisfying $0 < n \leq m - 1$, there exists $C_{p, n} > 0$ for which the estimates $$\rho_{p, n}(S(t)a) \leq C_{p, n} \, t^{-\frac{n}{m}} p(a)$$ are verified for all $t \in (0, 1]$ and $a \in \mathcal{A}$.}

\end{enumerate}

\textit{If $H_2 = \Delta := \sum_{k = 1}^d \delta_k^2$, then the $\Gamma$-dissipativity hypothesis on $H_2$ is superfluous, by the usual arguments. Also, if one of the above conditions is fulfilled, then $\eta$ has a dense set of projective analytic vectors.}\\

\textit{Theorem 3.4 remains valid if $\bm{\mathcal{A}}$ \textbf{is a complete locally convex $\bm{*}$-algebra} and $(\mathcal{A}_\infty, \mathfrak{g}, \eta)$ is a representation of $\mathfrak{g}$ by closed $*$-derivations, substituting the word ``automorphisms'' by ``$*$-automorphisms''.}

\subsection*{$\bullet$ Exponentiation in Locally C$\bm{^*}$-Algebras}
\addcontentsline{toc}{subsection}{Exponentiation in Locally C$^*$-Algebras}

\indent

In the special case where $\mathcal{A}$ is a locally C$^*$-algebra, an improvement of Theorem 3.4 may be obtained, at least for the case in which the operator dominating the representation is the negative of the Laplacian. Theorem 3.6, below, is in perfect analogy with the C$^*$ case (see \cite[Corollary 2]{brattelidissipative}), the only difference being that the result, below, is proved for arbitrary dense core domains.\\

\textbf{Definition 3.5 (Locally C$\bm{^*}$-Algebras):} \textit{A complete locally convex $*$-algebra $\mathcal{A}$ whose topology $\tau$ is generated by a directed family $\Gamma$ of seminorms satisfying the properties: $$(i) \, \, p(ab) \leq p(a) p(b), \quad (ii) \, \, p(a^*) = p(a), \quad (iii) \, \, p(a^*a) = p(a)^2, \qquad p \in \Gamma, a, b \in \mathcal{A},$$ is called a \textbf{locally C$\bm{^*}$-algebra} - a seminorm $p$ satisfying (iii) is called a C$^*$-seminorm.\footnote{As a consequence of \cite[Theorem 7.2, page 101]{fragoulopoulou}, a C$^*$-seminorm automatically satisfies properties (i) and (ii), above.} It is said to be unital if it possesses an identity: in other words, if it has an element 1 satisfying the property $1 \cdot a = a \cdot 1 = a$, for all $a \in \mathcal{A}$. Hence, the concept of a locally C$^*$-algebra generalizes that of a C$^*$-algebra.}\\

The terminology used to refer to locally C$^*$-algebras varies in the literature, and they are also called \textbf{pro-C$\bm{^*}$-algebras} or \textbf{LMC$\bm{^*}$-algebras} (see \cite[page 99]{fragoulopoulou} and \cite{phillips}, which calls them pro-C$^*$-algebras; the terminology ``locally C$^*$-algebra'' was apparently introduced by \cite{inoue}; \cite{becker} calls them LMC$^*$-algebras). If $\mathcal{A}$ is a locally C$^*$-algebra, the algebra of the so-called ``bounded elements'', defined as $$b(\mathcal{A}) := \left\{a \in \mathcal{A}: \sup_{p \in \Gamma} p(a) < \infty\right\},$$ is a C$^*$-algebra under the norm $$\|a\|_\infty := \sup_{p \in \Gamma} p(a),$$ and it is $\tau$-dense in $\mathcal{A}$ \cite[Proposition 1.11 (4)]{phillips}. For example, for the locally C$^*$-algebra $C(\mathbb{R}^d)$ of all continuous functions on $\mathbb{R}^d$, topologized by the family of seminorms defined by sup-norms over the compact subsets of $\mathbb{R}^d$, the corresponding C$^*$-algebra of bounded elements is the algebra $C_b(\mathbb{R}^d)$ of bounded continuous functions on $\mathbb{R}^d$.\\

Let $(\Lambda, \preceq)$ be a directed set and $\left\{\mathcal{H}_\lambda\right\}_{\lambda \in \Lambda}$ a family of Hilbert spaces such that $$\mathcal{H}_\lambda \subseteq \mathcal{X}_\nu, \quad \langle \cdot, \cdot \rangle_\lambda = \langle \cdot, \cdot \rangle_\nu|_{\mathcal{H}_\lambda},$$ whenever $\lambda, \nu \in \Lambda$ satisfy $\lambda \preceq \nu$. If $$i_{\lambda \nu} \colon \mathcal{H}_\lambda \longrightarrow \mathcal{H}_\nu, \quad \lambda \preceq \nu,$$ denotes the inclusion of $\mathcal{H}_\lambda$ into $\mathcal{H}_\nu$, then $(\mathcal{H}_\lambda, i_{\lambda \nu})$ is an inductive system of Hilbert spaces. The inductive limit vector space $$\mathcal{H} := \varinjlim \mathcal{H}_\lambda = \bigcup_{\lambda \in \Lambda} \mathcal{H}_\lambda,$$ when equipped with the inductive limit topology (in other words, the finest topology making the canonical inclusions $i_\lambda \colon \mathcal{H}_\lambda \longrightarrow \mathcal{H}$, $\lambda \in \Lambda$, continuous), is called a locally Hilbert space. Given a family of linear maps $\left\{T_\lambda\right\}_{\lambda \in \Lambda}$, where $T_\lambda \in \mathcal{L}(\mathcal{H}_\lambda)$, for every $\lambda \in \Lambda$ (a family having this property is called an inductive family), and such that $T_\nu|_{\mathcal{H}_\lambda} = T_\lambda$, whenever $\lambda \preceq \nu$, one can associate a unique continuous linear map $$T := \varinjlim T_\lambda \colon \mathcal{H} \longrightarrow \mathcal{H}, \quad T|_{\mathcal{H}_\lambda} := T_\lambda, \quad \lambda \in \Lambda.$$ Define $\mathcal{L}(\mathcal{H})$ as the set of all continuous linear maps of $\mathcal{H}$ into $\mathcal{H}$ such that $T = \varinjlim T_\lambda$, for some family of operators $\left\{T_\lambda\right\}_{\lambda \in \Lambda}$ satisfying $T_\nu|_{\mathcal{H}_\lambda} = T_\lambda$, whenever $\lambda \preceq \nu$, where $T_\lambda \in \mathcal{L}(\mathcal{H}_\lambda)$. Then, $\mathcal{L}(\mathcal{H})$ becomes an algebra with the natural operations induced by the operator algebras $\mathcal{L}(\mathcal{H}_\lambda)$, $\lambda \in \Lambda$. The involution $*$ on an operator $T$ in $\mathcal{L}(\mathcal{H})$ is defined as the inductive limit operator $T^* := \varinjlim T_\lambda^*$, and the topology on $\mathcal{L}(\mathcal{H})$ is generated by the family $\left\{\|T_\lambda\|_\lambda\right\}_{\lambda \in \Lambda}$ of C$^*$-seminorms, where $\|\, \cdot \,\|_\lambda$ denotes the usual operator norm on $\mathcal{L}(\mathcal{H}_\lambda)$. When equipped with this topology and these algebraic operations, $\mathcal{L}(\mathcal{H})$ becomes a locally C$^*$-algebra \cite[Example 7.6 (5), page 106]{fragoulopoulou}. By \cite[Theorem 8.5, page 114]{fragoulopoulou}, every locally C$^*$-algebra $\mathcal{A}$ is $*$-isomorphic to a closed $*$-subalgebra of $\mathcal{L}(\mathcal{H})$, where $\mathcal{H}$ is a locally Hilbert space. A commutative Gelfand-Naimark-type theorem is also available for locally C$^*$-algebras, just as in the C$^*$ normed setting (see \cite[Theorem 9.3, page 117]{fragoulopoulou} and \cite[Remarks 9.6, page 119]{fragoulopoulou}; see, also, \cite[Section 10 - Functional Calculus, page 119]{fragoulopoulou}).\\

Also, by \cite[Theorem 10.24, page 133]{fragoulopoulou}, every locally C$^*$-algebra $(\mathcal{A}, \Gamma)$ has the property that $\overline{\mathcal{A}/N_p} = \mathcal{A}/N_p$, for every $p \in \Gamma$. In other words, their quotient algebras are already complete, when seen as normed spaces.\\

\textbf{Theorem 3.6 (Exponentiation - Locally C$\bm{^*}$-Algebras, I):} \textit{Let $\mathcal{A}$ be a locally C$^*$-algebra, $\mathcal{D} \subseteq \mathcal{A}$ a dense $*$-subalgebra, $\mathfrak{g}$ a real finite-dimensional Lie algebra with an ordered basis $\mathcal{B} := \left(X_k\right)_{1 \leq k \leq d}$, $(\mathcal{D}, \mathfrak{g}, \eta)$ a representation of $\mathfrak{g}$ by closed $*$-derivations with $\delta_X := \eta(X)$, $\delta_k := \eta(X_k)$, $1 \leq k \leq d$, and $\Delta := \sum_{k = 1}^d \delta_k^2$. Assume that the following hypotheses are valid:}\\

\textit{$\Delta$ is an infinitesimal pregenerator of an equicontinuous semigroup $t \longmapsto S(t)$ satisfying $S(t)[\mathcal{A}] \subseteq \mathcal{D}$, for all $t \in (0, 1]$ and, being a pregenerator of an equicontinuous semigroup, let $\Gamma_2$ be a fundamental system of seminorms for $\mathcal{A}$ with respect to which the operator $\Delta$ has the (KIP), is $\Gamma_2$-dissipative and $S$ is $\Gamma_2$-contractively equicontinuous; also, suppose that, for every $p \in \Gamma_2$, there exists $C_p > 0$ for which the estimates $$(3.6.1) \qquad \rho_{p, 1}(S(t)a) \leq C_p \, t^{-\frac{1}{2}} p(a)$$ are verified, for all $t \in (0, 1]$ and $a \in \mathcal{A}$.}\\

\textit{Then, $(\mathcal{D}, \mathfrak{g}, \eta)$ exponentiates to a strongly continuous locally equicontinuous representation $\alpha \colon G \longrightarrow \mathcal{L}(\mathcal{A})$ \textbf{by $\bm{*}$-automorphisms} of a simply connected Lie group $G$, having $\mathfrak{g}$ as its Lie algebra, such that $d\alpha(X) = \overline{\eta(X)|_{\mathcal{A}_\infty}}$ is the generator of a $\Gamma_2$-isometrically equicontinuous one-parameter group, for all $X \in \mathfrak{g}$.}\\

\textbf{Proof of Theorem 3.6:} The only thing needed to prove is that the $*$-derivations $\eta(X)$, $X \in \mathfrak{g}$, are all $\Gamma_2$-conservative, because then, Theorem 2.12 becomes available. In order to obtain this result, reference \cite[Theorem 1]{brattelidissipative} will be used.\\

As in the beginning of the proof of Theorem 2.12, the estimates $$\rho_{p, 1}(a) \leq \epsilon \, p(\Delta(a)) + \frac{E_p}{\epsilon} \, p(a)$$ are valid for every $p \in \Gamma_2$, $0 < \epsilon \leq 1$ and $a \in \mathcal{D}$, which shows that each of the $*$-derivations $\left. \delta_k \right|_{\mathcal{D}}$ possesses the (KIP) with respect to $\Gamma_2$. By the triangle inequality, all $*$-derivations $\delta_X|_\mathcal{D}$, $X \in \mathfrak{g}$, have the (KIP) with respect to $\Gamma_2$. Hence, one can define the induced $*$-derivations $$\delta_{x, p} \colon [a]_p \longmapsto [\delta_X(a)]_p, \qquad [a]_p \in \pi_p[\mathcal{D}],$$ for all $p \in \Gamma_2$. The induced operator $(\Delta|_\mathcal{D})_p = \sum_{k = 1}^d \delta_{k, p}^2$ is dissipative on $\cap_{k = 1}^d \text{Dom }\delta_{k, p}^2 = \pi_p[\mathcal{D}]$, so \cite[Theorem 1]{brattelidissipative} mentioned above says that each induced $*$-derivation $$\delta_{k, p} \colon [a]_p \longmapsto [\delta_k(a)]_p, \qquad 1 \leq k \leq d, \, [a]_p \in \pi_p[\mathcal{D}],$$ is conservative, for every $p \in \Gamma_2$. Now, fix $X := \sum_{k = 1}^d c_k \, X_k \in \mathfrak{g}$, $p \in \Gamma$ and $[a]_p \in \pi_p[\mathcal{D}]$. If $f \in (\mathcal{A} / N_p)'$ is a tangent functional at $[a]_p$ - in other words, if $f$ is a continuous linear functional on $\mathcal{A}/N_p$ such that $|f([a]_p)| = \|f\| \, \|[a]_p\|_p$ - then by conservativity of the basis elements and \cite[Proposition 3.1.14, page 175]{bratteli1}, it follows that $$\text{Re }f(\delta_{x, p}([a]_p)) = \text{Re }f([\delta_X(a)]_p) = \sum_{k = 1}^d c_k \, \text{Re }f(\delta_{k, p}([a]_p)) = 0,$$ so $\delta_{x, p}$ is conservative. Since $X$ and $p$ are arbitrary, this proves $\delta_X|_\mathcal{D}$ is $\Gamma_2$-conservative, for every $X \in \mathfrak{g}$. By Theorem 2.12 and Lemma 3.2, this implies $(\mathcal{D}, \mathfrak{g}, \eta)$ exponentiates to a representation by $*$-automorphisms which are $\Gamma_2$-isometrically equicontinuous one-parameter groups. \hfill $\blacksquare$\\

In view of Theorem 3.6, a small upgrade of Theorem 3.4 can be done for locally C$^*$-algebras. More precisely, the hypothesis that the $*$-derivations $\left\{\delta_k\right\}_{1 \leq k \leq d}$ must be $\Gamma$-conservative is not needed, and a theorem in perfect analogy with the C$^*$ one \cite[Corollary 2]{brattelidissipative} can be obtained:\\

\textbf{Theorem 3.7 (Exponentiation - Locally C$\bm{^*}$-Algebras, II):} \textit{Let $\mathcal{A}$ be a locally C$^*$-algebra, $\mathfrak{g}$ a real finite-dimensional Lie algebra with an ordered basis $\mathcal{B} := \left(X_k\right)_{1 \leq k \leq d}$ and $(\mathcal{A}_\infty, \mathfrak{g}, \eta)$ a representation of $\mathfrak{g}$ by closed $*$-derivations, with $$\mathcal{\mathcal{A}}_\infty := \bigcap_{n = 1}^{+ \infty} \bigcap \left\{\text{Dom }[\delta_{i_1} \ldots \delta_{i_k} \ldots \delta_{i_n}]: 1 \leq k \leq n, 1 \leq i_k \leq d \right\},$$ where $\delta_X := \eta(X)$, $\delta_k := \eta(X_k)$, $1 \leq k \leq d$, and $\Delta := \sum_{k = 1}^d \delta_k^2$. Then, the following are equivalent:}

\begin{enumerate}

\item \textit{$(\mathcal{A}_\infty, \mathfrak{g}, \eta)$ exponentiates to a strongly continuous locally equicontinuous representation \linebreak $\alpha \colon G \longrightarrow \mathcal{L}(\mathcal{A})$ \textbf{by $\bm{*}$-automorphisms} of a simply connected Lie group $G$, having $\mathfrak{g}$ as its Lie algebra, such that $d\alpha(X) = \overline{\eta(X)|_{\mathcal{A}_\infty}}$ is the generator of a $\Gamma$-isometrically equicontinuous one-parameter group, for all $X \in \mathfrak{g}$;}
\item \textit{$\Delta$ is an infinitesimal pregenerator of a $\Gamma$-contractively equicontinuous semigroup $t \longmapsto S(t)$ satisfying $S(t)[\mathcal{A}] \subseteq \mathcal{A}_\infty$, for all $t \in (0, 1]$ and, for each $p \in \Gamma$, there exists $C_p > 0$ for which the estimates $$\rho_{p, 1}(S(t)a) \leq C_p \, t^{-\frac{1}{2}} p(a)$$ are verified for all $t \in (0, 1]$ and $a \in \mathcal{A}$.}

\end{enumerate}

\section*{Acknowledgements}
\addcontentsline{toc}{section}{Acknowledgements}

\indent

First of all, the author would like to express his most sincere thanks to his parents, Ione and Juarez, for their always unconditional support.

This work is part of a PhD thesis, so the author would like to thank his advisor, professor Michael Forger, for his insightful suggestions and for always sharing his scientific knowledge and points of view about Science in such a generous, sincere and inspiring way. His careful revision of this paper was also very much appreciated.

The author would also like to thank professor Severino T. Melo for many helpful discussions in the early stages of this work, and to professor Paulo D. Cordaro, for the clarifications on the topic of strongly elliptic operators.

\end{document}